\documentclass[11pt]{article}

\usepackage{latexsym,amsfonts,amsmath,amssymb,graphicx,hyperref,verbatim}

\newcommand{\Om}{\Omega}
\newcommand{\T}{{\cal T}}

\newcommand{\SA}{\mbox{\bf S}}

\newcommand{\St}{\mathbb{S}}
\newcommand{\EE}{\mbox{\bf E}\,}
\newcommand{\PP}{\mbox{\bf P}\,}

\newcommand{\R}{\mathbb{R}}
\newcommand{\C}{\mathbb{C}}
\newcommand{\Q}{\mathbb{Q}}
\newcommand{\HH}{\mathbb{H}}
\newcommand{\N}{\mathbb{N}}
\newcommand{\D}{\mathbb{D}}
\newcommand{\TT}{\mathbb{T}}
\newcommand{\Z}{\mathbb{Z}}
\newcommand{\A}{\mathbb{A}}

\newcommand{\B}{\mbox{\bf B}}

\newcommand{\pa}{\partial}

\newcommand{\F}{{\cal F}}
\newcommand{\no}{\noindent}

\newcommand{\mr}{\mathring}
\def\BGE{\begin{equation}}
\newcommand{\BGEN}{\begin{equation*}}
\def\EDE{\end{equation}}
\newcommand{\EDEN}{\end{equation*}}

\def\eps{\varepsilon}
\def\til{\widetilde}
\def\ha{\widehat}
\def\sem{\setminus}
\def\lin{\overline}
\def\vphi{\varphi}

\def\del{\delta}
\def\Del{\Delta}
\DeclareMathOperator{\ccap}{cap} \DeclareMathOperator{\rad}{rad}
\DeclareMathOperator{\sign}{sign} \DeclareMathOperator{\diam}{diam}
\DeclareMathOperator{\dist}{dist} \DeclareMathOperator{\dcap}{dcap}
\DeclareMathOperator{\hcap}{hcap} \DeclareMathOperator{\id}{id}
\DeclareMathOperator{\Imm}{Im } \DeclareMathOperator{\Ree}{Re }
\DeclareMathOperator{\Tip}{Tip} \DeclareMathOperator{\Set}{Set}

\def\dto{\stackrel{\rm Cara}{\longrightarrow}}
\def\h0{{\bf h}}
\def\luto{\stackrel{\rm l.u.}{\longrightarrow}}
\def\hto{\stackrel{\cal H}{\longrightarrow}}
\def\aequ{\stackrel{a}{\sim}}
\def\LERW{\mbox{LERW}}

\newtheorem{Lemma}{Lemma}[section]
\newtheorem{Theorem}{Theorem}[section]
\newtheorem{Definition}{Definition}[section]
\newtheorem{Corollary}{Corollary}[section]
\newtheorem{Proposition}{Proposition}[section]
\numberwithin{equation}{section}

\begin{document}
\title{\bf Continuous LERW Started\\ from Interior Points}
\date{\today}
\author{Dapeng Zhan\footnote{Supported by NSF grant 0963733} \\  Michigan State University}
\maketitle


\begin{abstract} We use the whole-plane Loewner equation to define a family of
continuous LERW in finitely connected domains that are started from interior points.
These continuous LERW satisfy conformal invariance,  preserve some continuous
local martingales, and are the scaling limits of the corresponding discrete LERW
on the discrete approximation of the domains.
\end{abstract}

\section{Introduction}
This paper is a follow-up of \cite{LERW}, in which we defined a family of
random curves called continuous LERW in finitely connected plane domains,
 and proved that they are the
scaling limits of the corresponding discrete LERW (loop-erased conditional random walk).

The continuous LERW defined in \cite{LERW} is a simple curve that grows from a
boundary point (or prime end, c.f.\ \cite{Ahl}), say $a$, of some domain, say $D$, and aims at a certain target, which could be
an interior point, a boundary arc  or another boundary point of $D$.
It is an SLE$_2$-type process that satisfies conformal invariance, which behaves locally like the SLE$_2$ process in simply
connected domains introduced by Oded Schramm (\cite{S-SLE}).
The special cases are when $D$ is a subdomain of the upper half plane $\HH=\{z\in\C:\Imm z>0\}$,
$a=0$, and the part of $\pa D$ near $a$ lies on $\R$. In this case, the LERW is the chordal Loewner evolution
driven by some semi-martingale, whose martingale part is $\sqrt 2$ times a Brownian motion, and
whose differentiable part contains the information of the domain and the target set. The continuous
LERW is first defined in the special cases, and then extended to general cases via conformal maps.

The corresponding discrete LERW is defined on the graph $D^\del$, which is the grid
approximation of $D$ by $\del\Z^2$ for some small $\del>0$. For the construction, we
first start a simple random walk on $D^\del$ from an interior vertex that is closest
to $a$, and stop it when it leaves the domain or hits a vertex that is closest to
the target. Then we condition this stopped random walk on the event that it ends
at a vertex that is closest to the target. Finally, we erase the loops on this conditional
random walk in the order they appear, and get the discrete LERW.

The convergence of the discrete LERW curves to the corresponding continuous LERW curves were proved
using the technique introduced by \cite{LSW-2}:
 first use Skorokhod's embedding theorem to prove the convergence of the driving function, and then use the tameness of
the discrete LERW curve to prove the convergence of the curves.

\vskip 3mm

This paper will consider the case when the start point $a$ is not a boundary point, but an interior point of $D$. It
is natural to define a discrete  LERW that starts from a vertex of $D^\del$ which is closest to $a$ and aims at a given target.
The motivation of this paper is to describe the scaling limit of this lattice path.
We will uses whole-plane Loewner equation (\cite{LawSLE}) to define a family of random curves,
which are still called continuous LERW, and prove that they are the scaling limits of the above discrete LERW.

For the definition of continuous LERW in the domain $D$ started from the interior point $a=0$ and aimed at
another interior point, say $z_e$, we solve an integral equation as below. For $\xi\in C((-\infty,T))$.
Let $K^\xi_t$ and $\vphi^\xi_t$, $-\infty<t<T$, be the whole-plane Loewner hulls and maps, respectively, driven by $\xi$
(c.f.\ Section 4.3 of \cite{LawSLE} or Section \ref{Section-Whole-Plane-Loewner-Equation} of this paper).
Suppose $K^\xi_t\subset D\sem\{z_e\}$ for $-\infty<t<T$.
Then for each $t\in(-\infty,T)$, $D\sem K^\xi_t$ is a finitely connected domain
containing $z_e$. Let
\BGE X^\xi(t)=(\pa_x\pa_y/\pa_y)[G(D\sem K^\xi_t,z_e;\cdot)\circ (\vphi^\xi_t)^{-1}\circ e^i\circ R_\R](\xi(t)),\label{X-1st-def}\EDE
where $G(D\sem K^\xi_t,z_e;\cdot)$ is Green's function in $D\sem K^\xi_t$ with pole at $z_e$,
 $e^i$ is the map $z\mapsto e^{iz}$, and $R_\R$ is the conjugate map $z\mapsto \lin{z}$.
Let $\kappa=2$, and
  $B^{(\kappa)}_\R(t)$, $-\infty<t<\infty$, be a driving function for whole-plane SLE$_\kappa$ (c.f.\ Section 6.6 of \cite{LawSLE}
  or Section \ref{Def-of-quasi} of this paper).
Let $\lambda=2$, and $\xi(t)$, $-\infty<t<T$, be the solution to the integral equation
\BGE \xi(t)=B^{(\kappa)}_\R(t) +\lambda\int_{-\infty}^t X^\xi(t) dt,\label{equation-LERW}\EDE
such that $(-\infty,T)$ is the maximal interval of the solution. It turns out that the solution exists,
and is a semi-martingale. So there is a random continuous curve $\beta(t)$, $-\infty\le t<T$, such that
$\beta(-\infty)=0$ and $K^\xi_t=\beta([-\infty,t])$, $-\infty<t<T$. Such $\beta$ is called the continuous LERW curve
in $D$ from $0$ to $z_e$. If the target is a boundary arc or another boundary point, we will use harmonic measure function
or Poisson kernel function instead of Green function in (\ref{X-1st-def}), and keep other formulas in the definition unchanged.

We then prove that these continuous LERW satisfy conformal invariance, and preserve
some continuous local martingales generated by generalized Poisson kernels.
Finally, we use the technique in \cite{LSW-2} and \cite{LERW}  to show
that these continuous LERW are the scaling limits of the corresponding discrete LERW.

\vskip 3mm
The continuous LERW defined in this paper turns out to be locally absolutely continuous w.r.t.\ the whole-plane SLE$_2$. In fact, if $U$ is a simply connected subdomain of $D$ that contains the initial point $0$, and is bounded away from $\pa D$ and the target, then the continuous LERW stopped at the time $\tau_U$ when it exits $U$ has a distribution absolutely continuous w.r.t.\ the whole-plane SLE$_2$ stopped at $\tau_U$. Moreover, there is a local martingale process $M(t)$ such that the above Radon-Nikodym derivative is $M(\tau_U)$. The formula of $M(t)$ will be given in Section \ref{distr}. So this gives an alternative way to define continuous LERW. First one may use whole-plane SLE$_2$ and the Radon-Nikodym derivative $M(\tau_U)$ to define a partial continuous LERW (stopped at $\tau_U$), say $\gamma_U$, for every $U$. Using the local martingale property of $M(t)$, one can check that these partial processes are consistent w.r.t.\ each other: $\gamma_{U_1}$ stopped at $\tau_{U_2}$ has the same distribution as $\gamma_{U_2}$ stopped at $\tau_{U_1}$. Then one may construct a complete continuous LERW $\gamma$ such that $\gamma$ stopped at any $\tau_U$ has the distribution of $\gamma_U$.

We prefer the definition using the driving function rather than Radon-Nikodym derivative. This is because when we prove the convergence of discrete LERW, the technique in \cite{LSW-2} and the Skorokhod's embedding theorem can be easily applied here without major modifications. If one uses the other definition, and tries to prove the convergence, he first has to work out the convergence of a particular  discrete LERW to the whole-plane SLE$_2$, and then show the convergence of the discrete Radon-Nikodym derivative (between discrete LERW) to the continuous Radon-Nikodym derivative $M(\tau_U)$. The first step requires no much less work than the other approach, while the second step seems very difficult to the author.

The Radon-Nikodym derivative approach is useful in other respects. For example, one may use the density functions together with the stochastic coupling technique introduced in \cite{reversibility} to prove the reversibility of continuous LERW without using discrete LERW. One may also use them to show that the continuous LERW is a loop-erasure of a plane Brownian motion restricted in the domain (\cite{LEBM}).

Unlike the SLE processes started from boundary, there are SLE$_\kappa$-type processes started from $0$, which are not locally absolutely continuous w.r.t.\ whole-plane SLE$_\kappa$ process. One example is the whole-plane Loewner process driven by $\xi(t)=B^{(\kappa)}_\R(t)+\sigma t$, where $\sigma$ is a nonzero real constant. Although this is not the case for continuous LERW, some care is required when dealing with the definition of SLE started from interior points.

\vskip 3mm

We expect that the definition of the continuous LERW started from interior points will shed some
light on the definition of some other random curves started from interior points, e.g., the
reversal of radial SLE curves, and the scaling limits of self-avoiding walks (SAW) that connect
two interior points. In particular, our result implies a description of the reversal of radial SLE$_2$.

This paper is organized in the following way. In Section \ref{prelim}, we review   some basic notation including the radial
Loewner equations and whole-plane Loewner equations. We also study the Carath\'eodory topology restricted to the
space of interior hulls.  In Section \ref{Continuous LERW},  we give the detailed definition of continuous
LERW started from interior points, and prove that such LERW satisfies conformal invariance, and preserves a family of
continuous local martingales generated by the generalized Poisson kernels. In Section \ref{Existence}, we prove that the
solution to (\ref{equation-LERW}) exists uniquely, and is a semi-martingale. In Section \ref{distr}, we prove that the
continuous LERW started from an interior point is locally absolutely continuous w.r.t.\ the whole-plane SLE$_2$ process. In the last
section, we introduce a family of discrete LERW defined on the discrete approximation of the domain, and
a sketch of a proof is given to show that the scaling limit of this discrete LERW is the continuous LERW
defined in this paper.

We will frequently cite notation and theorems from \cite{LERW}. The readers are suggested have a copy
of \cite{LERW} at hand. We will often use some basic properties of the SLE processes. The reader may refer
\cite{RS-basic} and \cite{LawSLE} for the background of SLE.

\section{Preliminary} \label{prelim}
\subsection{Some notation}
We adopt the notation in Section 2 of \cite{LERW} about finitely connected domain,
conformal closure, prime end, side arc, Green function, generalized Poisson kernel,
harmonic measure function, hull and Loewner chain, and etc. But now we call
the hull and Loewner chain in \cite{LERW} the boundary hull and
boundary Loewner chain, respectively, to distinguish them from the
interior hull and interior Loewner chain that will be defined in
this paper.

Throughout this paper, we use the following notation. Let $\ha\C
=\C\cup\{\infty\}$ denote Riemann sphere. Let $\HH$ be the
upper half plane $\{z\in\C:\Imm z>0\}$. Let $\D$ be the unit disc
$\{z\in\C:|z|<1\}$. Let $\TT$ be the
unit circle $\{z\in\C:|z|=1\}$.
Let $\St_h$ be the strip
$\{z\in\C:h>\Imm z>0\}$ for $h>0$. Let $\R_h$ be the line
$\{z\in\C:\Imm z=h\}$ for $h\in\R$. Then $\St_h$ is bounded by $\R$
and $\R_h$. Let $\A_h$ be the annulus $\{z\in\C:e^{-h}<|z|<1\}$ for
$h>0$. We define an almost-$\D$ domain to be a finitely
connected subdomain of $\D$ which contains $0$ and $\A_h$ for some $h>0$.

Let $e^i$ be the
map $z\mapsto e^{iz}$. Then $e^i$ is the covering map from $\HH$
onto $\D\sem\{0\}$, from $\St_h$ onto $\A_h$, and from $\R$ onto
$\TT$. Let $R_\R(z)=\lin z$ be the complex conjugate map. Let
$R_\TT(z)=1/\lin z$ be the reflection about $\TT$. Then $e^i\circ
R_\R=R_\TT\circ e^i$. For $w\in\C$, let $A_w$ denote the map
$z\mapsto w+z$; let $M_w$ denote the map $z\mapsto wz$. Then
$e^i\circ A_w=M_{e^i(w)}\circ e^i$.
Let $\B(z_0;r)$ be the ball $\{z\in\C:|z-z_0|<r\}$.
If $\sigma$ is a Jordan curve in $\C$, we use $U(\sigma)$ to
denote the bounded connected component of $\C\sem \sigma$, and let
$H(\sigma):=\lin{U(\sigma)}=U(\sigma)\cup\sigma$.
If $I$ is an interval on $\R$, let $C(I)$ denote the set of real
valued continuous functions on $I$. For $f\in C(I)$, if $[a,b]\subset I$, let
$\|f\|_{a,b}=\max\{|f(x)|:x\in[a,b]\}$; if
$(-\infty,a]\subset I$, let $\|f\|_a=\sup\{|f(x)|:x\le a\}$.

\subsection{Radial Loewner equation} \label{section-radial-Loewner}
If $H$ is a boundary hull in $\D$ such that $0\not\in H$, then we
say that $H$ is a boundary hull in $\D$ w.r.t.\ $0$. For such $H$,
there is a unique map $\psi_H$ that maps $\D\sem H$ conformally onto
$\D$ such that $\psi_H(0)=0$ and $\psi_H'(0)>0$. Then
$\dcap(H):=\ln(\psi_H'(0))\ge 0$ is called the capacity of $H$ in
$\D$ w.r.t.\ $0$. For example, $\emptyset$ is a boundary hull in $\D$
w.r.t.\ $0$, $\psi_\emptyset=\id_\D$, and $\dcap(\emptyset)=0$. From
Schwarz lemma, $|\psi_H(z)|\ge |z|$ for any $z\in\D\sem H$. If
$H_1\subset H_2$ are boundary hulls in $\D$ w.r.t.\ $0$, define
$H_2/H_1=\psi_{H_1}(H_2\sem H_1)$. Then $H_2/H_1$ is also a boundary
hull in $\D$ w.r.t.\ $0$, and we have $\psi_{H_2/H_1}=\psi_{H_2}\circ
\psi_{H_1}^{-1}$ and $\dcap(H_1)+\dcap(H_2/H_1)=\dcap(H_2)$. Thus,
 $|\psi_{H_2}(z)|\ge |\psi_{H_1}(z)|$ for any $z\in \D\sem H_2$.

The following proposition is the radial version of Lemma 2.8 in
\cite{LSW1}. The proof is similar. So we omit the proof.

\begin{Proposition} Let $\Xi$ be an open neighborhood of $x_0\in\TT$
in $\D$. Suppose $W$ maps $\Xi$ conformally into $\D$ such that,
 as $z\to\TT$ in $\Xi$,
$W(z)\to\TT$. Such $W$ extends conformally across $\TT$ near $x_0$
by Schwarz reflection principle. Then we have
\BGE \lim_{H\to x_0} \frac{\dcap(W(H))}{\dcap(H)}=|W'(x_0)|^2,\label{capacity}\EDE
where $H\to x_0$ means that $H$ is a nonempty hull in $\D$ w.r.t\ $0$, and $\diam(H\cup\{x_0\})\to 0$.
\label{dcap}
\end{Proposition}

Suppose $\xi\in C([0,T))$ for some $T\in(0,+\infty]$. The radial
Loewner equation driven by $\xi$ is as follows: \BGE
\pa_t\psi_t(z)=\psi_t(z)\,\frac{e^{i\xi(t)}+\psi_t(z)}
{e^{i\xi(t)}-\psi_t(z)},\qquad\psi_0(z)=z.\label{radial-Loewner}\EDE
For $0\le t<T$, let $L_t$ be the set of $z\in\D$ such that the
solution $\psi_s(z)$ blows up before or at time $t$. Then $L_t$ is a
boundary hull in $\D$ w.r.t.\ $0$, and $\psi_t=\psi_{L_t}$ for each
$t\in[0,T)$. We call $L_t$ and $\psi_t$, $0\le t<T$, the radial
Loewner hulls and maps, respectively, driven by $\xi$. We have the
following proposition.

\begin{Proposition} (a) Suppose $L_t$ and $\psi_t$, $0\le t<T$,
are the radial Loewner hulls and maps, respectively, driven by $\xi\in
C([0,T))$. Then $(L_t,0\le t<T)$ is a boundary Loewner chain in $\D$
avoiding $0$, and $\dcap(L_t)=t$ for any $0\le t<T$. Moreover, \BGE
\{e^{i\xi(t)}\}=\bigcap_{\eps\in(0,T-t)}\lin{L_{t+\eps}/
L_t},\quad 0\le t<T.\label{radial-driving}\EDE (b) Suppose $L_t$, $0\le t<T$, is a
boundary Loewner chain in $\D$ avoiding $0$, and $\dcap(L_t)=t$ for
any $0\le t<T$. Then there is $\xi\in C([0,T))$ such that $L_t$,
$0\le t<T$, are radial Loewner hulls driven by $\xi$.
\label{radial-Loewner-chain}
\end{Proposition}
{\bf Proof.} This is the main result in \cite{Pom-LDE}. $\Box$

\vskip 3mm

The covering radial Loewner equation driven by $\xi$ is: \BGE
\pa_t\til\psi_t(z)=\cot_2({\til\psi_t(z)-\xi(t)}),\quad\til\psi_0(z)=z.
\label{radial-Loewner-covering}\EDE In this paper, we use $\cot_2(z)$ to
denote the function $\cot(z/2)$. For $0\le t<T$, let $\til L_t$
be the set of $z\in\HH$ such that the solution $\til\psi_s(z)$ blows
up before or at time $t$. We call $\til L_t$ and $\til \psi_t$,
$0\le t<T$, the covering radial Loewner hulls and maps,
respectively,  driven by $\xi$. Then $\til\psi_t$ maps $\HH\sem\til
L_t$ conformally onto $\HH$, and satisfies
$\til\psi_t(z+2k\pi)=\til\psi_t(z)+2k\pi$ for any $k\in\Z$. Since
$\Imm\cot_2(z)<0$ for $z\in\HH$, so $\Imm\til\psi_t(z)$ decreases in
$t$. Suppose $L_t$ and $\psi_t$, $0\le t<T$, are the radial Loewner
hulls and maps, respectively, driven by $\xi$, then for any
$t\in[0,T)$, $\til L_t=(e^i)^{-1}(L_t)$, and $\psi_t\circ
e^i=e^i\circ\til\psi_t$.



\subsection{Interior Hulls and Interior Loewner
Chains}\label{sec-hull}
Suppose $D$ is a finitely connected domain. If $\emptyset\ne
F\subset D$ is compact and connected, and $D\sem F$ is also
connected, then we say that $F$ is an interior hull in $D$. If $F$
contains only one point, we say $F$ is degenerate; otherwise, $F$ is
non-degenerate. If $F$ is a non-degenerate interior hull in an
$n$-connected domain $D$, then $D\sem F$ is an $(n+1)$-connected
domain. If $H$ is another interior hull in $D$, and $F\subset H$,
then $H\sem F$ is a boundary hull in $D\sem F$.

Let $T\in(-\infty,+\infty]$.
We say the family $F(t)$, $-\infty<t<T$,  is an interior
Loewner chain in $D$ started from $z_0\in D$ if (i) for each $t\in
(-\infty,T)$, $F(t)$ is a non-degenerate interior hull in $D$; (ii)
$F(t_1)\subsetneqq F(t_2)$ for any $t_1<t_2<T$; (iii) for any
$t_0\in(-\infty,T)$, $(F(t_0+t)\sem F(t_0),0\le t<T-t_0)$ is a boundary
Loewner chain in $D\sem F(t_0)$; and (iv) $\bigcap_{-\infty <t<T} F(t)=\{z_0\}$. For
any $t_0<T$, if $(F(t_0+t)\sem F(t_0),0\le t<T-t_0)$ is started from
a prime end $w(t_0)$ of $D\sem F(t_0)$, then we say that $w(t_0)$ is
the prime end determined by $(F(t))$ at time $t_0$. Suppose $u$ is a
continuous (strictly) increasing function on $(-\infty,T)$, and
satisfies $u(-\infty)=-\infty$, that is, $\lim_{t\to
-\infty}u(t)=-\infty$. Let $u(T):=\lim_{t\to T}u(t)$. Then
$F(u^{-1}(t))$, $-\infty<t<u(T)$, is also an interior Loewner chain in
$D$ started from $z_0$. We call it the time-change of $(F(t))$ through
$u$. Suppose $\gamma :[-\infty,T)\to D$ is a simple curve. For
$t\in(-\infty,T)$, let $F(t)=\gamma([-\infty,t])$. Then $(F(t))$ is
an interior Loewner chain started from $\gamma(-\infty)$. We call
such $F$ the interior Loewner chain generated by $\gamma$. Then for
each $t<T$, $\gamma(t)$ is the prime end determined by $(F(t))$ at
time $t$.

If $F$ is an interior hull in $\ha\C$, and $\infty\not\in F$, then
we call $F$ a bounded interior hull. For example, if $\sigma$ is a
Jordan curve in $\C$, then $H(\sigma)$ is a bounded interior hull.
For any bounded interior hull $F$, there is a unique function
$\phi_F$ that maps $\ha\C\sem F$ conformally onto
$\ha\C\sem\lin{r\D}$ for some $r\ge 0$ such that
$\phi_F(\infty)=\infty$ and
$\phi_F'(\infty):=\lim_{z\to\infty}z/\phi_F(z)=1$. We call
$\rad(F):=r$ the radius of $F$, and $\ccap(F):=\ln(r)$ the capacity
of $F$ w.r.t.\ $\infty$. Here if $F$ contains only one point, say
$z_0$, then $\phi_F(z)=z-{z_0}$, so $\rad(F)=0$ and
$\ccap(F)=\ln(0)=-\infty$. If $F$ is non-degenerate, then
$\rad(F)>0$ and $\ccap(F)\in\R$, and we define
$\vphi_F:=M^{-1}_{\rad(F)}\circ \phi_F$. Then $\vphi_F$ maps
$\ha\C\sem F$ conformally onto $\ha\C\sem\lin\D$, and satisfies
$\vphi_F(\infty)=\infty$ and $\vphi_F'(\infty)>0$. Let
$\psi_F=R_\TT\circ \vphi_F\circ R_\TT$. Then $\psi_F$
maps $\ha\C\sem R_\TT(F)$ conformally onto $\D$, and satisfies
$\psi_F(0)=0$ and $\psi_F'(0)>0$.

The following results are well known (e.g., c.f.\ \cite{LawSLE}). If
$F$ is a bounded interior hull, $a,b\in\C$, then $\rad(aF+b)=|a|\rad
(F)$; $\rad(\lin{\B(z_0;r)})=r$ for any $z_0\in\C$ and $r>0$;
 $\rad(F)\ge \diam(F)/4$ for any bounded interior hull $F$, and the
equality holds if and only if $F$ is a line segment or a single
point. By taking logarithm, we get the corresponding results for
$\ccap(F)$. Suppose $F_1\subset F_2$ are two non-degenerate bounded
interior hulls. Then $\ccap(F_1)\le\ccap(F_2)$, where the equality
holds only if $F_1=F_2$. Let $F_2/F_1:=R_\TT\circ\vphi_{F_1}(F_2\sem
F_1)$. Then $F_2/F_1$ is a boundary hull in $\D$ w.r.t.\ $0$. Moreover, we have
\BGE\psi_{F_2/F_1}=R_\TT\circ\vphi_{F_2}\circ\vphi_{F_1}^{-1}\circ
R_\TT,\label{psi-vphi}\EDE and
$\dcap(F_2/F_1)=\ccap(F_2)-\ccap(F_1)$. Since
$|\psi_{F_2/F_1}(z)|\ge|z|$ for any $z\in\D\sem(F_2/F_1)$, so
$|\vphi_{F_1}(z)|\ge|\vphi_{F_2}(z)|$ for any $z\in \C\sem F_2$. If
$F_1\subset F_2\subset F_3$ are non-degenerate bounded interior
hulls, then $F_2/F_1\subset F_3/F_1$, and
$(F_3/F_1)/(F_2/F_1)=F_3/F_2$. Here $F_3/F_1$ and $F_2/F_1$ are
boundary hulls in $\D$ w.r.t.\ $0$, and the quotient between $F_3/F_1$ and $F_2/F_1$  uses the definition in the last subsection

Let ${\cal H}$ denote the set of all bounded interior hulls, and let
${\cal H}_0$ denote the set of $H\in{\cal H}$ such that $0\in
H$. From Proposition 3.30 in \cite{LawSLE}, there is an absolute
constant $C_{\cal H}\ge 3$ such that, for any $H\in{\cal H}_0$ with
$\rad(H)=1$, and any $z\in\C$ with $|z|>1$,
$$|\phi_H^{-1}(z)-z|\le C_{\cal H}.$$
Suppose $H\in{\cal H}_0$ is non-degenerate. Let
$H_0=H/\rad(H)\in{\cal H}_0$. Then $\rad(H_0)=1$. So for any
$z\in\C$ with $|z|>\rad(H)$,
$\phi_H^{-1}(z)=\rad(H)\phi_{H_0}^{-1}(z/\rad(H))$, which implies
$$|\phi_H^{-1}(z)-z|=|\phi_{H_0}^{-1}(z/\rad(H))-z/\rad(H)|\rad(H)\le
C_{\cal H}\rad(H).$$ If $H\in{\cal H}$ is non-degenerate, then there
is $z_0\in H$ with $|z_0|=\dist(0,H)$. Then $H_0=H-z_0\in{\cal
H}_0$, $\rad(H_0)=\rad(H)$, and
$\phi_H^{-1}=A_{z_0}\circ\phi_{H_0}^{-1}$. Thus, for any
$|z|>\rad(H)$, \BGE |\phi_H^{-1}(z)-z|\le
|z_0|+|\phi_{H_0}^{-1}(z)-z|\le \dist(0,H)+C_{\cal
H}\rad(H).\label{phi_H-1(z)-z}\EDE If $H=\{z_0\}$ is degenerate, (\ref{phi_H-1(z)-z})
still holds because $\phi_H^{-1}=A_{z_0}$ and $\dist(0,H)=|z_0|$.
For any interior hull $H$, Since $\phi_H^{-1}$ maps
$\{|z|>\rad(H)\}$ onto $\C\sem H$, so for any $z\in\C\sem H$,
\BGE|\phi_H(z)-z|\le
 \dist(0,H)+C_{\cal H}\rad(H).\label{phi_H(z)-z}\EDE

\subsection{Whole-plane Loewner equation}\label{Section-Whole-Plane-Loewner-Equation}
Suppose $F(t)$, $-\infty<t<T$, is an interior Loewner chain (in
$\ha\C$) avoiding $\infty$, that is, $\infty\not\in F(t)$ for any
$t<T$.  If $\ccap(F(t))=t$ for any $t<T$, we say $(F(t))$ is
parameterized by capacity. In the general case, $v(t):=\ccap(F(t))$ is a
continuous increasing function such that $v(-\infty)=-\infty$, and
the time-change of $(F(t))$ through $v$ is parameterized by
capacity. The following proposition is a combination of a theorem in \cite{LawSLE}
and its inverse statement.

\begin{Proposition}   (i) Suppose $\xi\in C((-\infty,T))$.
Then there is a unique interior Loewner chain $K_t$, $-\infty<t<T$, which is
started from $0$, avoids $\infty$, and is parameterized by capacity,
such that the followings hold. For $-\infty<t<T$, let $\vphi_t=\vphi_{K_t}$. Then
$\vphi_t$ satisfies \BGE
\pa_t\vphi_t(z)=\vphi_t(z)\,\frac{e^{i\xi(t)}+\vphi_t(z)}{e^{i\xi(t)}-\vphi_t(z)};
\label{whole-plane-equation}\EDE and for any $z_0\in\C\sem\{0\}$,
\BGE \lim_{t\to-\infty}e^t\vphi_t(z_0)=z_0.\label{near-infty}\EDE
(ii) Suppose $K_t$, $-\infty<t<T$, is an interior
Loewner chain started from $0$ avoiding $\infty$, and is
parameterized by capacity. Then there is $\xi\in C((-\infty,T))$ such that for
$\vphi_t=\vphi_{K_t}$, (\ref{whole-plane-equation}) and (\ref{near-infty}) both hold.
 \label{Whole-plne-Proposition-I}
\end{Proposition}
{\bf Proof.} (i) This is a special case of Proposition 4.21 in \cite{LawSLE}, where
$\mu_t=\delta_{e^i(\xi(t))}$.

(ii) Fix $t\in(-\infty,T)$. Since $K_t\in{\cal H}_0$,
$\rad(K_t)=e^t$, and $\phi_{K_t}=M_{e^t}\circ\vphi_t$, so from (\ref{phi_H(z)-z}),
\BGE |e^t\vphi_t(z)-z|<C_{\cal H}e^t,\quad z\in\C\sem K_t.\label{est-2}\EDE
Fix
$z_0\in\C\sem\{0\}$. Since $\{0\}=\bigcap_{-\infty<t<T} K_t$, so there is
$T_0\in(-\infty,T)$ such that $z_0\not\in K_t$ for $t\le T_0$,
which implies that $|e^t\vphi_t(z_0)-z_0|<C_{\cal H}e^t$ from (\ref{est-2}). Thus,
(\ref{near-infty}) holds.

Fix $b\in(-\infty,T)$. Since $(K_{b+t}/K_b,0\le t<T-b)$ is a
boundary Loewner chain in $\D$ avoiding $0$, and
$\dcap(K_{b+t}/K_b)=\ccap(K_{b+t})-\ccap(K_b)=t$ for $0\le t<T-b$,
so from Proposition \ref{radial-Loewner-chain}, there is $\xi_b\in
C([0,T-b))$ such that $K_{b+t}/K_b$, $0\le t<T-b$, are the radial
Loewner hulls driven by $\xi_b$, and for each $t\in[0,T-b)$, \BGE
\{e^{i\xi_b(t)}\}=\bigcap_{\eps\in(0,T-b-t)}\lin{(K_{b+t+\eps}/K_b)/
(K_{b+t}/K_b)}=\bigcap_{\eps\in(0,T-b-t)}\lin{K_{b+t+\eps}/
K_{b+t}}.\label{xi^b(t)}\EDE For $t\in(-\infty,T)$, choose
$b\in(-\infty,t]$, and let $\chi(t)=e^i(\xi_b(t-b))$. From
(\ref{xi^b(t)}), the value of $\chi(t)$ does not depend on the
choice of $b$. Since $\xi_b\in C([0,T-b))$ for each $b<T$, so $\chi$
is a $\TT$-valued continuous function. Thus, there is $\xi\in
C((-\infty,T))$ such that $\chi(t)=e^i(\xi(t))$ for $-\infty<t<T$.
Since
$e^i(\xi_b(t))=e^i(\xi(b+t))$ for $0\le t<T-b$, so $K_{b+t}/K_b$,
$0\le t<T-b$, are also the radial Loewner hulls driven by
$\xi(b+\cdot)$. Let $\psi^b_t$, $0\le t<T-b$, be the radial Loewner maps
driven by $\xi(b+\cdot)$. Then for $t\in[0,T-b)$,
$\psi^b_t=\psi_{K_{b+t}/K_b}=R_\TT\circ\vphi_{b+t}\circ\vphi_b^{-1}\circ
R_\TT$, and
$$\pa_t \psi^b_t(z)=\psi^b_t(z)\,\frac{e^{i\xi(b+t)}+\psi^b_t(z)}
{e^{i\xi(b+t)}-\circ\psi^b_t(z)}.$$ Since $R_\TT\circ
\psi^b_t=\vphi_{b+t}\circ\vphi_b^{-1}\circ R_\TT$, and
$\vphi_b^{-1}\circ R_\TT$ maps $\D\sem( K_{b+t}/K_b)$ onto $\C\sem
K_{b+t}$, so (\ref{whole-plane-equation}) holds for $t\in[b,T)$.
Since $b\in(-\infty,T)$ could be arbitrary,
so (\ref{whole-plane-equation}) holds for all $t\in(-\infty,T)$. $\Box$

\vskip 3mm

In the above proposition, $K_t$ and $\vphi_t$, $-\infty<t<T$,  are called the whole-plane Loewner
hulls and maps, respectively, driven by $\xi$. Since $e^i(\xi_b(t))=e^i(\xi(b+t))$ for $0\le t<T-b$, and $b\in(-\infty,T)$ is arbitrary, so
from (\ref{xi^b(t)}) we get a formula  similar to (\ref{radial-driving}), which is
\BGE \{e^{i\xi(t)}\}=\bigcap_{\eps\in(0,T-t)}
\lin{K_{t+\eps}/K_t}, \qquad -\infty<t<T.\label{xi(t)-whole-plane} \EDE

For $t\in(-\infty,T)$, let $L_t=R_\TT(K_t)$ and $\psi_t=R_\TT\circ\vphi_t\circ R_\TT$. Then
$\psi_t=\psi_{K_t}$, $\C\sem L_t$
is a simply connected domain that contains $0$, $\psi_t$ maps
$\C\sem L_t$ conformally onto $\D$, fixes $0$, and satisfies
\BGE
\pa_t\psi_t(z)=\psi_t(z)\,\frac{e^{i\xi(t)}+\psi_t(z)}
{e^{i\xi(t)}-\psi_t(z)}.\label{reflection-whole-plane-equation}\EDE
We call $L_t$ and $\psi_t$ the inverted whole-plane Loewner
hulls and maps, respectively, driven by $\xi$.

The covering whole-plane Loewner equation is defined as follows.
Let $\til K_t=(e^i)^{-1}(K_t)$, $-\infty<t<T$. Suppose $\til \vphi_t$, $-\infty<t<T$, satisfy that for each $t$, $\til \vphi_t$
maps $\C\sem \til K_I(t)$ conformally onto $-\HH$, $e^i\circ\til \vphi_t=\vphi_t\circ e^i$, and the following differential
equation holds:
\BGE \pa_t\til\vphi_t(z)=\cot_2(\til \vphi_t(z)-\xi(t)),\quad -\infty<t<T;\label{covering-whole-eqn}\EDE
\BGE \lim_{t\to -\infty} (\til \vphi_t(z)-it)=z,\quad z\in\C.\label{covering-whole-limit}\EDE
Then we call $\til K_t$ and $\til \vphi_t$ the covering whole-plane Loewner hulls and maps, respectively, driven by $\xi$.
Such family of $\til \vphi_t$ exists and is unique. In fact, for each $t\in(-\infty,T)$, we can find some $\til \vphi_t$
that maps $\C\sem \til K_I(t)$ conformally onto $-\HH$ such that $e^i\circ\til g_I(t,\cdot)=g_I(t,\cdot)\circ e^i$.
Such $\til \vphi_t$ is not unique. Since $  \vphi_t$ is differentiable in $t$, so one may choose $\til \vphi_t$ such
that it is also differentiable in $t$. From (\ref{whole-plane-equation}) we conclude that (\ref{covering-whole-eqn}) must hold.
From (\ref{near-infty}) we conclude that $\lim_{t\to -\infty} (\til \vphi_t(z)-it)=z+i2n\pi$ for some $n\in\Z$, and such $n$
is the same for every $z$. Now we replace $\til \vphi_t$ by $\til \vphi_t-i2n\pi$. Then (\ref{covering-whole-eqn})
and (\ref{covering-whole-limit}) still hold. So we have the existence of $\til g_I(t,\cdot)$. The uniqueness follows from
the same argument.

For $-\infty<t<T$, let $\til L_t=R_{\R}(\til K_t)$ and $\til \psi_t=R_{\R}\circ \til \vphi_t\circ R_{\R}$.
Then we have $\til L_t=(e^i)^{-1}(L_t)$, $e^i\circ \til \psi_t=\psi_t\circ e^i$,  $\til\psi_t$ maps $\C\sem \til L_t$ conformally
onto $\HH$,  and satisfies \BGE \pa_t\til\psi_t(z)=\cot_2({\til\psi_t(z)-\xi(t)}),\quad -\infty<t<T.\label{covering-whole}\EDE
We call $\til L_t$ and $\til \psi_t$ the inverted covering whole-plane Loewner hulls and maps, respectively, driven by $\xi$.
It is easily seen that for $-\infty<t<T$, the whole-plane Loewner objects driven by $\xi$ at time $t$, such as
$K_t$, $\vphi_t$, $L_t$, $\psi_t$, $\til K_t$, $\til\vphi_t$, $\til L_t$, $\til \psi_t$, are all determined by $e^i(\xi(s))$, $-\infty<s\le t$.

\vskip 3mm

From (\ref{covering-whole-eqn}) and (\ref{covering-whole-limit}), for any $z\in\C\sem \til L_t$, we have
$$R_\R(\til \psi_t(z)-(z-it))=\til \vphi_t(\lin z)-(\lin z+it)=\int_{-\infty}^t (\cot_2({\til\vphi_s(\lin z)-\xi(s)})-i)ds$$
\BGE =\int_{-\infty}^t i\Big(\frac{e^i(\til\vphi_s(\lin z))+e^i(\xi(s))}{e^i(\til\vphi_s(\lin z))-e^i(\xi(s))}-1\Big)ds=
\int_{-\infty}^t \frac{2i e^i(\xi(s)) }{\vphi_s(e^i(\lin{z}))-e^i(\xi(s))}\,ds.\label{est-4}\EDE
Suppose that $(1+C_{\cal H})e^t|e^{iz}|\le 1/2$. Since $|e^i(\lin{z})|=1/|e^i(z)|$, so for any $s\in(-\infty,t]$,
\BGE e^{-s}|e^i(\lin z)|-(1+C_{\cal H})\ge \frac{e^{-s}}{2|e^i(z)|}.\label{est-3}\EDE
Note that $e^i(\lin z)\in\C\sem K_t\subset \C\sem K_s$, $-\infty<s\le t$. From (\ref{est-2}), for $s\in(-\infty,t]$, we have
$|\vphi_s(e^i(\lin{z}))-e^{-s}e^i(\lin z)|\le C_{\cal H}$, which together with (\ref{est-3}) implies that
$$|\vphi_s(e^i(\lin{z}))-e^i(\xi(s))|\ge |\vphi_s(e^i(\lin{z}))|-1\ge |e^{-s}e^i(\lin z)|-C_{\cal H}-1\ge \frac{e^{-s}}{2|e^i(z)|}.$$
From (\ref{est-4}) and the above formula, we have
\BGE
|\til\psi_t(z)-(z-it)|\le 4(1+C_{\cal H})e^t|e^{iz}|,\quad\mbox{
if }(1+C_{\cal H})e^t|e^{iz}|\le 1/2.\label{tilpsi-near-infty}\EDE

\subsection{Carath\'eodory topology}
The following definition is about the convergence
of domains in Carath\'eodory topology.

\begin{Definition} Suppose $D_n$ is a sequence of domains and $D$ is
a domain. We say that $(D_n)$ converges to $D$, denoted by $D_n\dto
D$, if for every $z\in D$, $\dist^\#(z,\pa^\# D_n)\to \dist^\#
(z,\pa^\#
D)$. This is equivalent to the followings: \\
(i) every compact subset of $D$ is contained in all but finitely
many $D_n$'s; and\\ (ii) for every point $z_0\in\pa^\# D$,
$\dist^\#(z_0,\pa^\# D_n)\to 0$ as $n\to\infty$.
\end{Definition}

A sequence of domains may converge to two different domains. For
example, let $D_n=\C\sem((-\infty,n])$. Then $D_n\dto\HH$, and
$D_n\dto -\HH$ as well. But two different limit domains of the same
domain sequence must be disjoint from each other, because if they
have nonempty intersection, then one contains some boundary point of
the other, which implies a contradiction.

Suppose $D_n\dto D$, and for each $n$, $f_n$ is a $\ha\C$ valued
function on $D_n$, and $f$ is a $\ha\C$ valued function on $D$. We
say that $f_n$ converges to $f$ locally uniformly in $D$, or
$f_n\luto f$ in $D$, if for each compact subset $F$ of $D$, $f_n$
converges to $f$ in the spherical metric uniformly on $F$. If every
$f_n$ is analytic (resp.\ harmonic), then $f$ is also analytic
(resp.\ harmonic).

\begin{Lemma} Suppose $D_n\dto D$, $f_n$ maps $D_n$ conformally onto
some domain $E_n$ for each $n$, and $f_n\luto f$ in $D$. Then either
$f$ is constant on $D$, or $f$ maps $D$ conformally onto some domain
$E$. And in the latter case, $E_n\dto E$ and $f_n^{-1}\luto f^{-1}$
in $E$. \label{domain convergence}
\end{Lemma}

This lemma is similar to Theorem 1.8, the Carath\'eodory kernel
theorem, in \cite{Pom-bond}, and the proof is also similar.

Recall that $\cal H$ is the set of all bounded interior hulls in $\C$.
For every sequence $(H_n)$ in $\cal H$, there is at most one $H\in\cal H$ such that
$\C\sem H_n\dto \C\sem H$ because if we also have $\C\sem H_n\dto \C\sem H'$ for some $H'\in\cal H$,
then from $(\C\sem H')\cap(\C\sem H)\ne \emptyset$ we conclude that $\C\sem H'=\C\sem H$, and so
 $H'=H$. We write $H_n\hto H$ for $\C\sem H_n\dto \C\sem H$.
 We will define a metric $d_{\cal H}$ on $\cal H$ such that $H_n\to H$ w.r.t.\ $d_{\cal H}$ iff
 $H_n\hto H$.

Recall that for each $H\in\cal H$, $\phi_H$ maps $\ha\C\sem H$
conformally onto $\ha\C\sem\{|z|\le \rad(H)\}$ such that $\phi_H(\infty)=\infty$ and $\phi_H'(\infty)=1$.
So $\phi_H^{-1}$ is defined on $\{|z|>\rad(H)\}$. For $H_1,H_2\in\cal H$, let
$$d_{\cal H}^\vee(H_1,H_2)=|\rad(H_1)-\rad(H_2)|+$$
\BGE +\sum_{m=1}^\infty 2^{-m}\sup\{|\phi_{H_1}^{-1}(z)-\phi_{H_2}^{-1}(z)|: |z|\ge({\rad(H_1)}\vee
\rad(H_2))+\frac 1m\}.\label{dvee}\EDE
It is clear that $d_{\cal H}^\vee (H_1,H_2)=d_{\cal H}^\vee (H_2,H_1)\ge 0$, and
$d_{\cal H}^\vee (H_1,H_2)=0$ iff $H_1=H_2$. From (\ref{phi_H-1(z)-z}) we have $d_{\cal H}^\vee(H_1,H_2)<\infty$.
But $d_{\cal H}^\vee$ may not satisfy the triangle inequality.

We now define a metric $d_{\cal H}$ from $d_{\cal H}^\vee$
such that for $H_1,H_2\in\cal H$,
\BGE d_{\cal H}(H_1,H_2)=\inf\{\sum_{k=1}^n d_{\cal
H}^\vee (F_{k-1},F_k): F_0=H_1, F_n=H_2, F_k\in {\cal H}, 0\le k\le n,
n\in\N\}.\label{dcal}\EDE
It is clear that  $0\le d_{\cal H} (H_1,H_2)=d_{\cal H} (H_2,H_1)\le d_{\cal H}^\vee (H_1,H_2)<\infty$
and $d_{\cal H} $   satisfies the triangle inequality. We need to check that $d_{\cal H}(H_1,H_2)=0$ if and only if $H_1=H_2$.
The ``if'' part is clear because $d_{\cal H} (H_1,H_2) \le d_{\cal H}^\vee (H_1,H_2)$. For the ``only if'' part, we prove by
contradiction. Suppose $H_1\ne H_2$ and $d_{\cal H}(H_1,H_2)=0$. If there are $F_k\in\cal H$, $0\le k\le n$, such
that $F_0=H_1$ and $F_n=H_2$, then from (\ref{dvee}) we have
$$\sum_{k=1}^n d_{\cal H}^\vee (F_{k-1},F_k)\ge \sum_{k=1}^n |\rad(F_{k-1})-\rad(F_k)|\ge |\rad(H_1)-\rad (H_2)|.$$
So we have $|\rad(H_1)-\rad (H_2)|\le d_{\cal H}(H_1,H_2)=0$, which implies that $\rad(H_1)=\rad(H_2)$. Let  $r=\rad(H_1)$.
Since  $H_1\ne H_2$, so $\phi_{H_1}^{-1}\ne \phi_{H_2}^{-1}$ on $\{|z|>r\}$. Thus, there is $m\in\N$ such that
\BGE \sup\{|\phi_{H_1}(z)-\phi_{H_2}(z)|:|z|\ge r+\frac 1m\}\ge \frac 1m.\label{1m}\EDE
Since $d_{\cal H}(H_1,H_2)=0$, so there are $F_k\in\cal H$, $0\le k\le n$, such
that $F_0=H_1$ and $F_n=H_2$, and
\BGE \sum_{k=1}^n d_{\cal H}^\vee (F_{k-1},F_k)<\frac{2^{-2m}}{2m}.\label{contra-1}\EDE
For any $1\le j\le n$, since
$$\sum_{k=1}^j |\rad(F_{k-1})-\rad(F_k)|\le \sum_{k=1}^j d_{\cal H}^\vee (F_{k-1},F_k)<\frac{2^{-2m}}{2m}
<\frac 1{2m},$$
so $\rad(F_j)\le \rad(F_0)+\frac 1{2m}=r+\frac 1{2m}$. Thus, from (\ref{dvee}) and (\ref{1m})
$$\sum_{k=1}^n d_{\cal H}^\vee (F_{k-1},F_k)\ge
\sum_{k=1}^n 2^{-2m} \sup\{|\phi_{F_{k-1}}^{-1}(z)-\phi_{F_k}^{-1}(z)|: |z|\ge({\rad(F_{k-1})}\vee
\rad(F_k))+\frac 1{2m}\}$$
$$\ge 2^{-2m} \sum_{k=1}^n \sup\{|\phi_{F_{k-1}}^{-1}(z)-\phi_{F_k}^{-1}(z)|: |z|\ge r+\frac 1m\}$$
$$\ge  2^{-2m}\sup\{|\phi_{H_1}^{-1}(z)-\phi_{H_2}^{-1}(z)|: |z|\ge r+\frac 1m\}\ge \frac{2^{-2m}}m,$$
which  contradicts (\ref{contra-1}). Thus, $d_{\cal H}$ is a metric on $\cal H$.

Suppose $H_n\to H$ w.r.t.\ $d_{\cal H}$. Then we have $\rad(H_n)\to \rad (H)$ and $\phi_{H_n}^{-1}$ converges
to $\phi_{H}^{-1}$ uniformly on $\{|z|\ge \rad(H)+\eps\}$ for any $\eps>0$. Thus, $\{|z|> \rad(H_n)\}\dto
\{|z|>\rad(H)\}$ and $\phi_{H_n}^{-1}\luto\phi_{H}^{-1}$ in $\{|z|>\rad(H)\}$. From Lemma \ref{domain convergence},
we have $\C\sem H_n=\phi_{H_n}^{-1}(\{|z|> \rad(H_n)\})\dto \phi_{H }^{-1}(\{|z|> \rad(H )\})=\C\sem H$, i.e.,
$H_n\hto H$. On the other hand, suppose $\C\sem H_n\dto \C\sem H$. We will show that $H_n\to H$ w.r.t.\ $d_{\cal H}$.
For this purpose, we will derive a stronger result.

\vskip 3mm

Recall that ${\cal H}_0$ is the set of hulls in $\cal H$ that
contains $0$. Since $d_{\cal H}(H_n,H_0)\to 0$ implies $
H_n\hto H_0$, so ${\cal H}_0$ is a closed subset of $({\cal H},d_{\cal H})$.
For any $F\in{\cal H}$, let ${\cal H}(F)=\{H\in{\cal H}: H\subset
F\}$,  ${\cal H}_0(F)={\cal H}(F)\cap{\cal H}_0$, and
 ${\cal H}_0^x(F)=\{H\in {\cal H}_0(F):\ccap(F)\ge x\}$, $x\in\R$. Then
${\cal H}_0(F)$ and ${\cal H}_0^x(F)$ are closed subsets of ${\cal H}(F)$
because $\rad$ is continuous w.r.t.\ $d_{\cal H}$. The hulls in ${\cal H}_0^x(F)$
are non-degenerate because they have finite capacities.  If
$\sigma$ is a Jordan curve in $\C$, we write ${\cal H}(\sigma)$,
${\cal H}_0(\sigma)$, and ${\cal H}_0^x(\sigma)$  for ${\cal
H}(H(\sigma))$, ${\cal H}_0(H(\sigma))$, and ${\cal
H}_0^x(H(\sigma))$, respectively.

\begin{Lemma} ${\cal H}(F)$, ${\cal H}_0(F)$ and ${\cal H}_0^x(F)$
are all compact subsets of $({\cal H},d_{\cal H})$.\label{compact}
\end{Lemma}
{\bf Proof.} Let $r_F=\max\{|z|:z\in F\}$. Then for any $H\in{\cal
H}(F)$,  $\rad(H)\le \rad(F)\le r_F$.
Suppose $(H_n)$ is a sequence in ${\cal H}(F)$. By passing to a
subsequence, we may assume that $\rad(H_{n})\to r_0\in[0,r_F]$. For
each $n\in\N$, let $g_n(z)=\phi_{H_n}^{-1}({z})-z$ for
$\ha\C\sem\lin{\B(0;\rad(H_n))}$. Then $g_n$ is analytic. From
(\ref{phi_H-1(z)-z}), $|g_n|$ is bounded by $C_F:=(1+C_{\cal
H})r_F$. Since $\rad(H_n)\to r_0$, so
$\ha\C\sem\lin{\B(0;\rad(H_n))}\dto\ha\C\sem\lin{\B(0;r_0)}$. Since
$(g_n)$ is a normal family, by passing to a subsequence, we may
assume that $(g_{n})\luto g_0$   in
$\ha\C\sem\lin{\B(0;r_0)}$. Then $|g_0|$ is also bounded by $C_F$ on
$\ha\C\sem\lin{\B(0;r_0)}$. Let $f_0(z)=g_0(z)+z$ for $|z|>r_0$. Then
$f_0(z)-z$ is bounded, and $\phi_{H_{n}}^{-1}(z)=g_{n}(z)+z\to f_0(z)$
uniformly on $\{|z|\ge
 r\}$ for any $r>r_0$. From Lemma \ref{domain convergence}, $f_0$ is either
constant or a conformal map on $\ha\C\sem\lin{\B(0;r_0)}$. Since
$f_0(z)-z$ is bounded, so $f_0$ cannot be constant. Thus, $f_0$ is a
conformal map, and $\ha\C\sem H_n\dto f_0(\ha\C\sem\lin{\B(0;r_0)})$.
Since $f_0(z)-z$ is bounded, so $\infty=f_0(\infty)\in
f_0(\ha\C\sem\lin{\B(0;r_0)})$ and $f_0'(\infty)=\lim_{n\to\infty}(\phi_{H_{n}}^{-1})
'(\infty)= 1$. Since $f_0(\ha\C\sem\lin{\B(0;r_0)})$ is
simply connected, so its complement in $\ha\C$ is some $H_0\in\cal
H$. Thus, $f_0=\phi_{H_0}^{-1}$ and $\rad(H_0)=r_0$.

We now have proved that, by passing to a subsequence, we have $\rad(H_{n})\to
\rad(H_0)$ and $\phi_{H_{n}}^{-1}\luto \phi_{H_0}^{-1}$
in $\ha\C\sem\lin{\B(0;r_0)}$. Moreover, for any $|z|>\rad(H_{n})\vee \rad(H_0)$,
$$|\phi_{H_{n}}^{-1}(z)-\phi_{H_0}^{-1}(z)|\le |\phi_{H_{n}}^{-1}(z)-z|
+|\phi_{H_0}^{-1}(z)-z|\le 2C_F.$$ Given $\eps>0$, there is $M\in\N$
such that $2^{-M}(2C_F)<\eps/3$. There is $N\ge M$ such that, for
$n\ge N$, $|\rad(H_{n})-\rad(H_0)|<(\eps/3)\wedge (1/N)$, and
$|\phi_{H_{n}}^{-1}(z)-\phi_{H_0}^{-1}(z)|\le \eps/3$ for any
$|z|\ge \rad(H_0)+ 1/N$. Thus,
$$d_{\cal H}(H_{n},H_0)\le d_{\cal H}^\vee (H_{n},H_0)=
|\rad(H_{n})-\rad(H_0)|+$$
$$+\sum_{m=1}^\infty 2^{-m}\sup\{|\phi_{H_{n}}^{-1}(z)-\phi_{H_0}^{-1}(z)|
:|z|\ge (\rad(H_{n})\vee\rad(H_0))+1/m\}$$
$$<\frac\eps 3+\frac \eps 3\sum_{m=1}^N
2^{-m}+2C_F\sum_{m=N+1}^\infty 2^{-m}<\frac\eps 3+\frac\eps
3+\frac\eps 3=\eps.$$ So we have $d_{\cal H}(H_{n}, H_0)\to 0$. Thus,
${\cal H}(F)$ is compact. The rest part of the lemma follows from
the facts that ${\cal H}_0$ and $\{H\in{\cal H}:\ccap(H)\ge x\}$ are
closed. $\Box$

\vskip 3mm

Suppose $H_n\hto H$. Choose $r_0\in(0,\infty)$ such that $H\subset \{|z|<r_0\}$.
Then $\{|z|=r_0\}$ is a compact subset of $\C\sem H$.
Let $\del=\dist(H,\{|z|=r_1\})>0$. Choose $z_0\in\pa H=\pa(\C\sem H)$.
Since $\C\sem H_n\dto\C\sem H$, so there is $N\in\N$ such that if $n\ge N$,
then $\{|z|=r_0\}\subset \C\sem H_n$ and $\dist(z_0,\pa(\C\sem H_n))<\del$,
which implies that $H_n\cap\{|z|=r_0\}=\emptyset$ and
$H_n\cap \{|z|<r_0\}\ne \emptyset$. Since $H_n$ is connected, so
$H_n\subset  \{|z|<r_0\}$ if $n\ge N$. Thus, $\{H_n:n\ge N\}\subset {\cal H}(\{|z|\le r_0\})$.
From Lemma \ref{compact}, $\{H_n:n\in\N\}$ is a pre-compact set. Assume that $H_n\not\to H$
w.r.t.\ $d_{\cal H}$. Then there is $\eps>0$ and a subsequence $(H_{n_k})$ of $(H_n)$ such
that $d_{\cal H}(H_{n_k},H)\ge \eps$ for any $k\in\N$. By passing to a subsequence, we may
assume that $H_{n_k}\to H'$ w.r.t.\ $d_{\cal H}$.
Then $H'\ne H$ and $H_{n_k}\hto H'$. Since $H_n\hto H$,
so the subsequence $H_{n_k}\hto H$ as well. Then we must have $H'=H$, which is
a contradiction. So $H_n\to H$ w.r.t.\ $d_{\cal H}$. So the topology on $\cal H$ generated by $d_{\cal H}$
agrees with Carath\'eodory topology. From (\ref{dvee}) and (\ref{dcal}) we see that if $H_n\to H$ w.r.t.\ $d_{\cal H}$,
then $\rad(H_n)\to \rad (H)$ and  $\phi_{H_n}^{-1}\luto \phi_H^{-1}$  in $\{|z|>\rad(H)\}$.
From Lemma \ref{domain convergence} we have $\phi_{H_n}\luto \phi_H$ in $\C\sem H$. Recall that for
every non-degenerate interior hull $H$, $\vphi_H=\rad(H)^{-1}\phi_H$ maps $\C\sem H$ conformally
onto $\{|z|>1\}$ and $\psi_H=R_\TT\circ \vphi_H\circ R_\TT$ maps $\C\sem R_TT(H)$ conformally onto $\D$.

\begin{Lemma} Let $\alpha$ be a Jordan curve, $F$ be a compact subset of $\C\sem H(\alpha)$, and $b\in\R$.
If $(H_n)_{n=1}^\infty$ is a sequence in ${\cal H}^{b}(\alpha)$, then there is $H\in {\cal H}^{b}(\alpha)$
and a subsequence $(H_{n_k})$ of $(H_n)$ such that
$\vphi_{H_{n_k}}\to \vphi_H$ uniformly on $F$, and $\psi_{H_{n_k}}\to\psi_H$ uniformly on $R_\TT(F)$.
\label{compact-2}
\end{Lemma}
{\bf Proof.} From Lemma \ref{compact},  there is $H\in {\cal H}^{b}(\alpha)$
and a subsequence $(H_{n_k})$ of $(H_n)$ such that $H_{n_k}\hto H$. Then  $\rad(H_{n_k})\to \rad(H)$ and
$\phi_{H_{n_k}}\luto \phi_{H}$ in $\C\sem H$. Since $H\subset H(\alpha)$ and
$F$ is a compact subset of $\C\sem H(\alpha)$, so $F$ is also a compact subset of $\C\sem H $.
Thus, $\phi_{H_{n_k}}\to \phi_H$ uniformly on $F$. Since
$\rad(H_{n_k})\to \rad(H)\ge e^{b}$, so
$\vphi_{H_{n_k}}\to \vphi_H$ uniformly on $F$, and $\psi_{H_{n_k}}\to\psi_H$ uniformly on $R_\TT(F)$. $\Box$

\section{Continuous LERW} \label{Continuous LERW}
\subsection{Continuous boundary LERW} \label{section-radial-LERW} Let
$\Om$ be an almost-$\D$ domain, and $p\in\Om$. Let $\til\Om=(e^i)^{-1}(\Om)$ and $\til
p=(e^i)^{-1}(p)$. For $\xi\in C([0,T))$, let $\psi^\xi_t$ (resp.\
$\til\psi^\xi_t$) and $L^\xi_t$ (resp.\ $\til L^\xi_t$), $0\le t<T$,
denote the radial (resp.\ covering radial) Loewner maps and hulls,
respectively, driven by $\xi$. Suppose $L^\xi_t\subset\Om\sem\{p\}$,
that is, $\til L^\xi_t\subset \til\Om\sem\til p$. Then $\Om\sem
L^\xi_t$ is a finitely connected subdomain of $\Om$, and contains
$p$. Let $\Om^\xi_t=\psi^\xi_t(\Om\sem L^\xi_t)$,
$\til\Om^\xi_t=(e^i)^{-1}(\Om^\xi_t)=\til\psi^\xi_t(\til\Om\sem\til
L^\xi_t)$, $p^\xi_t=\psi^\xi_t(p)$, and $\til
p^\xi_t=\til\psi^\xi_t(\til p)$. Then $\Om^\xi_t$ is also an almost-$\D$
domain, $p^\xi_t\in\Om^\xi_t$, and $\til
p^\xi_t\subset\til\Om^\xi_t$. For a finitely connected domain $D$ and $z_0\in D$,
let $G(D,z_0;\cdot)$ denote the Green function in $D$ with the pole at $z_0$.
Let \BGE J^\xi_t=G(\Om\sem
L^\xi_t,p;\cdot)\circ(\psi^\xi_t)^{-1}=G(\Om^\xi_t,p^\xi_t;\cdot),\label{jxit}\EDE and $\til
J^\xi_t=J^\xi_t\circ e^i$. Then $\til J^\xi_t$ is
harmonic on $\til\Om^\xi_t\sem\til p^\xi_t$, and vanishes on $\R$,
so can be extended harmonically across $\R$ by Schwarz reflection
principle. Let $X^\xi(t)=(\pa_x\pa_y/\pa_y)\til J^\xi_t(\xi(t))$.
The following theorem is similar to Theorem 3.1 in
\cite{LERW}. The difference is that here we use radial Loewner equation.
We will prove the theorem in Section \ref{Section-The-radial-equation}.

\begin{Theorem} (i) For any $f\in C([0,\infty))$ and $\lambda\in\R$,
the equation \BGE \xi(t)=f(t)+\lambda\int_0^t
X^\xi(s)ds\label{pre-LERW-eqn-radial}\EDE has a solution $\xi(t)$ on
$[0,a]$ for some $a>0$.\\ (ii) If for $j=1,2$, $\xi_j$ solves
(\ref{pre-LERW-eqn-radial}) for $0\le t<T_j$, and $T_j>0$, then
there is $S>0$ such that $\xi_1(t)=\xi_2(t)$ for $0\le t\le S$.
\label{Thm-pre-LERW-eqn-radial}
\end{Theorem}

\no {\bf Remark} The statement of the above theorem is enough for
the use of this paper. In fact, the followings are true. Equation
(\ref{pre-LERW-eqn-radial}) has a unique maximal solution $\xi_f(t)(t)$,
$0\le t<T_f$, for some $T_f>0$. Here we call a solution maximal if it
can not be extended. Moreover,
for any $a\ge 0$, $\{f\in C([0,\infty)):T_f>a\}$
is open w.r.t.\ $\|\cdot\|_{0,a}$, and $f\mapsto \xi_f$ is
$(\|\cdot\|_{0,a},\|\cdot\|_{0,a})$ continuous on $\{T_f>a\}$. Let
$\lambda=2$ and $f(t)=\sqrt 2 B(t)$, where $B(t)$ is a Brownian
motion. Let $\xi(t)$, $0\le t<T$, be the maximal solution to
(\ref{pre-LERW-eqn-radial}). For $0\le t<T$, let
$$u(t)=\int_0^t \pa_y\til J^\xi_s(\xi(s))^2ds.$$
One can prove that $(L^\xi_{u^{-1}(t)},0\le t<u(T))$ has the same
distribution as the continuous $\LERW(\Om;1\to p)$ defined in
\cite{LERW}. The proof is similar to that of Theorem 3.2 in
\cite{LERW}. So the radial Loewner equation plays an equivalent role
as chordal Loewner equation in defining a continuous boundary LERW.

\subsection{Continuous interior LERW}\label{Def-of-quasi} Let $D$ be
a finitely connected domain that contains $0$. Fix $z_e\in
D\sem\{0\}$. Let $\Om=R_\TT(D)$, $p=R_\TT(z_e)$,
$\til\Om=(e^i)^{-1}(\Om)$, and $\til p=(e^i)^{-1}(p)$. Let $\xi\in
C((-\infty,T))$. We use $K^\xi_t$ (resp.\ $L^\xi_t$, $\til L^\xi_t$) and
$\vphi^\xi_t$ (resp.\ $\psi^\xi_t$, $\til\psi^\xi_t$), $0\le t<T$, to
denote the whole-plane (resp.\ inverted whole-plane, inverted
covering whole-plane) Loewner hulls and maps, respectively, driven
by $\xi\in C((-\infty,T))$. Recall that if $\xi\in C([0,T))$, we use
$\psi^\xi_t$, $\til\psi^\xi_t$, $L^\xi_t$ and $\til L^\xi_t$ to
denote the radial Loewner objects driven by $\xi$. But this will not cause
ambiguity.

If for some $t< T$, $K^\xi_t\subset D\sem\{z_e\}$,
that is, $L^\xi_t\subset \Om\sem\{p\}$ or $\til
L^\xi_t\subset\til\Om\sem\til p$, then let
$\Om^\xi_t=\psi^\xi_t(\Om\sem L^\xi_t)$, $p^\xi_t=\psi^\xi_t(p)$,
$\til \Om^\xi_t=(e^i)^{-1}(\Om^\xi_t)=\til\psi^\xi_t(\til\Om\sem\til L^\xi_t)$, and
$\til p^\xi_t=(e^i)^{-1}(p^\xi_t)=\til\psi^\xi_t(\til p)$. Then
$\Om^\xi_t$ is an almost-$\D$ domain that contains $p^\xi_t$.

Let \BGE J^\xi_t=G(\Om\sem
L^\xi_t,p;\cdot)\circ(\psi^\xi_t)^{-1}=G(\Om^\xi_t,p^\xi_t;\cdot),\label{jxit-whole}\EDE and
$\til J^\xi_t=J^\xi_t\circ e^i$. Then $\til J^\xi_t$ is a
positive harmonic function in $\til \Om^\xi_t\sem\til p^\xi_t$, and
vanishes on $\R$. From Schwarz reflection principle, $\til J^\xi_t$
extends harmonically across $\R$.  Let
$X^\xi(t)=(\pa_x\pa_y/\pa_y)\til J^\xi_t(\xi(t))$.
Recall that $\psi^\xi_t=R_{\TT}\circ \vphi^\xi_t\circ R_{\TT}$. It is easy to check that
the $X^\xi(t)$ here agrees with that in (\ref{X-1st-def}).

For $a\in\R$, let $\T_a$ denote the topology on $C((-\infty,a])$
generated by $\|\cdot\|_{b,a}$, $b\le a$. For $f_1,f_2\in
C((-\infty,a])$, we write $f_1\aequ f_2$ if
$e^i(f_1(t))=e^i(f_2(t))$ for any $t\le a$. Let $\T^\TT_a$ be the
set of $S\in \T_a$ such that $\pi_a^{-1}(\pi_a(S))=S$, where $\pi_a$
is the projection map from $C((-\infty,a])$ onto
$C((-\infty,a])/\aequ$. Then $\T^\TT_a$ is also a topology on
$C(\R)$. Let ${\cal F}^0_a$ be the $\sigma$-algebra generated by
$\T^\TT_a$. Then ${\cal F}^0_a$ agrees with the $\sigma$-algebra
generated by the functions $f\mapsto e^i(f(t))$, $t\in(-\infty, a]$.
The proposition and theorem below will be proved in Section \ref{section-proof-whole-plane}.

\begin{Proposition} If $L^\xi_a\subset \Om\sem\{p\}$, then the improper integral
$\int_{-\infty}^a X^\xi(t)dt$ converges. \label{improper}
\end{Proposition}

\begin{Theorem} Fix $\lambda\in\R$. For any $f\in C(\R)$, the equation
\begin{equation}\xi(t)=f(t)+\lambda\int_{-\infty}^t X^\xi(s)ds\label{pre-interior-LERW-eqn}
\end{equation}
has a unique maximal solution $\xi_f\in C((-\infty,T_f))$ for some
$T_f\in(-\infty,+\infty]$. Moreover, \begin{itemize}
\item[$\mathrm{(i)}$] for any $a\in\R$,
$\{f\in C(\R):T_f>a\}\in\T^\TT_a$, and $f\mapsto \xi_f $ is
$(\T^\TT_a,\T^\TT_a)$-continuous on $\{f\in C(\R):T_f>a\}$;
\item[$\mathrm{(ii)}$] there
does not exist a Jordan curve $\alpha$ such that $\bigcup_{t<T_f}
K^{\xi_f}_t\subset H(\alpha)\subset D\sem\{z_e\}$.
\end{itemize}
 \label{existance-and-uniqueness-Theorem-Interior}
\end{Theorem}

Let $B_+(t)$ and $B_-(t)$, $0\le t<\infty$, be two independent
Brownian motions. Let $\bf x$ be a random variable that is uniformly
distributed on $[0,2\pi)$, and independent of $B_\pm(t)$. For $\kappa>0$
and
$t\in\R$, let $B^{(\kappa)}_\R(t)={\bf x}+\sqrt \kappa B_{\sign(t)}(|t|)$.
Then the whole-plane Loewner hulls driven by $B^{(\kappa)}_\R(t)$ are called
the whole-plane SLE$_\kappa$ hulls. We will be particularly interested
in the case that $\kappa=2$.

Let
$(\F_t)$ be the usual augmentation of $(\F^{0}_t)$ w.r.t.\ the
distribution of $B^{(2)}_\R$. So $(\F_t)$ is right-continuous.
Let $\F_\infty=\vee_{t\in\R} \F_t$.
Suppose $S$ is a finite $(\F_t)_{t\in\R}$-stopping time. Then for
any $t\ge 0$, $S+t$ is an $(\F_t)_{t\in\R}$-stopping time. So we
have a filtration $({\cal F}_{S+t})_{t\ge 0}$. For $t\ge 0$, let
$B_S(t):=(B^{(2)}_\R(S+t)-B^{(2)}_\R(S))/\sqrt 2$. It is well known
that $(B_S(t), t\ge 0)$ is an $({\cal F}_{S+t})_{t\ge 0}$-Brownian
motion.

 Suppose $\xi\in C((-\infty,T))$ is the maximal solution to
(\ref{pre-interior-LERW-eqn}) with $f=B^{(2)}_\R$ and $\lambda=2$.
Then we call $(K^\xi_t,0\le t<T)$ a continuous interior LERW process in $D$
from $0$ to $z_e$, and let it be denoted by $\LERW(D;0\to z_e )$.
From Theorem \ref{existance-and-uniqueness-Theorem-Interior} (i),
$T$ is an $(\F_t)_{t\in\R}$-stopping time, and $(e^i(\xi(t)))$ is
$(\F_t)$-adapted. So for any fixed $a\in\R$, $(\xi(a+t)-\xi(a),0\le
t<T-a)$ is $(\F_{a+t})_{t\ge 0}$-adapted. Since
$K^\xi_t,L^\xi_t,\til L^\xi_t,\vphi^\xi_t,\psi^\xi_t,\til\psi^\xi_t$
are determined by $e^i\circ\xi(s)$, $-\infty<s\le t$, so they are all $(\F_t)_{t\in\R}$-adapted.
Note that in general $(\xi(t))$ is  not $(\F_t)_{t\in\R}$-adapted.

Let $R=\dist(0;\pa D\cup\{z_e\})>0$. Fix $r\in (0,R)$. From Theorem
\ref{existance-and-uniqueness-Theorem-Interior} (ii), there is
$t_0\in(-\infty,T)$ such that $K^\xi_{t_0}\not\subset \B(0;r)$. Then
$T>t_0=\ccap(K^\xi_{t_0})\ge \ln(r/4)$. So $T\ge \ln(R/4)$. Fix
$a\in(-\infty,\ln(R/4))$. Then $a<T$. Let $T_a=T-a$ and
$\xi_a(t)=\xi(a+t)-\xi(a)$ for $0\le t<T_a$. Then $T_a$ is an
$({\cal F}_{a+t})_{t\ge 0}$-stopping time, $(\xi_a(t))$ and
$(X^\xi(a+t))$ are $({\cal F}_{a+t})_{t\ge 0}$-adapted. Recall that
$B_a(t)=(B^{(2)}_\R(a+t)-B^{(2)}_\R(a))/\sqrt 2$ is an
$({\cal F}_{a+t})_{t\ge 0}$-Brownian motion, so
$\xi_a$ solves the $({\cal F}_{a+t})_{t\ge 0}$-adapted SDE: \BGE d\xi_a(t)=\sqrt
2dB_a(t)+2X^\xi({a+t})dt,\quad 0\le t<T_a.\label{SDE-1}\EDE

From Girsanov's theorem (\cite{RY}) and the existence of the radial SLE$_2$ trace, one can easily show
that the interior Loewner chain $K^\xi_t$, $-\infty<t<T$, is a.s.\ generated by a simple curve $\beta(t)$,
$-\infty\le t<T$, with $\beta(-\infty)=0$. We call such $\beta$ an $\LERW(D;0\to z_e )$ curve.

Suppose  $z_0\ne z_e\in D$. If $z_0\in\C$, we define $\LERW(D;z_0\to
z_e)$ to be the image of $\LERW(A_{z_0}^{-1}(D);0\to
A_{z_0}^{-1}(z_e))$ under the map $A_{z_0}$. If $z_0=\infty$, we define
$\LERW(D;z_0\to z_e)$ to be the image of $\LERW(W(D);0\to W(z_e))$
under the map $W(z)=1/z$.

\vskip 3mm

\no {\bf Remark} A continuous $\LERW(\ha\C:0\to\infty)$ has the same
distribution as a whole-plane SLE$_2$. This can be seen from the fact that
 $X^\xi(t)\equiv 0$.

\subsection{Conformal invariance} \label{section-conformal}
\begin{Theorem} Let $D$ be a finitely connected domain, and $z_0,z_e\in D$ with
$z_0\ne z_e$. Let $(K_t, -\infty<t<T)$ be an $\LERW(D;z_0\to z_e )$ process.
Suppose $V$ maps $D$ conformally onto another finitely connected domain $D^*$. Then after a
time-change, $(V(K_{t}),-\infty<t<T)$ has the same distribution as $(K^*_t,-\infty<t<T^*)$, which is
an $\LERW(D^*;z^*_0\to z^*_e)$ process, where $z^*_0=V(z_0)$ and
$z^*_e=V(z_e)$. \label{thm-conformal-invariance}
\end{Theorem}
{\bf Proof.} WLOG, assume $z_0= z^*_0=0$. Let $\kappa=2$.
From the definition, $K_t=K^\xi_t$ for $-\infty<t<T$, where
$\xi(t)$, $-\infty<t<T$, is the maximal solution to the equation
\BGE\xi(t)=B^{(\kappa)}_\R(t)+\Big(3-\frac \kappa 2\Big)\int_{-\infty}^tX^\xi_sds\label{LERW-int-equation}.\EDE
Since $\kappa\le 4$, so from the property of SLE$_\kappa$ (c.f.\ \cite{RS-basic}), a.s.\  $V^{-1}(\infty)\not\in K_t$ for
any $t<T$. Since $V(0)=0$, so $(V(K_t),-\infty<t<T)$ is a.s.\ an interior Loewner
chain started from $0$  avoiding $\infty$. Let $u(t)=\ccap(V(K_t))$
for $-\infty<t<T$, and $T^*=u(T)$. Let $v(t)=u^{-1}(t)$ and
$K^*_t=V(K_{v(t)})$ for $-\infty<t<T^*$. So $(K^*_t)$ is a
time-change of $(V(K_{t}))$. We will prove that $(K^*_t,-\infty
<t<T^*)$ has the same distribution as an $\LERW(D^*;0\to z^*_e)$.
Since $(K^*_t)$ is parameterized by capacity, so from Proposition
\ref{Whole-plne-Proposition-I}, there is $\xi^*\in C((-\infty,
T^*))$ such that $K^*_t=K^{\xi^*}_t$ for $-\infty<t< T^*$. For
simplicity, we omit the superscripts $\xi$, and replace the
superscripts $\xi^*$ by $*$ for the whole-plane Loewner objects driven by $\xi$ or $\xi^*$,
respectively, in the rest of this proof.

Recall $\Om=R_\TT(D)$,
$\til\Om=(e^i)^{-1}(\Om)$, $p=R_\TT(z_e)$, $\til p=(e^i)^{-1}(p)$,  $\Omega_t
=\psi_t(\Omega\sem L_t)$, and $\til \Omega_t=(e^i)^{-1}(\Omega_t)$. We can define $\Om^* $, $\til\Om^* $,
$p^* $, $\til p^* $,  $\Omega^*_t $, and $\til \Omega^*_t$, similarly for $D^*$ and the driving function $\xi^*$.
Let $W=R_\TT\circ V\circ R_\TT$. Then $W$ maps $\Om$ conformally onto $\Om^*$, and
$W(p)=p^*$. There is $\til W$ that maps $\til\Om$ conformally onto
$\til\Om^*$ such that $e^i\circ \til W=W\circ e^i$.
Let \BGE W_t=\psi^*_{u(t)}\circ W\circ
\psi_{t}^{-1},\quad \til W_t=\til \psi^*_{u(t)}\circ \til W\circ
\til \psi_{t}^{-1},\quad -\infty<t<T. \label{Wt}\EDE Then $e^i\circ \til
W_t=W_t\circ e^i$, and $W_t$ (resp.\ $\til W_t$) maps $\Om_{t}$
(resp.\ $\til\Om_{t}$) conformally onto $\Om^*_{u(t)}$ (resp.\
$\til\Om^*_{u(t)}$). Since $W_t(z)\to\TT$ as $\Omega_t\ni z\to\TT$, and $\til W_t(z)\to\R$
as $\til\Omega_t\ni z\to\R$, so
from Schwarz reflection principle, $W_t$ (resp.\
$\til W_t$) extends conformally across $\TT$ (resp.\ $\R$). Since
 $W_t (K_{t+\eps}/K_{t})=K^*_{u(t+\eps)}/K^*_{u(t)}$ for
$-\infty<t<t+\eps< T$. So from (\ref{xi(t)-whole-plane}),
$W_t(e^i(\xi(t)))=e^i(\xi^*(u(t)))$. Since
$\dcap(K_{t+\eps}/K_{t})=\eps$ and $\dcap( K^*_{u(t+\eps)}/
K^*_{u(t)})=u(t+\eps)-u(t)$, so from (\ref{capacity}) we have, \BGE
u'(t)=|W_t'(e^i(\xi(t)))|^{2}=\til
W_t'(\xi(t))^{2},\quad -\infty<t<T.\label{v'}\EDE Now $e^i\circ \til
W_t(\xi(t))=W_t\circ e^i(\xi(t))=e^i(\xi^*(u(t)))$, so $\til
W_{t}(\xi(v(t)))$ is also the driving function of $(K^*_t)$. So we
may choose $\xi^*$ such that,  \BGE \xi^*(u(t))=\til
W_{t}(\xi(t)),\quad -\infty<t<T.\label{xi=W(xi)}\EDE

Differentiate the equality $\til W_t\circ
\til\psi_{t}(z)=\til\psi^*_{u(t)}\circ \til W(z)$ w.r.t\ $t$ for
$t\in (-\infty,T)$ and $z\in \til\Om\sem\til L_{t}$. From
(\ref{covering-whole}), (\ref{v'}), and (\ref{xi=W(xi)}), we have
$$\pa_t\til W_t(\til\psi_{t}(z))+\til
W_t'(\til\psi_{t}(z))\cot_2({\til\psi_{t}(z)-\xi(t)})$$
$$=u'(t)\cot_2({\til\psi^*_{u(t)}\circ
\til W(z)-\xi^*(u(t))})=\til W_t'(\xi(t))^2\cot_2({\til W_t\circ
\til\psi_{t}(z)-\til W_t(\xi(t))}).$$ Since $\til\psi_{t}$
maps $\til\Om\sem \til L_{t}$ onto $\til\Om_{t}$, so
 for any $w\in\til\Om_{t}$, $$\pa_t\til W_t(w)=\til W_t'(\xi(t))^2\cot_2({\til
W_t(w)-\til W_t(\xi(t))})-\til
W_t'(w)\cot_2({w-\xi(t)}).$$ Letting $w\to \xi(t)$ in
$\til\Om_t$,  we get \BGE \pa_t\til
W_t(\xi(v(t)))=-3\til W_t''(\xi(v(t))).\label{-3}\EDE

Since $W$ maps $\Om\sem L_{v(t)}$ conformally onto $\Om^*\sem
L^*_{t}$, and $W(p)=p^*$, so $G(\Om\sem
L_{t},p;\cdot)=G(\Om^*\sem L^*_{u(t)}, p^*;\cdot)\circ W$. Thus,
$J_{t}=J^*_{u(t)}\circ W_{t}$, and so $\til J_{t}=\til
J^*_{u(t)}\circ\til W_{t}$. Since $X(t)=(\pa_x\pa_y/\pa_y)\til
J_{t}(\xi(t))$, $X^*({u(t)})=(\pa_x\pa_y/\pa_y)\til
J^*_{u(t)}(\xi^*(u(t)))$, so from (\ref{xi=W(xi)}),
\BGE X(t)=\til W_{t}''(\xi(t))/\til W_{t}'(\xi(t))+\til
W_{t}'(\xi(t))X^*({u(t)}), \quad -\infty<t<T.\label{X(v(t))-mrX(t)}\EDE

We now want to apply It\^o's formula. The following non-rigorous argument illustrate the
idea of the proof. From (\ref{LERW-int-equation}), $\xi(t)$, $-\infty<t<T$, satisfies the SDE
\BGE d\xi(t)=d B^{(\kappa)}_\R(t)+\Big(3-\frac\kappa 2\Big)X(t) dt.\label{infty-SDE-1}\EDE
One may think of $B^{(\kappa)}_\R(t)$ as $\sqrt\kappa B(t)$. From (\ref{xi=W(xi)}), (\ref{-3}), and  It\^o's formula,
we have
$$ d\xi^*(u(t))=\til W_t'(\xi(t))d\xi(t)+\pa_t \til W_t(\xi(t))dt +\frac\kappa 2 \til W_t''(\xi(t))dt$$
\BGE =\til W_t'(\xi(t))d B^{(\kappa)}_\R(t)+\Big(3-\frac\kappa 2\Big)\Big(\til W_t'(\xi(t))X(t) dt
- \til W_t''(\xi(t))\Big)dt.\label{infty-SDE-2}\EDE
From (\ref{X(v(t))-mrX(t)}) we then have
$$d\xi^*(u(t))=\til W_t'(\xi(t))d B^{(\kappa)}_\R(t)+\Big(3-\frac\kappa 2\Big)\til W_t'(\xi(t))^2 X^*({u(t)})dt.$$
Finally, we use (\ref{v'}) to conclude that there is another copy of $B^{\kappa}_\R(t)$ such that
$$d\xi^*(t)=d B^{(\kappa)}_\R(t)+\Big(3-\frac\kappa 2\Big) X^*(t)dt,\quad -\infty< t<T^*.$$
So   $\xi^*(t)$ is a driving function for continuous LERW in $D^*$ from $0$ to $z_e^*$.
The  argument is not rigorous because $B^{(\kappa)}_\R(t)$ is not a Brownian motion in the usual sense, and
 It\^o's formula does not directly apply to   time intervals of the form $(-\infty,T)$. We have a way to solve these problems,
 which is to truncate the time-interval.

We will use the filtration $\F_t$, $t\in\R$, in Section \ref{Def-of-quasi}.
Suppose that $a$ is a finite $(\F_t)$-stopping time such that $a<T$ always holds.
Let $\F^a_t=\F_{a+t}$, $0\le t<\infty$. Then we have a new filtration $(\F^a_t)_{t\ge 0}$.
Let $T_a=T-a>0$. Then $T_a$ is an $(\F^a_t) $-stopping time.
Let $B_a(t)=(B^{(\kappa)}_\R(a+t)-B^{(\kappa)}_\R(a))/\sqrt\kappa$, $0\le t<\infty$.
Then $B_a(t)$ is an $(\F^a_t) $-Brownian motion. Let $ \xi_a(t)=\xi(a+t)-\xi(a)$
and $X_a(t)=X(a+t)$. Then $(\xi_a)$ and $(X_a)$ are both  $(\F^a_t) $-adapted, and
$ \xi_a(t)$, $0\le t<T_a$, satisfies the  $(\F^a_t) $-adapted SDE:
\BGE d \xi_a(t)=\sqrt\kappa dB_a(t)+\Big(3-\frac\kappa 2\Big)X_a(t)dt.\label{a-SDE-1}\EDE
Let $u_a(t)=u(a+t)-u(a)$, $0\le t<T_a$. Then $u_a$ is continuous and increasing on $[0,T_a)$,
and $u_a(0)=0$. Let $\xi^*_b(t)=\xi^*(b+t)-\xi^*(b)$ for $b\in(-\infty,T^*)$ and $t\in[0,T^*-b)$.
Let $$\til W_{a,t}=A_{\xi^*(u(a))}^{-1}\circ \til W_{a+t}\circ
A_{\xi(a)}.$$ Then $(\til W_{a,t})$ is also $(\F^a_t) $-adapted.
From (\ref{xi=W(xi)}), (\ref{v'}), and (\ref{-3}) we have
\BGE\xi^*_{u(a)}(u_a(t))=\til W_{a,t}(\xi_a(t)).\label{xi=W(xi)-tr}\EDE
\BGE u_a'(t)=\til W_{a,t}'(\xi_a(t))^2.\label{v'-tr}\EDE
\BGE \pa_t \til W_{a,t}(\xi_a(t))=-3\til W_{a,t}''(\xi_a(t)).\label{-3-tr}\EDE
Now we apply It\^o's formula to the $(\F^a_t) $-adapted SDE. From (\ref{a-SDE-1}), (\ref{xi=W(xi)-tr}), and
(\ref{-3-tr}), we have
$$d\xi^*_{u(a)}(u_a(t))=\til W_{a,t}'(\xi_a(t))\sqrt\kappa dB_a(t)+
\Big(3-\frac\kappa 2\Big)\Big(\til W_{a,t}'(\xi_a(t))X_a(t) dt- \til W_{a,t}''(\xi_a(t))\Big)dt.$$
From (\ref{X(v(t))-mrX(t)}) we have
\BGE d\xi^*_{u(a)}(u_a(t))=\til W_{a,t}'(\xi_a(t))\sqrt\kappa dB_a(t)+\til W_{a,t}'(\xi_a(t))^2
\Big(3-\frac\kappa 2\Big) X^*_{u(a)}(u_a(t))dt,\label{SDE-xi*-v}\EDE
where $X^*_b(t)=X^*(b+t)$ for $b\in(-\infty,T^*)$ and $t\in[0,T^*-b)$.  Now we apply some time-change.
Recall that $u_a$ is continuously increasing, and maps $[0,T_a)$ onto $[0,T^*-u(a))$. So its inverse map, say
$v_a$ is defined on $[0,T^*-u(a))$. We extend $v_a$ to be defined on $[0,\infty)$ such that if $t>T^*-u(a)$ then
$v_a(t)=T_a$. Since $(u_a)$ is  $(\F^a_t) $-adapted, so for every $t\in[0,\infty)$, $v_a(t)$ is an  $(\F^a_t) $-stopping
time. Let $\F^{a,v}_t=\F^a_{v_a(t)}$, $0\le t<\infty$. Then we have a new filtration $(\F^{a,v}_t)_{0\le t<\infty}$.
From (\ref{v'-tr}) and (\ref{SDE-xi*-v}) we see that there is a stopped $(\F^{a,v}_t)$-Brownian motion $B_{a,v}(t)$,
$0\le t<T^*-u(a)$, such that
$\xi^*_{u(a)}(t)$  satisfies the $(\F^{a,v}_t)$-adapted SDE:
\BGE d\xi^*_{u(a)}(t)= \sqrt\kappa dB_{a,v}(t)+
\Big(3-\frac\kappa 2\Big) X^*_{u(a)}(t)dt,\quad 0\le t<T^*-u(a).\label{SDE-xi*}\EDE
Using Proposition \ref{improper}, we may define
\BGE B^*(t)=\xi^*(t)-\Big(3-\frac\kappa 2\Big)\int_{-\infty}^t X^*(s)ds, \quad -\infty<t<T^*.\label{def-B*}\EDE
From (\ref{SDE-xi*}) we have
\BGE \sqrt\kappa B_{a,v}(t)=B^*(u(a)+t)-B^*(u(a)),\quad 0\le t<T^*-u(a).\label{B-B/kappa}\EDE
From (\ref{xi=W(xi)}) and  (\ref{X(v(t))-mrX(t)}) we know that $(e^i(\xi^*(u(t))))$ and $(X^*(u(t)))$
are both $(\F_t)$-adapted. So $(e^i(B^*(u(t))))$ is also  $(\F_t)$-adapted. Especially, $e^i(B^*(t))$, $-\infty<t\le u(a)$,
are $\F_a$-measurable.  Since $\F_a=\F^a_0=\F^{a,v}_0$, so from (\ref{B-B/kappa}),
$B_{a,v}(t)=(B^*(u(a)+t)-B^*(u(a)))/\sqrt\kappa$, $0\le t<T^*-u(a)$, is a stopped Brownian motion independent of
$e^i(B^*(t))$, $-\infty<t\le u(a)$.


Recall that in the above argument, we need that $a$ is a finite $(\F_t)$-stopping time such that $T>a$ always holds.
Let $R=\dist(0,\C\sem(D^*\sem\{z^*_e,V(\infty)\}))$. From Theorem
\ref{existance-and-uniqueness-Theorem-Interior} (ii), for any $r\in
(0,R)$, there is $t_r<T$ such that $K_{t_r}\not\subset
V^{-1}(\B(0;r))$, so $K^*_{u(t_r)}=V(K_{t_r})\not\subset\B(0;r)$.
Thus, $T^*>u(t_r)=\ccap(K^*_{u(t_r)})\ge \ln(r/4)$. So $T^*\ge
\ln(R/4)$. Fix any deterministic number $b\in(-\infty,\ln(R/4))$. Then $T^*>b$ always holds.
Let $a=u^{-1}(b)$. Then $a$ is a finite stopping time such that $T>a$ always holds, and
$u(a)=b$ is a deterministic number. From the last paragraph, we then conclude that
$(B^*(b+t)-B^*(b))/\sqrt\kappa$, $0\le t<T^*-b$, is a stopped Brownian motion independent of $e^i(B^*(t))$, $-\infty<t\le b$.
Since this holds for any deterministic number $b\in(-\infty,\ln(R/4))$, so we may extend $B^*(t)$ to be defined on
  $\R$ such that $(e^i(B^*(t)))$ has the same distribution as $(e^i(B^{(\kappa)}_\R(t)))$. This means that there is
  an integer valued random variable ${\bf n}$ such that $(B^*(t)-2{\bf n}\pi)$ has the same distribution as
  $(B^{(\kappa)}_\R(t))$. Since $\xi^*(t)$ and $\xi^*(t)-2{\bf n}\pi$ generate the same whole-plane Loewner objects,
  so by replacing $\xi^*(t)$ by $\xi^*(t)-2{\bf n}\pi$, we may assume that $(B^*(t))$ has the same distribution as
  $(B^{(\kappa)}_\R(t))$. From (\ref{def-B*}), $\xi^*(t)$ solves
\BGE \xi^*(t)=B^*(t)+\Big(3-\frac\kappa 2\Big)\int_{-\infty}^t X^*(s)ds, \quad -\infty<t<T^*.\label{integral-1}\EDE
So we can conclude that $K^*_t=V(K_{v(t)})$, $-\infty<t<T^*$, is a stopped LERW process in $D^*$
from $0$ to $z_e^*$. To finish the proof, we need to show that $(-\infty,T^*)$ is a.s.\ the maximal interval of the solution
to (\ref{integral-1}) for the extended function $B^*(t)$, which is now defined on $\R$.

Assume that $(-\infty,T^*)$ is not the maximal interval of the solution. So we have $K^*_{T^*}$, which is an
interior hull in $D^*\sem\{\infty,z_e^*\}$ that contains $K^*_t$ for all $t\in(-\infty,T^*)$. Since $\kappa\le 4$, so
a.s.\ $V(\infty)\not\in K^*_{T^*}$. Excluding a null event, we may assume that $K^*_{T^*}\subset
D^*\sem\{\infty,V(\infty),z_e^*\}$. We can find a Jordan curve
$\sigma^*$ in $\C$ such that $K^*_{T^*}\subset H(\sigma^*)\subset D^*\sem \{\infty, V(\infty),z_e^*\}$.
So $V(K_t)=K^*_{u(t)}\subset H(\sigma^*)$ for $-\infty<t<T$. Let $\sigma=V^{-1}(\sigma^*)$. Then
$\sigma$ is a Jordan curve in $\C$, and $H(\sigma)=V^{-1}(H(\sigma^*))\subset D\sem\{V^{-1}\infty,\infty,z_e\}$.
We have $K_t\subset H(\sigma)$ for $-\infty<t<T$, which contradicts Theorem \ref{existance-and-uniqueness-Theorem-Interior}
(ii). So $(-\infty,T^*)$ is a.s.\ the maximal interval of the solution, and the proof is finished. $\Box$

\vskip 3mm

\no {\bf Remark.} The ideas behind  (\ref{v'}), (\ref{xi=W(xi)}), and (\ref{-3}) first appeared in \cite{LSW1},
which were used there to show that SLE$_6$ satisfies locality property.
From the above  proof we see that for any $\kappa\in(0,4]$, the above theorem still holds if $\xi(t)$
is the solution to (\ref{pre-interior-LERW-eqn}) with $f(t)=B^{(\kappa)}_\R(t)$ and $\lambda=3-\frac\kappa 2$.
If $\kappa>4$, the statement should be modified. We can conclude that after a time-change,
$(V(K_t),-\infty<t<S)$ has the same distribution as $(K^*_t,-\infty<t<S^*)$, where $S\in(-\infty,T]$ is the biggest number
such that $K_t\subset D\sem \{V^{-1}(\infty)\}$ for $t\in(-\infty,S)$, and $S^*\in(-\infty,T^*]$ is the biggest number
such that $K^*_t\subset D^*\sem \{V(\infty)\}$ for $t\in(-\infty,S^*)$.

\subsection{Local martingales} \label{Section-martingale}
Let $D$ be a finitely connected domain, $0\in D$, and $z_e\in
D\sem\{0\}$. Let $p=R_\TT(z_e)$ and $\Om=R_\TT(D)$. For $\xi\in
C((-\infty,T))$, let $L^\xi_t$ (resp.\ $\til L^\xi_t$) and $\psi^\xi_t$
(resp.\ $\til\psi^\xi_t$)
be the inverted whole-plane (resp.\ covering whole-plane)
Loewner hulls and maps driven by $\xi$.
Suppose $\bigcup_{t<T}L^\xi_t\subset\Om\sem\{p\}$.
For each $t\in(-\infty,T)$ and $x\in\R$, let $P^\xi(t,x,\cdot)$ be the
generalized Poisson kernel in $\Om^\xi_t$ with the pole at $e^{ix}$,
normalized by $P^\xi(t,x,\psi^\xi_t(p))=1$, and let $\til
P^\xi(t,x,\cdot)= P^\xi(t,x,\cdot)\circ e^i$. It
is standard to check that both $P^\xi$ and $\til P^\xi$ are $C^{1,2,h}$
differentiable, where ``$h$'' means harmonic.

\begin{Lemma} For any $t\in(-\infty,T)$ and $z\in\til \Om\sem \til L^\xi_t$, we
have  $\til{\cal V}_t(z)=0$, where
$$
\til {\cal V}_t(z)=\pa_1 \til P^\xi (t,\xi(t),\til\psi^\xi_t(z))+2\pa_2 \til
P^\xi(t,\xi(t),\til \psi^\xi_t(z))
X^{\xi}_t+\pa_2^2 \til P^\xi(t,\xi(t),\til
\psi^\xi_t(z))$$$$+2\Ree({\pa_{3,z}\til
P^\xi(t,\xi(t),\til\psi^\xi_t(z))}\cot_2
({\til\psi^\xi_t(z)-\xi(t)})). $$
Here $\pa_1$ and $\pa_2$ are partial derivatives w.r.t.\ the first
two (real) variables, and $\pa_{3,z}=(\pa_{3,x}-i\pa_{3,y})/2$ is
the partial derivative w.r.t.\ the third (complex) variable.
\label{vanish*}
\end{Lemma}
{\bf Proof.}  For simplicity, we assume that $\pa\Om$ is smooth, so every boundary
point of $\Om$ or $\Om^\xi_t$ corresponds to a prime end.
In the general case, we have to work on the conformal closure of $\Om$.
For any $t\in(-\infty,T)$ and $z\in \Om\sem  L^\xi_t$, let
$${\cal V}_t(z)=\pa_1   P^\xi (t,\xi(t), \psi^\xi_t(z))+2\pa_2
P^\xi(t,\xi(t), \psi^\xi_t(z)) X^{\xi}_t+\pa_2^2  P^\xi(t,\xi(t),
\psi^\xi_t(z))$$
\BGE+2\Ree\Big({\pa_{3,z}
P^\xi(t,\xi(t),\psi^\xi_t(z))}\psi^\xi_t(z) \frac{e^{i\xi(t)}+\psi^\xi_t(z)}
{e^{i\xi(t)}-\psi^\xi_t(z)}\Big).\label{V_t}\EDE
It is easy to check that ${\cal V}_t\circ e^i=\til{\cal V}_t$.
For $t\in(-\infty,T)$, $x\in\R$ and $z\in\pa\Om$, since
$\psi^\xi_t(z)\in\pa \Om^\xi_t\sem\TT$, so
$P^\xi(t,x,\psi^\xi_t(z))=0$, which implies that
$\pa_2  P^\xi=\pa_2^2   P^\xi=0$ at $(t,x, \psi^\xi_t(z))$, and
$$\pa_1   P^\xi(t,x, \psi^\xi_t(z))+2\Ree\Big({\pa_{3,z}
P^\xi(t,\xi(t),\psi^\xi_t(z))}\psi^\xi_t(z) \frac{e^{i\xi(t)}+\psi^\xi_t(z)}
{e^{i\xi(t)}-\psi^\xi_t(z)}\Big)=0.$$ Thus, ${\cal V}_t$
vanishes on $\pa\Om $ for $t\in [0,T)$.

Let ${\cal W}_t={\cal
V}_t\circ (\psi^\xi_t)^{-1}$ and $\til{\cal W}_t={\cal W}_t\circ e^i$. Then ${\cal W}_t$ vanishes on
$\pa\Om^\xi_t\sem\TT$ for $t\in (-\infty,T)$, and $\til{\cal W}_t=\til{\cal V}_t\circ (\til\psi^\xi_t)^{-1}$.
Thus, for $t\in(-\infty,T)$ and
$w\in \til\Om^\xi_t$,
$$\til{\cal W}_t(w)= \pa_1  \til P^\xi(t,\xi(t),w)+2\pa_2
\til P^\xi(t,\xi(t),w)X^\xi_t$$\BGE+\pa_2^2
\til P^\xi(t,\xi(t),w)+2\Ree ({\pa_{3,z}\til P^\xi(t,\xi(t),w)} \cot_2(w-\xi(t)) ).\label{W_t}\EDE
Since $\til P^\xi(t,\xi(t),\cdot)$ vanishes on $\R\sem\{\xi(t)+2n\pi:n\in\Z\}$,
and $\cot_2(w-\xi(t))$ is real on $\R\sem\{\xi(t)+2n\pi:n\in\Z\}$, so $\til {\cal W}_t$
vanishes on $\R\sem\{\xi(t)+2n\pi:n\in\Z\}$, which implies that
${\cal W}_t$ vanishes on $\TT\sem\{e^{i\xi(t)}\}$. So ${\cal W}_t$ vanishes on $\pa\Om^\xi_t\sem \{e^{i\xi(t)}\}$.

Since $\til P^\xi(t,x,\cdot)$ has period $2\pi$, and has  simple poles at $x+2n\pi$, $n\in\Z$, so there are
 $c(t,x)\in\R$ and some analytic
function $F(t,x,\cdot)$ defined in some neighborhood of $\R$ such
that in that neighborhood,
$P^\xi(t,x,w)=\Imm(F(t,x,w)+{c(t,x)}\cot_2(w-x))$. Then we have
$$\pa_1 P^\xi(t,\xi(t),w)=\Imm( \pa_1 F(t,\xi(t),w)+\pa_1
c(t,\xi(t))\cot_2(w-\xi(t))).$$
$$\pa_2 \til P^\xi(t,\xi(t),w)=\Imm
\Big(\pa_2 F(t,\xi(t),w)+{\pa_2
c(t,\xi(t))}\cot_2({w-\xi(t)})+\frac{c(t,\xi(t))}{2\sin_2{(w-\xi(t))}^{2}}\Big).$$
$$\pa_2^2 \til P^\xi(t,\xi(t),w)=\Imm \Big(\pa_2^2
F(t,\xi(t),w)+{\pa_2^2
c(t,\xi(t))}\cot_2({w-\xi(t)})$$$$+\frac{2\pa_2c(t,\xi(t))}{2\sin_2{(w-\xi(t))}^{2}}+
\frac{c(t,\xi(t))\cos_2{(w-\xi(t))}}{2\sin_2{(w-\xi(t))}^{3}}\Big).$$
$$2\Ree(\pa_{3,z}P^\xi(t,\xi(t),w)\cot_2({w-\xi(t)}))
=\Imm\Big({2F'(t,\xi(t),w)}\cot_2({w-\xi(t)})$$$$-\frac{c(t,\xi(t))\cos_2{(w-\xi(t))}}{2\sin_2{(w-\xi(t))}^{3}}\Big).$$
From (\ref{W_t}) and the above formulas, $\til {\cal W}_t(w)$ equals the imaginary part of
$$\pa_1 F(t,\xi(t),w)+{\pa_1
c(t,\xi(t))}\cot_2({w-\xi(t)})$$
$$+2\Big(\pa_2 F(t,\xi(t),w)+{\pa_2
c(t,\xi(t))}\cot_2({w-\xi(t)}) +\frac{c(t,\xi(t))}{2\sin_2{(w-\xi(t))}^{2}}\Big)X^\xi_t$$
$$+ \pa_2^2
F(t,\xi(t),w)+{\pa_2^2
c(t,\xi(t))}\cot_2({w-\xi(t)}) +\frac{\pa_2c(t,\xi(t))}{\sin_2{(w-\xi(t))}^{2}}$$$$+
\frac{c(t,\xi(t))\cos_2{(w-\xi(t))}}{2\sin_2{(w-\xi(t))}^{3}} +{2F'(t,\xi(t),w)}\cot_2({w-\xi(t)})
-\frac{c(t,\xi(t))\cos_2{(w-\xi(t))}}{2\sin_2{(w-\xi(t))}^{3}}$$
$$=G_t(w)+{A_1(t)}\cot_2({w-\xi(t)})+\frac{A_2(t)}{\sin_2{(w-\xi(t))}^{2}}$$
for some function $G_t$, which is analytic near $\R$, and real valued functions $A_1(t)$ and $A_2(t)$, where
$A_2(t)=c(t,\xi(t))X^\xi_t+\pa_2 c(t,\xi(t))$.

Since $J^\xi_t=G(\Om^\xi_t,\vphi^\xi_t(p);\cdot)$, so for $x\in\R$,
$\pa_{\bf n}J^\xi_t(e^{ix})$ equals the value at $\vphi^\xi_t(p)$ of the
(usual) Poisson kernel in $\Om^\xi_t$ with the pole at $e^{ix}$. Comparing
the residues of $\pa_{\bf n}J^\xi_t(e^{ix})$ and $P^\xi(t,x,\vphi^\xi_t(p))$
at $e^{ix}$, we conclude that
$$\pa_{\bf n}J^\xi_t(e^{ix})/(-1/\pi)=P^\xi(t,x,\vphi^\xi_t(p))/(2c(t,x))=1/(2c(t,x)).$$
It is clear that $\pa_{\bf n}J^\xi_t(e^{ix})=\pa_y \til J^\xi_t(x)$.
Thus, $c(t,x)\pa_y\til J^\xi_t(x)=-1/(2\pi)$ for any $x\in\R$. Differentiating
this equality w.r.t.\ $x$, we get
 $$0=c(t,\xi(t))\pa_x\pa_y \til J^\xi_t(\xi(t))+\pa_2
c(t,\xi(t))\pa_y \til J^\xi_t(\xi(t))=A_2(t)\pa_y \til J^\xi_t(\xi(t)).$$ Thus,
$A_2(t)$ vanishes. So $\til {\cal W}_t(w)$ equals the imaginary part of some
analytic function plus $ {A_1(t)}\cot_2({w -\xi(t)})$ near $\R$.
Hence, $ {\cal W}_t(w)$ equals the imaginary part of some
analytic function plus $-i {A_1(t)}\frac{e^{i\xi(t)}+w}{e^{i\xi(t)}-w}$ near $\TT$.
Since ${\cal W}_t$ is harmonic in $\Om^\xi_t$, and vanishes at every
prime end of $\Om^\xi_t$ other than $e^{i\xi(t)}$, so ${\cal
W}_t=C(t)P^\xi(t,\xi(t),\cdot)$ for some $C(t)\in\R$. Since
$P^\xi(t,x,\psi^\xi_t(p))=1$ for any $t\in(-\infty,T)$ and $x\in\R$, so from
(\ref{V_t}), we have ${\cal V}_t(p)=0$. Thus, ${\cal W}_t(\psi^\xi_t(p))=0$. So for $t\in(-\infty,T)$, we have
$C(t)=0$, which implies that ${\cal W}_t$ vanishes on $\Om^\xi_t$,
and so $\til{\cal V}_t={\cal W}_t\circ \psi^\xi_t\circ e^i$ vanishes on $\til\Om\sem \til L^\xi_t$. $\Box$

\begin{Theorem} Let $\beta(t)$, $-\infty\le t<T$, be an $\LERW(D;0\to z_e)$ curve.
For each $t\in(-\infty,T)$, let $P_t$ be the generalized Poisson kernel in $D\sem \beta([-\infty,t])$ with
the pole at $\beta(t)$, normalized by $P_t(z_e)=1$. Then for any $z\in
D\sem\{0\}$, $(P_t(z))$ is a continuous local martingale.
\label{SLE2martingale}
\end{Theorem}
{\bf Proof.} We may assume that the driving function
$\xi(t)$, $-\infty<t<T$, is the maximal solution to (\ref{pre-interior-LERW-eqn})
with $f(t)=B^{(2)}_\R(t)$ and $\lambda=2$. Then
 $\bigcup_{t<T} L^\xi_t\subset \Om\sem\{p\}$. Let $P^\xi$ be  defined as at the beginning of
 this subsection.
Then  $P_t\circ R_\TT\circ e^i(z)=\til P^\xi(t,\xi(t),\til\psi^\xi_t(z))$.
 Let $(\F_t)$ be the filtration generated by $(e^i(B^{(2)}_\R(t)))$.  Then $(e^i(\xi(t)))$, $(\psi^\xi_t)$,
$(\Om^\xi_t)$ and $(X^{\xi}_t)$  are all $(\F_t)$-adapted. Let $R=\dist(0;\pa D\cup\{z_e\})$.
Fix a constant $a\in(-\infty,\ln(R/4))$. Then $a$ is always less than $T$. Let $T_a=T-a$ and
$\xi_a(t)=\xi(a+t)-\xi(a)$ for $0\le t<T_a$.
Let $B_a(t)=(B^{(2)}_\R(a+t)-B^{(2)}_\R(a))/\sqrt 2$
for $t\ge 0$. Then $B_a(t)$ is an $(\F_{a+t})_{t\ge 0}$-Brownian
motion, and $\xi_a(t)$ satisfies the $(\F_{a+t})_{t\ge 0}$-adapted
SDE: \BGE d\xi_a(t)=\sqrt 2 dB_a(t)+2 X^{\xi}_{a+t}dt,\quad  0\le t<T_a.\label{SDE-zeta}\EDE

For each $t\in[0,T_a)$ and $x\in\R$, let $Q(t,x,\cdot)$ be the
generalized Poisson kernel in $M_{e^{i\xi(a)}}^{-1}(\Om^\xi_{a+t})$ with the pole at $e^{ix}/{e^{i\xi(a)}}$,
normalized by $Q(t,x,\psi^\xi_{a+t}(p)/{e^{i\xi(a)}})=1$, and let $\til
Q(t,x,\cdot)= Q(t,x,\cdot)\circ e^i$. It is clear that $\til Q(t,x,z)=\til P^\xi(a+t,x+\xi(a),z+\xi(a))$
for $0\le t<T_a$, $x\in\R$ and $z\in A_{e^{i\xi(a)}}^{-1}(\til\Om^\xi_{a+t})$.
Since $e^{i\xi(a)}$ is $\F_a$-measurable, and $\Om^\xi_{a+t}$ is
$\F_{a+t}$-measurable, so $(Q(t,\cdot,\cdot))$ is $(\F_{a+t})_{t\ge 0}$-adapted, and so is
 $(\til Q(t,\cdot,\cdot))$.

For $0\le t<T_a$ and $z\in\til\Om\sem \til L^\xi_{a+t}$, let $\til g_t(z)=\til\psi^\xi_{a+t}(z)-\xi(a)$. Then $(\til g_t)$ is
$(\F_{a+t})_{t\ge 0}$-adapted, and satisfies $\pa_t \til g_t(z)=\cot_2(\til g_t(z)-\xi_a(t))$.
From Lemma \ref{vanish*}, we have that
$$
\pa_1 \til Q (t,\xi_a(t),\til g_t(z))+2\pa_2 \til
Q(t,\xi_a(t),\til g_t(z))
X^{\xi}_{a+t}+\pa_2^2 \til Q(t,\xi_a(t),\til
g_t(z))$$$$+2\Ree({\pa_{3,z}\til
Q(t,\xi_a(t),\til g_t(z))}\cot_2
({\til g_t(z)-\xi_a(t)}))=0.$$
Since $P_{a+t}\circ R_\TT\circ e^i(z)=\til Q(t,\xi_a(t),\til g_t(z))$, so from It\^o's formula, the above formula
and that $\pa_t \til g_t(z)=\cot_2(\til g_t(z)-\xi_a(t))$, we conclude that for any $z\in \til\Om$,
$(P_{a+t}\circ R_\TT\circ e^i(z),0\le t<T_a)$ is a continuous local martingale. Since $R_\TT\circ e^i$ maps
$\til\Om$ onto $D\sem\{0\}$, so for any $z\in D\sem\{0\}$,
$(P_t(z),a\le t<T)$ is a continuous local martingale. Since this holds for any $a\in(-\infty,\ln(R/4))$, so
the proof is completed. $\Box$

\vskip 3mm
\no{\bf Remark.} The similar local martingales first appear in \cite{LSW-2}, which was used to prove the
convergence of LERW to radial SLE$_2$.
For the process in the case $(\kappa,\lambda)\ne (2,2)$,  so far
we do not know any local martingale generated by harmonic functions.

\subsection{Other kinds of targets}

Suppose $D$ is a finitely connected
domain that contains $0$, and $I_e$ is a side arc of $D$.
Then $R_\TT(I_e)$ is a side arc of $\Omega=R_\TT(D)$.
Now we change the definition of $J^\xi_t$ in (\ref{jxit}) by replacing $G(\Om\sem
L^\xi_t,p;\cdot)$  by $H(\Om\sem L^\xi_t,R_\TT(I_e);\cdot)$, which is the harmonic measure of $R_\TT(I_e)$
in $\Om\sem L^\xi_t$, and still let $\til J^\xi_t=J^\xi_t\circ e^i$ and $X^\xi_t=(\pa_x\pa_y/\pa_y)
\til J^\xi_t(\xi(t))$. Let everything else in Section \ref{Def-of-quasi}
be unchanged. Then Theorem \ref{existance-and-uniqueness-Theorem-Interior} still
holds. For the new meaning of $X^\xi_t$, let $\xi\in C((-\infty,T))$ be
the maximal solution to (\ref{pre-interior-LERW-eqn}) with $f=B^{(2)}_\R$ and $\lambda=2$.
Let $K^\xi_t$, $-\infty<t<T$, be the whole-plane Loewner hulls driven by $\xi$.
Then we call the interior Loewner chain $K^\xi_t$, $0\le t<T$, a continuous interior LERW in $D$
from $0$ to $I_e$. Let it be denoted by $\LERW(D;0\to I_e )$. Such
Loewner chain is almost surely generated by a random
simple curve started from $0$, which is called an $\LERW(D;0\to I_e)$ curve. Through conformal maps,
we can then define continuous LERW from any interior point to a side arc. Then
we can prove that this kind of continuous LERW is conformally invariant up to
a time-change.

Let $\beta(t)$, $0\le t<T$, denote an $\LERW(D;0\to I_e)$ curve. For each $t\in[0,T)$, let
$P_t$ be the generalized Poisson kernel in $D\sem \beta([0,t])$ with the pole at $\beta(t)$,
normalized by $\int_{I_e} \pa_{\bf n} P_t(z)ds(z)=1$, where $\bf n$ is the inward unit
normal vector, and $ds$ is the measure of length. Then for any fixed $z\in D$, $(P_t(z))$ is
a continuous local martingale.

\vskip 3mm

\no {\bf Remark} After a time-change, a continuous $\LERW(\D;0\to \TT)$ has
the same distribution as a standard disc SLE$_2$   defined in \cite{Zhan}.

\vskip 3mm

Now let $w_e$ be a prime ends of $D$. Then $R_\TT(w_e)$ is a prime end
of $\Om$. Choose $h$ that maps a neighborhood
$U$ of $R_\TT(w_e)$ in $\ha \Om$ conformally onto a neighborhood $V$ of $0$
in $\lin\HH$ such that $h(R_\TT(w_e))=0$ and $h(U\cap\ha\pa \Om)\subset\R$.
Here $\ha\Om$ and $\ha\pa\Om$ are the conformal closure and conformal boundary,
respectively, of $\Om$ as defined in \cite{LERW}.
 Change the definition of $J^\xi_t$ by replacing $G(\Om\sem
L^\xi_t,p;\cdot)$  by $P(\Om\sem L^\xi_t,R_\TT(w_e),h;\cdot)$ in
(\ref{jxit}), where we use $P(\Om\sem L^\xi_t,R_\TT(w_e),h;\cdot)$ to
denote the generalized Poisson kernel $P$ in $\Om\sem L^\xi_t$ with the pole
at $R_\TT(w_e)$, normalized by $P\circ h^{-1}(z)=-\Imm(1/z)+O(1)$, as $z\to 0$ in $\HH$.
We still let $\til J^\xi_t=J^\xi_t\circ e^i$ and
$X^\xi_t=(\pa_x\pa_y/\pa_y)\til J^\xi_t(\xi(t))$. For the new meaning of $X^\xi_t$, let $\xi\in C((-\infty,T))$ be
the maximal solution to (\ref{pre-interior-LERW-eqn}) with $f=B^{(2)}_\R$ and $\lambda=2$.
Let $K^\xi_t$, $-\infty<t<T$, be the whole-plane Loewner hulls driven by $\xi$.
Then we call the interior Loewner chain $K^\xi_t$, $0\le t<T$, a continuous interior LERW in $D$
from $0$ to $w_e$. Let it be denoted by $\LERW(D;0\to w_e )$. Such
Loewner chain is almost surely generated by a random
simple curve started from $0$, which is called an $\LERW(D;0\to w_e)$ curve. Through conformal maps,
we can then define continuous LERW from any interior point to a prime end. Then
we can prove that this kind of continuous LERW is conformally invariant up to
a time-change.

Let $\beta(t)$, $0\le t<T$, denote an $\LERW(D;0\to w_e)$ curve. Fix $h$ that
maps a neighborhood $U$ of $w_e$ in $\ha D$ conformally into $\HH$ such that $h(w_e)=0$
and $h(U\cap\ha\pa D)\subset\R$. For each $t\in[0,T)$, let
$P_t$ be the generalized Poisson kernel in $D\sem \beta([0,t])$ with the pole at $\beta(t)$,
normalized by $\pa_y (P_t\circ h^{-1})(0)=1$. Then for any fixed $z\in D$, $(P_t(z))$ is
a continuous local martingale.

\section{Existence and Uniqueness} \label{Existence}
\subsection{The radial equation}\label{Section-The-radial-equation}
In this subsection, we will prove Theorem
\ref{Thm-pre-LERW-eqn-radial}. We will use the notation in Section
\ref{section-radial-LERW}, and use $\cot_2(z)$, $\sin_2(z)$, $\coth_2(z)$, $\sinh_2(z)$ and
$\cosh_2(z)$ to denote the functions $\cot(z/2)$, $\sin(z/2)$, $\coth(z/2)$, $\sinh(z/2)$ and
$\cosh(z/2)$, respectively.

\begin{Lemma} Let $\xi\in C([0,T))$. Suppose
$a\in[0,T)$ and $H>0$ satisfy $\cosh_2(H)>e^{a/2}$. Then for any
$z\in\C$ with $\Imm z\ge H$, $\til\psi^\xi_a(z)$ is meaningful, and $\cosh_2(\Imm
\til\psi^\xi_a(z))\ge \cosh_2(H)/e^{a/2}$. \label{radial-height}
\end{Lemma}
{\bf Proof.} Let $h>0$ be the solution of $\cosh_2(h)=\cosh_2(H)/e^{a/2}$.
Suppose $z\in\C$ and $\Imm z\ge H$. Let $b\in(0,a]$ be the maximal number such that
$\til\psi^\xi_t(z)$ exists for $t\in[0,b)$. Let
$h(t)=\Imm\til\psi^\xi_t(z)$ for $t\in[0,b)$. From (\ref{radial-Loewner-covering}) we see that there is some real
valued function $\theta(t)$ such that
$$h'(t)=\Imm\cot_2({\til\psi^\xi_t(z)-\xi(t)})=\Imm\cot_2({\theta(t)+i
h(t)})\ge -\coth_2({h(t)}),$$ which implies that
$\tanh_2(h(t))h'(t)/2\ge -1/2$. So for $t\in[0,b)$,
$$\ln\cosh_2(h(t))-\ln\cosh_2(\Imm z)=\ln\cosh_2(h(t))-\ln\cosh_2(h(0))\ge -t/2.$$
Thus, $\cosh_2(h(t))\ge \cosh_2(\Imm
z)/e^{t/2}\ge\cosh_2(H)/e^{a/2}=\cosh_2(h)$, and so $h(t)\ge h$
for $t\in[0,b)$. Since $h>0$, so $\til\psi^\xi_t(z)$ does not blow
up at $b$. Thus, $b=a$, and $\Imm\til\psi^\xi_a(z)=\lim_{t\to
a^-}h(t)\ge h$. So we have $\cosh_2(\Imm
\til\psi^\xi_a(z))\ge \cosh_2(H)/e^{a/2}$.
$\Box$

\begin{Lemma} Let $a,h>0$ be such that $\cosh_2(h)>e^{a/2}$. There
is $C>0$ such that, for any $\eta,\xi\in C([0,a])$, $b\in[0,a]$, and
$z\in\{\Imm z\ge h\}$,
$\til\psi_b^\eta\circ(\til\psi_b^\xi)^{-1}(z)$ is meaningful, and
\BGE |z-\til\psi_b^\eta\circ(\til\psi_b^\xi)^{-1}(z)|\le
C\|\eta-\xi\|_{0,b}.\label{z-psi-zeta-psi-eta}\EDE
\label{Lemma-1-E-U}
\end{Lemma}
\vskip -4mm {\bf Proof.} Suppose $\eta,\xi\in C([0,a])$,
$b\in[0,a]$, and $\Imm z\ge h$. Since $\Imm\til\psi^\xi_t(w)$
decreases in $t$, so $\Imm(\til\psi_b^\xi)^{-1}(z)\ge \Imm z\ge h$.
From Lemma \ref{radial-height}, we see that for $0\le t\le b$,
$\til\psi_t^\xi\circ(\til\psi_b^\xi)^{-1}(z)$ and
$\til\psi_t^\eta\circ(\til\psi_b^\xi)^{-1}(z)$ are meaningful, and
$\cosh_2^2(\Imm \til\psi_t^\eta\circ(\til\psi_b^\xi)^{-1}(z)), \cosh_2^2(\Imm
\til\psi_t^\xi\circ(\til\psi_b^\xi)^{-1}(z))\ge
\cosh_2^2(h)/e^a$, which implies that
 \BGE \sinh_2^2(\Imm \til\psi_t^\eta\circ(\til\psi_b^\xi)^{-1}(z)),\,\sinh_2^2(\Imm
\til\psi_t^\xi\circ(\til\psi_b^\xi)^{-1}(z))\ge \cosh_2^2(h)/e^a-1.\label{height-below} \EDE
Since $|\cot_2'(z)|=\frac 12 |\sin_2^{-2}(z)|\le \frac 12\sinh_2^{-2}(\Imm z)$,
so if $\Imm z_1,\Imm z_2\ge H>0$, then
\BGE |\cot_2(z_1)-\cot_2(z_2)|\le \frac 12 \sinh_2^{-2}(H) |z_1-z_2|.\label{gradient}\EDE
Let
$$g(t)=|\til\psi^\xi_t\circ(\til\psi_b^\xi)^{-1}(z)-\til\psi^\eta_t\circ(\til\psi_b^\xi)^{-1}(z)|,\quad 0\le t\le b.$$
  From (\ref{radial-Loewner-covering}), (\ref{height-below}), and (\ref{gradient}), we see that for any $0\le t\le b$,
$$g(t)\le \int_0^t\Big|\cot_2({\til\psi^\xi_s\circ(\til\psi_b^\xi)^{-1}(z)-\xi(s)})-
\cot_2({\til\psi^\eta_s\circ(\til\psi_b^\xi)^{-1}(z)-\eta(s)})\Big|dt$$
\BGE \le \int_0^t C_1(g(s)+|\xi(s)-\eta(s)|)dt\le C_1 \int_0^t
(g(s)+\|\xi-\eta\|_{0,b})ds,\label{inequality-1}\EDE
where $$C_1=\frac{1/2}{\cosh_2^2(h)/e^a-1}>0.$$
 Let $C=e^{aC_1}-1$. Solving (\ref{inequality-1}), we get
\begin{equation*} |z-\til\psi_b^\eta\circ(\til\psi_b^\xi)^{-1}(z)|= g(b)\le
(e^{bC_{1}}-1)\|\eta-\xi\|_{0,b}\le C\|\eta-\xi\|_{0,b}.\quad \Box
\end{equation*}

\begin{Lemma} There are $a_0,C>0$ such that, for any $t\in[0,a_0]$, and
$\zeta,\eta\in C([0,t])$, we have $L^\eta_t,L^\zeta_t\subset\Omega\sem\{p\}$, and
$|X^\eta_t -X^\zeta_t| \le
C\|\eta-\zeta\|_{0,t}$.\label{Lemma-2-E-U}
\end{Lemma}
{\bf Proof.} There is $H>0$ such that $\St_H\subset\til\Om\sem\til
p$. Choose $a_0>0$ such that $e^{a_0}<\cosh_2(H)$. Let $h>0$ be such
that $\cosh_2(h) =\cosh_2(H)/e^{a_0/2}$. Then
$\cosh_2(h)^2/e^{a_0}>1$. Fix $t\in[0,a_0]$. Suppose $\zeta,\eta\in
C([0,t])$. From Lemma \ref{radial-height}, for any $z\in\C$ with
$\Imm z\ge H$, $\til\psi^\zeta_t(z)$ and $\til\psi^\eta_t(z)$ are meaningful, and
$\Imm \til\psi^\zeta_t(z),\Imm \til \psi^\eta_t(z)\ge h$. Thus, $\til
L^\zeta_t,\til L^\eta_t\subset\St_H\subset\til\Om\sem\til p$, and
$\St_h\subset \til\Om^\zeta_t\sem\til p^\zeta_t,\til\Om^\eta_t\sem\til p^\eta_t$. So we have
$L^\zeta_t,L^\eta_t\subset\Om\sem\{p\}$.

 Choose
$h_0>h_{0.5}>h_1>h_2\in(0,h)$ such that $\cosh_2(h_2)>e^{a_0/2}$.
Let $C_0>0$ be the $C$ given by Lemma \ref{Lemma-1-E-U} with $a=a_0$
and $h=h_2$. Let $C_*>1$ be the number depending only on $h$ and $h_0$
 such that, if $f$ is positive and harmonic
in $\St_h$, and has period $2\pi$, then for any $x_1,x_2\in\R_{h_0}$, $f(x_1)\le
C_* f(x_2)$. Let \BGE\delta=\min\Big\{\frac{h_{0.5}-h_1}{C_0},
\frac{(h_0-h_{0.5})(h_1-h_2)}{8h_0C_0C_*}\Big\}>0.\label{delta-1}\EDE

Suppose first that $\|\eta-\zeta\|_{0,t}<\delta$.
Let $m=\inf\{\til J^\zeta_t(z):z\in\R_{h_0}\}$, $M=\sup\{\til
J^\zeta_t(z):z\in\R_{h_0}\}$, and $D_\nabla=\sup\{|\nabla \til
J^\zeta_t(z)|: z\in\St_{h_{0.5}}\}$.
Since $\St_h\subset
\til\Om^\zeta_t\sem\til p^\zeta_t$, so $\til J^\zeta_t$ is positive
and harmonic in $\St_h$, and vanishes on $\R$. After a reflection
about $\R$, $\til J^\zeta_t$ is harmonic in $\{|\Imm z| <h\}$, and
$|\til J^\zeta_t|$ is bounded by $M$ on $\{|\Imm z|\le h_0\}$.
Moreover, $\til J^\zeta_t$ has period $2\pi$. Thus, $M\le C_*m$.
From Harnack's
inequality, we have \BGE D_\nabla\le 2M/(h_{0}-h_{0.5}),\label{nabla}\EDE and for any
$x\in\R$, \BGE \pa_y \til J^\zeta_t(x)\ge m/h_0,\quad|\pa_x\pa_y
\til J^\zeta_t(x)|\le 4M/h_0^2,\quad |\pa_x^2\pa_y \til
J^\zeta_t(x)|\le 12M/h_0^3.\label{pa-J-xi}\EDE

For $j=1,2$, let $\rho_j=(\til\psi^\eta_t)^{-1}(\R_{h_j})$. Then
$\rho_1$ and $\rho_2$ lie in $\til\Om\sem\til p\sem\til L^\eta_t$, and
$\rho_2$ disconnects $\rho_1$ from $\til L^\eta_t$. Since $\cosh_2(h_2)>e^{a_0/2}$
and $t\in[0,a_0]$, so from Lemma \ref{Lemma-1-E-U}, for any
$z\in\C$ with $\Imm z\ge h_2$,
$\til\psi^\zeta_t\circ(\til\psi^\eta_t)^{-1}(z)$ is meaningful, so
$\rho_1$ and $\rho_2$ lie in $\HH\sem\til L^\zeta_t$, and $\rho_2$
disconnects $\rho_1$ from $\til L^\zeta_t$. Thus, $\rho_1$ and
$\rho_2$ lie in $\til\Om\sem\til p\sem(\til L^\zeta_t\cup\til L^\eta_t)$, and
$\rho_2$ disconnects $\rho_1$ from $\til L^\zeta_t\cup\til
L^\eta_t$. For $\xi\in C([0,t])$, let $G^\xi_t=G(\Om\sem
L^\xi_t,p;\cdot)$ and $\til G^\xi_t=G^\xi_t\circ e^i$. Then $\til
J^\xi_t=\til G^\xi_t\circ(\til\psi^\xi_t)^{-1}$. For $j=1,2$, define
\BGE N_j=\sup_{z\in\R_{h_j}}\{|\til J^\eta_t(z)- \til
J^\zeta_t(z)|\}=\sup_{z\in\R_{h_j}}\{|\til
G^\eta_t\circ(\til\psi^\eta_t)^{-1}(z) -\til G^\zeta_t\circ(\til
\psi^\zeta_t)^{-1}(z)|\};\label{Nj}\EDE \BGE N_j'=\sup_{w\in
\rho_j}\{|\til G^\eta_t(w)-\til
G^\zeta_t(w)|\}=\sup_{z\in\R_{h_j}}\{|\til
G^\eta_t\circ(\til\psi^\eta_t)^{-1}(z) -\til
G^\zeta_t\circ(\til\psi^\eta_t)^{-1}(z)|\}.\label{Nj'}\EDE

Note that $\til J^\eta_t-\til J^\zeta_t$ is harmonic in $\St_h$, and
vanishes on $\R$, since both $\til J^\eta_t$ and $\til J^\zeta_t$
satisfy these properties. Since the probability that a plane
Brownian motion started from a point on $\R_{h_2}$ visits $\R_{h_1}$
before $\R$ is $h_2/h_1$, so \BGE N_2\le
(h_2/h_1)N_1.\label{N2-le-N1}\EDE Since every $z\in\til p$ is a
removable singularity of $\til G^\eta_t-\til G^\zeta_t$, so after an
extension, $\til G^\eta_t-\til G^\zeta_t$ is harmonic in
$\til\Omega\sem (\til L^\zeta_t\cup\til L^\eta_t)$. Since $\rho_1$
and $\rho_2$ lie in $\til\Om\sem(\til L^\zeta_t\cup\til L^\eta_t)$,
and $\rho_2$ disconnects $\rho_1$ from $\til L^\zeta_t\cup\til
L^\eta_t$, so from the maximum principle, we have \BGE N_1'\le
N_2'.\label{N1-le-N2}\EDE

Fix $j\in\{1,2\}$ and $z_0\in \R_{h_j}$. Since $\Imm z_0\ge h_2$, $\cosh_2(h_2)>e^{a_0/2}$, and
$t\in[0,a_0]$, so from Lemma \ref{Lemma-1-E-U}, the choice of $C_0$,  (\ref{delta-1}), and that
$\|\eta-\zeta\|_{0,t}<\delta$, we have
\BGE |z_0-\til\psi^\zeta_t\circ(\til\psi^\eta_t)^{-1}(z_0)|\le
C_0\|\eta-\zeta\|_{0,t}<C_0\delta\le h_{0.5}-h_1.\label{BD-1}\EDE
Thus, $\Imm \til\psi^\zeta_t\circ(\til\psi^\eta_t)^{-1}(z_0)< \Imm
z_0+ h_{0.5}-h_1\le h_{0.5}$. On the other hand, since $
\til\psi^\zeta_t\circ(\til\psi^\eta_t)^{-1}(z_0)\in\HH$, so $\Imm
\til\psi^\zeta_t\circ(\til\psi^\eta_t)^{-1}(z_0)>0$. Thus,
$[z_0,\til\psi^\zeta_t\circ(\til\psi^\eta_t)^{-1}(z_0)]\subset\St_{h_{0.5}}$.
So from (\ref{nabla}) and (\ref{BD-1}),
$$|\til G^\zeta_t\circ(\til \psi^\zeta_t)^{-1}(z_0)-\til G^\zeta_t\circ(\til
\psi^\eta_t)^{-1}(z_0)| =|\til J^\zeta_t(z_0)-\til
J^\zeta_t\circ\til\psi^\zeta_t\circ(\til\psi^\eta_t)^{-1}(z_0)|$$
\BGE \le \sup_{z\in\St_{h_{0.5}}}\{|\nabla\til J^\zeta_t(z)|\}\cdot
|z_0-\til\psi^\zeta_t\circ(\til\psi^\eta_t)^{-1}(z_0)|\le D_\nabla
C_0\|\eta-\zeta\|_{0,t}\le
\frac{2MC_0\|\eta-\zeta\|_{0,t}}{h_0-h_{0.5}}.\label{BD-2}\EDE Let
\BGE \Delta=\frac{2MC_0\|\eta-\zeta\|_{0,t}}{h_0-h_{0.5}}.\label{Delta}\EDE Then from
(\ref{Nj}), (\ref{Nj'}), (\ref{BD-2}), and (\ref{Delta}), we have \BGE
|N_j'-N_j|\le\sup_{z\in \R_{h_j}}\{|\til G^\zeta_t\circ(\til
\psi^\zeta_t)^{-1}(z)-\til G^\zeta_t\circ(\til \psi^\eta_t)^{-1}(z)|
\}\le \Delta,\quad j=1,2.\label{Nj'-Nj}\EDE

 From (\ref{N2-le-N1}), (\ref{N1-le-N2}), and (\ref{Nj'-Nj}), we have
$$N_1\le N_1'+\Delta\le
N_2'+\Delta\le N_2+2\Delta\le (h_2/h_1)N_1+2\Delta.$$ Thus, \BGE N_1\le
2h_1\Delta/(h_1-h_2).\label{N_1le}\EDE From Harnack's inequality, for any $x\in\R$,
\BGE |\pa_y\til J^\eta_t(x)-\pa_y\til J^\zeta_t(x)|\le N_1/h_1\le
2\Delta/(h_1-h_2);\label{pa-J-eta-pa-J-xi-1}\EDE \BGE|\pa_x\pa_y\til
J^\eta_t(x)-\pa_x\pa_y\til J^\zeta_t(x)|\le 4N_1/h_1^2\le
8\Delta/(h_1(h_1-h_2))\label{pa-J-eta-pa-J-xi-2}.\EDE
From (\ref{delta-1}), (\ref{pa-J-xi}), (\ref{Delta}), (\ref{pa-J-eta-pa-J-xi-1}), $M\le C_*m$,
and that $\|\eta-\zeta\|_{0,t}<\delta$, for any $x\in\R$, \BGE
\pa_y \til J^\eta_t(x)\ge \frac {m}{h_0}-\frac{2\Delta}{h_1-h_2}>\frac
{m}{h_0}-\frac{4MC_0\delta}{(h_0-h_{0.5})(h_1-h_2)}\ge \frac
{m}{2h_0}.\label{lower bound}\EDE
From 
(\ref{pa-J-eta-pa-J-xi-2}) and (\ref{lower bound}),
we have \BGE |\pa_x\pa_y\til J^\eta_t(\eta(t))/\pa_y\til
J^\eta_t(\eta(t))-\pa_x\pa_y\til J^\zeta_t(\eta(t))/\pa_y\til
J^\eta_t(\eta(t))|\le
\frac{16 h_0\Delta}{m h_1(h_1-h_2)}.\label{pax-pay-J-1}\EDE
 From (\ref{pa-J-xi}), (\ref{pa-J-eta-pa-J-xi-1}) and (\ref{lower bound}), we have \BGE |\pa_x\pa_y\til
J^\zeta_t(\eta(t))/\pa_y\til J^\eta_t(\eta(t))-\pa_x\pa_y\til
J^\zeta_t(\eta(t))/\pa_y\til J^\zeta_t(\eta(t))|\le \frac{16M\Delta}{m^2(h_1-h_2)}.\label{X(eta)-X(eta)}\EDE
From (\ref{pa-J-xi}) and that $M\le C_*m$, for any $x\in\R$,
$$|\pa_x(\pa_x\pa_y/\pa_y)\til J^\zeta_t(x)|\le |(\pa_x^2\pa_y/\pa_y)\til J^\zeta_t(x)|
+|(\pa_x\pa_y/\pa_y)\til J^\zeta_t(x)|^2$$ $$\le
({12M/h_0^3})/({m/h_0})+({4M/h_0^2})^2/({m/h_0})^2\le
28C_*^2/h_0^2.$$ Thus, \BGE |(\pa_x\pa_y/\pa_y)\til
J^\zeta_t(\eta(t))-(\pa_x\pa_y/\pa_y)\til J^\zeta_t(\zeta(t))|\le
28C_*^2/h_0^2\|\eta-\zeta\|_{0,t}.\label{X(xi)-X(eta)}\EDE From (\ref{Delta}),
(\ref{pax-pay-J-1})-(\ref{X(xi)-X(eta)}), and that $M\le C_* m$, we get
$|X^\eta_t-X^\zeta_t|\le C\|\eta-\zeta\|_{0,t}$, where $$
C:=\frac{32C_*C_0h_0/h_1}{(h_0-h_{0.5})(h_1-h_2)}+\frac{32C_*^2C_0}{(h_0-h_{0.5})(h_1-h_2)}+\frac{28C_*^2}{h_0^2}>0.
$$

In the above argument, we assumed that $\|\eta-\zeta\|_{0,t}<\del$. In the general case,
we can find $n\in\N$ such that $\|\eta-\zeta\|_{0,t}/n<\del$.
Let $\zeta_k=\zeta+k(\eta-\zeta)/n$, $0\le k\le n$. Then
$$|X^\eta_t-X^\zeta_t|\le \sum_{k=1}^n
|X^{\zeta_{k}}_t-X^{\zeta_{k-1}}_t|\le \sum_{k=1}^n
C\|\zeta_{k}-\zeta_{k-1}\|_{0,t}= C\|\eta-\zeta\|_{0,t}.\quad \Box$$ %

\no{\bf Proof of Theorem \ref{Thm-pre-LERW-eqn-radial}.} (i) Let
$a_0,C>0$ be given by Lemma \ref{Lemma-2-E-U}. Use the method of
 Picard iteration to define a sequence of functions $(\xi_n(t))$ in $C([0,a_0])$ such that
 $\xi_0(t)=f(t)$, $0\le t\le a_0$, and for $n\in\N$, \BGE \xi_{n}(t)=f(t)+\lambda\int_0^t
X^{\xi_{n-1}}_sds,\quad 0\le t\le a_0.\label{induction-radial}\EDE Then for $a\in[0,a_0]$, if
 $0\le t\le a$, then
$$|\xi_{n+1}(t)-\xi_n(t)|\le
|\lambda|\int_0^t|X^{\xi_n}_s-X^{\xi_{n-1}}_s|ds\le C|\lambda| t
\|\xi_{n}-\xi_{n-1}\|_{0,a}.$$ Thus, $\|\xi_{n+1}-\xi_n\|_{0,a}\le C
|\lambda| a \|\xi_{n}-\xi_{n-1}\|_{0,a}$. Choose $a\in(0,a_0)$ such
that $C|\lambda|a<1/2$. Then $(\xi_n)$ is a Cauchy sequence w.r.t.\
$\|\cdot\|_{0,a}$. Let $\xi\in C([0,a])$ be the limit of this
sequence. Let $n\to\infty$ in (\ref{induction-radial}), then $\xi$
solves (\ref{pre-LERW-eqn-radial}) for $0\le t\le a$. \vskip 3mm

\no (ii) Suppose for $j=1,2$, $\xi_j$  solves
(\ref{pre-LERW-eqn-radial}) for $0\le t<T_j$ for some $T_j>0$.
Choose $S\in(0,T_1\wedge T_2\wedge a_0)$ such that $C|\lambda|S\le 1/2$. Then
$$\|\xi_1-\xi_2\|_{0,S}\le
C|\lambda|S\|\xi_1-\xi_2\|_{0,S}
\le\|\xi_1-\xi_2\|_{0,S}/2,$$ which implies
$\|\xi_1-\xi_2\|_{0,S}=0$. Thus, $\xi_1(t)=\xi_2(t)$ for $0\le t\le
S$. $\Box$

\subsection{The whole-plane equation} \label{section-proof-whole-plane}
In this section, we will prove Proposition \ref{improper} and
Theorem \ref{existance-and-uniqueness-Theorem-Interior}. We use the
notation in Section \ref{Def-of-quasi}. Let $R=\dist(0,\pa
D\cup\{z_e\})>0$ throughout this subsection.

\begin{Lemma} Suppose $t<\ln(R)-\ln(1+C_{\cal H})$. Let
$h=\ln(R/e^t-C_{\cal H})>0$. Then for any $\xi\in C((-\infty,t])$,
we have $K^\xi_t\subset D\sem\{z_e\}$ and $\St_h\subset
\til\Om^\xi_t\sem\til p^\xi_t$. \label{strip-contains}
\end{Lemma}
{\bf Proof.} Suppose $\xi\in C((-\infty,t])$. From (\ref{est-2}), if $1<|z|<e^h$, then
$|(\vphi^\xi_t)^{-1}(z)-e^tz|\le C_{\cal H}e^t$, and so
$|(\vphi^\xi_t)^{-1}(z)|< e^t(|z|+C_{\cal H})<R$. Thus,
$(\vphi^\xi_t)^{-1}(\{1<|z|<e^h\})\subset \{|z|<R\}\subset D\sem\{z_e\}$,
which implies that $\{1<|z|<e^h\}\subset \vphi^\xi_t(D\sem K^\xi_t\sem\{z_e\})$,
and so $\St_h\subset \til\Om^\xi_t\sem\til p^\xi_t$.
Since $K^\xi_t$ is
surrounded by $(\vphi^\xi_t)^{-1}(\{1<|z|<e^h\})$, so
$K^\xi_t\subset\{|z|<R\}\subset D\sem\{z_e\}$. $\Box$

\begin{Lemma} There are non-increasing functions $S_1$ and $S_2$
defined on $(0,\infty)$ with $S_j(h)=O(he^{-h})$ as $h\to \infty$, $j=1,2$, such that
for any $h>0$,  if $J(z)$ is positive and harmonic in $\St_h$, vanishes on
$\R$, and has period $2\pi$, then for $j=1,2$, \BGE
|(\pa_x^j\pa_y/\pa_y)J(x)|\le S_j(h),\quad\mbox{for any }
x\in\R.\label{HH}\EDE
\label{|X|'}
\end{Lemma}
\vskip -3mm
{\bf Proof.} There is a positive measure $\mu$ on $[0,2\pi)$ such
that $J(z)=\int \SA_h(z-x)d\mu(x)$,  for any $z\in\St_h$, where
\BGE\SA_h(z):=\Imm\Big(\frac{ z}h+\frac 1i\mbox{ P.V.}\sum_{
{n\in\Z, 2\nmid n}} \frac{e^{nh}+e^{iz}}
{e^{nh}-e^{iz}}\Big).\label{SA*}\EDE Such $\SA_h$ is positive and
harmonic on $\St_h$, has period $2\pi$, and vanishes on
$\R\cup\R_h\sem\{2m\pi+hi:m\in\Z\}$. And $2m\pi+hi$ is a simple pole
of $\SA_h$ for each $m\in\Z$. In fact, $\St_h\circ e^i$ is a Poisson kernel in
$\A_h$ with the pole at $e^{-h}$. Let
$$S_j(h)=\sup\{|(\pa_x^j\pa_y/\pa_y)\SA_h(x)|:x\in\R\},\quad j=1,2.$$
Then for $j=1,2$, $S_j(h)\in(0,\infty)$, and $|\pa_x^j\pa_y\SA_h(x)|\le S_j(h)\pa_y
\SA_h(x)$ for any $x\in\R$. Since $J(z)=\int \SA_h(z-x)d\mu(x)$, so
$|\pa_x^j\pa_y J(x)|\le S_j(h)\pa_y  J(x)$ for any $x\in\R$, which
implies (\ref{HH}). If $h'>h$, applying (\ref{HH}) to $J=\SA_{h'}$,
we find that $S_j(h')\le S_j(h)$. So $S_j(h)$ is non-increasing.

Now for $x\in\R$, $j=0,1,2$,
$$\pa_x^j\pa_y\SA_h(x)=\frac{d^{j+1}}{dz^{j+1}}\Big(\frac{
z}h+\frac 1i\mbox{ P.V.}\sum_{ {n\in\Z, 2\nmid n}}
\frac{e^{nh}+e^{iz}} {e^{nh}-e^{iz}}\Big)\Big|_{z=x}.$$ So $$\pa_y
\SA_h(x)=\frac 1h +\sum_{ {n\in\Z, 2\nmid n}} \frac{2e^{nh}e^{ix}}
{(e^{nh}-e^{ix})^2}=\frac 1h +\sum_{n\in\N, 2\nmid n}\Big(
\frac{2e^{nh}e^{ix}} {(e^{nh}-e^{ix})^2}+ \frac{2e^{-nh}e^{ix}}
{(e^{-nh}-e^{ix})^2}\Big)$$ \BGE =\frac 1h +\Ree \sum_{n\in\N,
2\nmid n} \frac{4e^{nh}e^{ix}} {(e^{nh}-e^{ix})^2}\ge\frac 1h
-\sum_{n\in\N, 2\nmid n} \frac{4e^{nh}} {(e^{n h}-1)^2}\ge\frac 1h
-\sum_{n=1}^\infty \frac{4e^{nh}} {(e^{n h}-1)^2};\label{pay}\EDE
\BGE |\pa_x\pa_y \SA_h(x)|=\Big|\sum_{ {n\in\Z, 2\nmid
n}}\frac{4ie^{nh}e^{ix}(e^{nh}+e^{ix})}{(e^{nh}-e^{ix})^3}\Big|\le\sum_{
{n\in\N, 2\nmid n}}\frac{8e^{2nh}}{(e^{nh}-1)^3}\le
\sum_{n=1}^\infty \frac{8e^{2nh}}{(e^{nh}-1)^3};\label{paxy}\EDE
\BGE |\pa_x^2\pa_y\SA_h(x)|=\Big|\sum_{ {n\in\Z, 2\nmid
n}}\frac{-4e^{nh}e^{ix}(e^{2nh}+4e^{nh}e^{ix}+e^{i2x})}{(e^{nh}-e^{ix})^4}\Big|
\le \sum_{n=1}^\infty \frac {24 e^{3nh}}{(e^{nh}-1)^4}
.\label{paxxy}\EDE Thus, $1/\pa_y\SA_h(x)=O(h)$ and
$\pa_x^j\pa_y\SA_h(x)=O(e^{-h})$ for $j=1,2$, as $h\to\infty$,
uniformly in $x\in\R$. So for $j=1,2$, we have $S_j(h)=O(he^{-h})$
as $h\to\infty$. $\Box$

\vskip 3mm

For $t\le\ln(R)$, let \BGE E_0(t)=e^{t-\ln(R)},\quad
E_1(t)=(\ln(R)-t)e^{t-\ln(R)},\quad E_2(t)= E_0(t)+
E_1(t).\label{Del(t)}\EDE Then $ \lim_{t\to -\infty}E_j(t)=0$, $j=0,1,2$; and
 $\int_{-\infty}^t E_j(s)ds= E_{j+1}(t)$, $j=0,1$.

\begin{Lemma} There are absolute constants $M,C>0$ such that, if $t\le \ln(R)-M$ then
for any $\xi\in C((-\infty,t])$,  $K^\xi_t\subset D\sem\{z_e\}$ and $|X^\xi_t|\le
CE_1(t)$. \label{|X|}\end{Lemma} {\bf Proof.} Since
$S_1(x)=O(x e^{-x})$ as $x\to\infty$, so there are $h_0,C_0>0$ such that, if $x\ge
h_0$ then $S_1(x)\le C_0xe^{-x}$. Let $M=\ln(2C_{\cal H}+e^{h_0})$.
Suppose $t\le \ln(R)-M$ and $\xi\in C((-\infty,t])$. Let
$h=\ln(R/e^t-C_{\cal H})$. Then $h>h_0$ and $\ln(R)-t\ge h\ge
\ln(R/2)-t$. Let $C=2C_0$. Then $S_1(h)\le C_0he^{-h}\le CE_1(t)$.
From Lemma \ref{strip-contains}, we have $K^\xi_t\subset D\sem\{z_e\}$ and
$\St_h\subset \til\Om^\xi_t\sem\til p^\xi_t$. Since
$X^\xi_t=(\pa_x\pa_y/\pa_y)\til J^\xi_t(\xi(t))$, and $\til J^\xi_t$
is positive and harmonic in $\til\Om^\xi_t\sem\til p^\xi_t$,
vanishes on $\R$, and has period $2\pi$, so from Lemma \ref{|X|'},
we have $|X^\xi_t|\le S_1(h)\le CE_1(t)$.   $\Box$

\vskip 3mm

\no{\bf Proof of Proposition \ref{improper}.} It is easy to check that
$X^\xi_t$ is continuous in $t$. So from the above lemma, the improper integral converges. $\Box$

\begin{Lemma} Suppose $\xi\in C((-\infty,t])$, $z\in\C\sem\til L^\xi_t$, and $s\in(-\infty,t]$.
Then \BGE e^s\sinh_2^2({\Imm \til\psi^\xi_s(z)})\ge
e^t\sinh_2^2({\Imm \til\psi^\xi_t(z)});\label{tilpsi_ab}\EDE
\BGE\exp({\Imm z})/4\ge e^t\sinh_2^2({\Imm
\til\psi^\xi_t(z)}).\label{tilpsi_a}\EDE
\end{Lemma}
{\bf Proof.} Let $h(r)=\Imm\til\psi^\xi_r(z)$ for $r\in(-\infty,t]$. From (\ref{covering-whole}), there is
a real valued function $\theta$ on $(-\infty,t]$ such that for $r\in(-\infty,t]$,
$$h'(r)=\Imm\cot_2({\til\psi^\xi_r(z)-\xi(r)})=\Imm\cot_2({\theta(r)+i
h(r)})\le -\tanh_2({h(r)}),$$ which implies that
$\coth_2(h(r))h'(r)\le -1$. So we have
$$2\ln\sinh_2(h(t))-2\ln\sinh_2(h(s))= \int_s^t \coth_2({h(r)})h'(r)dr\le
-(t-s).$$ This immediately implies (\ref{tilpsi_ab}). Now let $t$ be
fixed and let $s\to -\infty$. Since $\til\psi^\xi_s(z)-(z-is)\to 0$,
so $e^s\sinh_2^2({\Imm \til\psi^\xi_s(z)})\to \exp({\Imm z})/4$,
which implies (\ref{tilpsi_a}). $\Box$

\begin{Lemma} Let $h_2>h_3>0$, $s\le t$, and $\zeta,\eta\in C((-\infty,t])$. Let
\BGE C_0=\frac1{2\sinh_2^2(h_3)}\mbox{, }\,A_0=\frac{2(1+C_{\cal
H})e^{s-t}}{\sinh_2^2(h_2)}\mbox{, }\,
\Delta_0=e^{C_0}A_0+(e^{C_0}-1)\|\eta-\zeta\|_{s,t}.\label{A,C_0}\EDE
Assume that
\BGE A_0\le 1,\quad \Delta_0\le h_2-h_3.\label{assume-0}\EDE
Then for any $z\in\C$
with $\Imm z\ge h_2$, $\til\psi^\eta_t\circ(\til\psi^\zeta_t)^{-1}(z)$ is
meaningful, and
$$|\til\psi^\eta_t\circ(\til\psi^\zeta_t)^{-1}(z)-z|<
\Delta_0.$$\label{Lemma-1-E-U*}
\end{Lemma} \vskip -7mm
{\bf Proof.} 
Fix $z\in\C$ with $\Imm z\ge h_2$. From (\ref{tilpsi_a}), we have
$$\exp(\Imm(\til\psi^\zeta_t)^{-1}(z))\ge4 e^t\sinh_2^2(h_2).$$
Then we have  $(1+C_{\cal
H})e^s\exp(-\Imm(\til\psi^\zeta_t)^{-1}(z))\le A_0/2\le 1/2$. Thus, from
 (\ref{tilpsi-near-infty}),
$$|\til\psi^\eta_s\circ(\til\psi^\zeta_t)^{-1}(z)-((\til\psi^\zeta_t)^{-1}(z)-is)|\le
4(1+C_{\cal H})e^s\exp(-\Imm(\til\psi^\zeta_t)^{-1}(z))\le  A_0/2.$$
Similarly,
$|\til\psi^\zeta_s\circ(\til\psi^\zeta_t)^{-1}(z)-((\til\psi^\zeta_t)^{-1}(z)-is)|\le
A_0/2$. Thus, \BGE
|\til\psi^\eta_s\circ(\til\psi^\zeta_t)^{-1}(z)-\til\psi^\zeta_s\circ(\til\psi^\zeta_t)^{-1}(z)|\le
A_0.\label{diff-at-b}\EDE Note that $A_0<\Delta_0\le h_2-h_3$. Let $t_0$ be the
maximal number in $(s,t]$ such that, for $r\in[s,t_0)$,
$\til\psi^\eta_r\circ(\til\psi^\zeta_t)^{-1}(z)$ is meaningful, and
\BGEN g(r):=|\til\psi^\eta_r\circ(\til\psi^\zeta_t)^{-1}(z)-
\til\psi^\zeta_r\circ(\til\psi^\zeta_t)^{-1}(z)|<h_2-h_3.\EDEN
From (\ref{diff-at-b}), $g(s)\le A_0$. Since
$\Imm\til\psi^\zeta_r(w)$ decreases in $r$, so for $r\le t$,
$$\Imm\til\psi^\zeta_r\circ(\til\psi^\zeta_t)^{-1}(z)\ge \Imm\til\psi^\zeta_t\circ(\til\psi^\zeta_t)^{-1}(z)
=\Imm z\ge h_2.$$ Thus, for $r\in[s,t_0)$,
$$\Imm\til\psi^\eta_r\circ(\til\psi^\zeta_t)^{-1}(z)\ge\Imm
\til\psi^\zeta_r\circ(\til\psi^\zeta_t)^{-1}(z)-(h_2-h_3) \ge h_3.$$
So $\til\psi^\eta_r\circ(\til\psi^\zeta_t)^{-1}(z)$ does not blow up
at $r=t_0$, and $\Imm
\til\psi^\eta_{t_0}\circ(\til\psi^\zeta_t)^{-1}(z)\ge h_3$. From
(\ref{tilpsi_ab}), for $r\in[s,t_0]$ and $\xi=\zeta$ or $\eta$,
$$\sinh_2^2(\Imm\til\psi^\xi_r\circ(\til\psi^\zeta_t)^{-1}(z))\ge
e^{t_0-r}\sinh_2^2(\Imm\til\psi^\xi_{t_0}\circ(\til\psi^\zeta_t)^{-1}(z))\ge
e^{t_0-r}\sinh_2^2(h_3).$$ Thus, for any $r\in[s,t_0]$ and
$w\in [\til\psi^\zeta_r\circ(\til\psi^\zeta_t)^{-1}(z)-\zeta(r),\til\psi^\eta_r\circ(\til\psi^\zeta_t)^{-1}(z)-\eta(r)]$,
we have $\sinh_2^2(\Imm w)\ge e^{t_0-r}\sinh_2^2(h_3)$.
Since $| \cot_2'(w)|\le \frac 12\sinh_2^{-2}(\Imm w)$ for
$w\in\HH$, so for $r\in[s,t_0]$,
$$\Big|\cot_2({\til\psi^\eta_r\circ(\til\psi^\zeta_t)^{-1}(z)-\eta(r)})
-\cot_2({\til\psi^\zeta_r\circ(\til\psi^\zeta_t)^{-1}(z)-\zeta(r)})\Big|
\le \frac{g(r)+|\eta(r)-\zeta(r)|}{2e^{t_0-r}\sinh_2(h_3)^2}.$$ From
(\ref{covering-whole}) and the above formula, for $r\in[s,t_0]$,
$$g(r)\le
g(s)+\int_s^r\frac{g(u)+|\eta(u)-\zeta(u)|}{2e^{t_0-u}\sinh_2(h_3)^2}\,du
\le A_0+C_0e^{-t_0}\int_s^r e^{u}(g(u)+\|\eta-\zeta\|_{s,t})du.$$
Solving this inequality, we get
$$g(t_0) \le A_0e^{C_0(1-e^{s-t_0})}+ \|\eta-\zeta\|_{s,t}(e^{C_0(1-e^{s-t_0})}-1)< \Delta_0\le h_2-h_3.$$
From the choice of $t_0$, we have $t_0=t$, and so $\til\psi^\eta_t\circ(\til\psi^\zeta_t)^{-1}(z)$
is meaningful, and
$$|\til\psi^\eta_t\circ(\til\psi^\zeta_t)^{-1}(z)-z|=
|\til\psi^\eta_t\circ(\til\psi^\zeta_t)^{-1}(z)-\til\psi^\zeta_t\circ(\til\psi^\zeta_t)^{-1}(z)|
=g(t)< \Delta_0.\quad\Box$$

\vskip 3mm

Suppose $f$ is positive and harmonic in $\{a-H<\Imm z<a+H\}$ for
some $a\in\R$ and $H>0$, and has period $2\pi$. From Harnack's
inequality, there is $C_*>0$ depending only on $H$ such that \BGE
\sup\{f(z):z\in\R_a\}\le C_*\inf\{f(z):z\in\R_a\}.\label{Mm*}\EDE
Let $S_1(h)$ and $S_2(h)$ be given by Lemma \ref{|X|'}. Let
$S_3(h)=S_2(h)+S_1(h)^2$. Then $S_3(h)$ is non-increasing, and
$S_3(h)= O(he^{-h})$ as $h\to\infty$.

\begin{Lemma} Let $s\le t\in\R$ and $\zeta,\eta\in C((-\infty,t])$.
Let  $h>0$ and $H\in(0,h/8]$. Let $h_\lambda=h-(1+\lambda)H$ for
$\lambda=0,0.5,1,2,3$. Let $C_*$ be given by (\ref{Mm*}). Let $C_0,A_0,\Del_0$ be given by
(\ref{A,C_0}). Suppose $L^\zeta_t\subset \Om\sem \{p\}$,
$\St_h\subset \til \Om^\zeta_t\sem\til p^\zeta_t$, and
 \BGE A_0\le 1,\quad\Delta_0\le H^2/({16
C_*h_0}).\label{assume}\EDE Then $L^\eta_t\subset \Om\sem \{p\}$,
and  \BGE
|X^\eta_t-X^\zeta_t|<{192C_*^2}\Delta_0/H^2+S_3(h)|\eta(t)-\zeta(t)|.
\label{C*}\EDE  \label{Lem-main-est**}
\end{Lemma} \vskip -4mm
{\bf Proof.} This lemma is similar to Lemma \ref{Lemma-2-E-U}. The difference is
that this lemma is about the whole-plane Loewner objects, while Lemma \ref{Lemma-2-E-U}
is about the radial Loewner objects.
Recall that $X^\zeta_t=(\pa_x\pa_y/\pa_y)\til J^\zeta_t(\zeta(t))$,
$\til J^\zeta_t$ is positive and harmonic in $\til\Om^\zeta_t\sem
\til p^\zeta_t$, and vanishes on $\R$. Since $\St_h\subset \til
\Om^\zeta_t\sem\til p^\zeta_t$, so after a reflection, $\til
J^\zeta_t$ is harmonic in $\{|\Imm z| <h\}$. Let $m=\inf\{\til
J^\zeta_t(z):z\in\R_{h_0}\}$, $M=\sup\{\til
J^\zeta_t(z):z\in\R_{h_0}\}$, and $D_\nabla=\sup\{|\nabla \til
J^\zeta_t(z)|: z\in\St_{h_{0.5}}\}$. Since
$\{h_0-H<\Imm z<h_0+H\}\subset \St_h$, and $J^\zeta_t$ has period $2\pi$,
so from (\ref{Mm*}), $M\le C_* m$. From Harnack's inequality, we find that
(\ref{nabla}) and (\ref{pa-J-xi}) also hold here. So we have $D_\nabla\le  4M/H$.

From (\ref{assume}), we have $\Del_0<H=h_2-h_3$. So (\ref{assume-0}) holds. From Lemma
\ref{Lemma-1-E-U*}, we see that for any $z\in\C$ with $\Imm z\ge h_2$,
both $\til\psi^\eta_t\circ(\til\psi^\zeta_t)^{-1}(z)$
and $\til\psi^\zeta_t\circ(\til\psi^\eta_t)^{-1}(z)$ are meaningful, and
$|\til\psi^\eta_t\circ(\til\psi^\zeta_t)^{-1}(z)-z|<\Del_0$. Fix
$w\in\C\sem(\til\Om\sem\til p)$. Let $z=(\til\psi^\zeta_t)(w)\in
\HH\sem (\til \Om^\zeta_t\sem \til p^\zeta_t)$. Since $\St_h\subset \til
\Om^\zeta_t\sem\til p^\zeta_t$, so $\Imm z\ge h\ge
h_2$. Thus, $|\til\psi^\eta_t(w)-z|< \Delta_0$, which implies that
$\Imm\til\psi^\eta_t(w)> \Imm z-\Del_0> h-H=
h_0$. Since this holds for any
$w\in\C\sem(\til\Om\sem\til p)$, so $L^\eta_t\subset
\Om\sem\{p\}$ and $\St_{h_0}\subset \til\Om^\eta_t\sem\til
p^\eta_t$. On the other hand, since $h_0<h$, so
$\St_{h_0}\subset\St_h\subset\til \Om^\zeta_t\sem\til p^\zeta_t$.
Thus, $\til J^\zeta_t$ and $\til J^\eta_t$ are both harmonic in
$\St_{h_0}$.

For $j=1,2$, let $\rho_j=(\til\psi^\eta_t)^{-1}(\R_{h_j})$. Then
$\rho_1$ and $\rho_2$ lie in $\til\Om\sem\til L^\eta_t$, and
$\rho_2$ disconnects $\rho_1$ from $\til L^\eta_t$. Since for any
$z\in\C$ with $\Imm z\ge h_2$,
$\til\psi^\zeta_t\circ(\til\psi^\eta_t)^{-1}(z)$ is meaningful, so
$\rho_1$ and $\rho_2$ lie in $\C\sem\til L^\zeta_t$, and $\rho_2$
disconnects $\rho_1$ from $\til L^\zeta_t$. Thus, $\rho_1$ and
$\rho_2$ lie in $\til\Om\sem(\til L^\zeta_t\cup\til L^\eta_t)$, and
$\rho_2$ disconnects $\rho_1$ from $\til L^\zeta_t\cup\til
L^\eta_t$. For $\xi\in C((-\infty,t])$, let $G^\xi_t=G(\Om\sem
L^\xi_t,p;\cdot)$ and $\til G^\xi_t=G^\xi_t\circ e^i$. Then $\til
J^\xi_t=\til G^\xi_t\circ(\til\psi^\xi_t)^{-1}$. For $j=1,2$, define
$N_j$ and $N_j'$ by (\ref{Nj}) and (\ref{Nj'}). Then the same argument can be used
to derive (\ref{N2-le-N1}) and (\ref{N1-le-N2}).

Fix $j\in\{1,2\}$ and $z_0\in \R_{h_j}$. Since $\Imm z_0\ge h_2$, so
from Lemma \ref{Lemma-1-E-U*},
$|z_0-\til\psi^\zeta_t\circ(\til\psi^\eta_t)^{-1}(z_0)|<\Del_0$.
Thus, $\Imm \til\psi^\zeta_t\circ(\til\psi^\eta_t)^{-1}(z_0)< \Imm
z_0+ \Del_0\le h_1+H/2=h_{0.5}$. On the other hand, we have $
\til\psi^\zeta_t\circ(\til\psi^\eta_t)^{-1}(z_0)\in\HH$, so  the line
segment $[z_0,\til\psi^\zeta_t\circ(\til\psi^\eta_t)^{-1}(z_0)]$
lies in $\St_{h_{0.5}}$. So
$$|\til G^\zeta_t\circ(\til \psi^\zeta_t)^{-1}(z_0)-\til G^\zeta_t\circ(\til
\psi^\eta_t)^{-1}(z_0)| =|\til J^\zeta_t(z_0)-\til
J^\zeta_t\circ\til\psi^\zeta_t\circ(\til\psi^\eta_t)^{-1}(z_0)|$$
$$\le \sup_{z\in\St_{h_{0.5}}}\{|\nabla\til J^\zeta_t(z)|\}\cdot
|z_0-\til\psi^\zeta_t\circ(\til\psi^\eta_t)^{-1}(z_0)|< D_\nabla
\Delta_0\le {4M}\Delta_0/{H} .$$ Let $\Delta=4M\Delta_0/H$.
From (\ref{Nj}), (\ref{Nj'}), and the above formula, we find that
 (\ref{Nj'-Nj}) also holds here, which together with
 (\ref{N2-le-N1}) and (\ref{N1-le-N2}) implies (\ref{N_1le}). Thus,
(\ref{pa-J-eta-pa-J-xi-1}) and (\ref{pa-J-eta-pa-J-xi-2}) both hold here.
From
 (\ref{pa-J-xi}), (\ref{pa-J-eta-pa-J-xi-1}), (\ref{assume}),  $\Delta=4M\Delta_0/H$,
  $M\le C_*m$ and $h_1-h_2=H$,
we find that, for any $x\in\R$,
\BGEN \pa_y \til J^\eta_t(x)\ge \frac m{h_0}-\frac{2\Delta}{h_1-h_2}>\frac
{m}{h_0}-\frac{8M\Delta_0}{H(h_1-h_2)}\ge \frac
{m}{2h_0}. \EDEN This is similar to
(\ref{lower bound}). Then (\ref{pax-pay-J-1}) and (\ref{X(eta)-X(eta)}) both hold here.

Using (\ref{pax-pay-J-1}), (\ref{X(eta)-X(eta)}),  $\Delta=4M\Delta_0/H$, $M\le C_*m$, $h_1-h_2=H$ and $h_0\le 2h_1$, we get
\BGE |\pa_x\pa_y\til
J^\eta_t(\eta(t))/\pa_y\til J^\eta_t(\eta(t))-\pa_x\pa_y\til
J^\zeta_t(\eta(t))/\pa_y\til J^\eta_t(\eta(t))|\le
{128C_*\Delta_0}/{H^2} .\label{pax-pay-J-1*}\EDE
\BGE |\pa_x\pa_y\til J^\zeta_t(\eta(t))/\pa_y\til
J^\eta_t(\eta(t))-\pa_x\pa_y\til J^\zeta_t(\eta(t))/\pa_y\til
J^\zeta_t(\eta(t))|\le
{64C_*^2\Delta_0}/{H^2}.\label{X(eta)-X(eta)*} \EDE From
Lemma \ref{|X|'} and the definition of $S_3(h)$, for any $x\in\R$
$$|\pa_x(\pa_x\pa_y/\pa_y)\til J^\zeta_t(x)|\le |(\pa_x^2\pa_y/\pa_y)\til J^\zeta_t(x)|
+|(\pa_x\pa_y/\pa_y)\til J^\zeta_t(x)|^2\le S_3(h).$$ Thus, \BGE
|(\pa_x\pa_y/\pa_y)\til J^\zeta_t(\eta(t))-(\pa_x\pa_y/\pa_y)\til
J^\zeta_t(\zeta(t))|\le
S_3(h)|\eta(t)-\zeta(t)|.\label{X(xi)-X(eta)*}\EDE Then (\ref{C*})
follows from (\ref{pax-pay-J-1*})-(\ref{X(xi)-X(eta)*}). $\Box$

\begin{Lemma} For $j=0,1,2$, Let $E_j(t)$ be as in (\ref{Del(t)}).
There are absolute constants $M,C\ge 1$ such that the followings hold. \\
(i) For any $s\le
t\le \ln(R)-M$, if $\zeta,\eta\in C((-\infty,t])$ and
$\|\eta-\zeta\|_{s,t}\le 1$, then \BGE |X^\eta_t -X^\zeta_t| \le C
(E_0(s)+E_1(t) \|\eta-\zeta\|_{s,t}).\label{paxpa/pay-paxpay/pay*}\EDE
(ii) For any $t\le \ln(R)-M$ and $\zeta,\eta\in
C((-\infty,t])$,\BGE |X^\eta_t -X^\zeta_t| \le C E_1(t)
\|\eta-\zeta\|_{t}.\label{paxpa/pay-paxpay/pay*'}\EDE\label{Lem-main-est*-1}
\end{Lemma}\vskip -7mm
{\bf Proof.} (i) Let $C_*>0$ be the $C_*$ in (\ref{Mm*}) with
$H=1$. Let \BGE C_1=\max\{20e^3(1+C_{\cal H})\exp(5/(2e^4)),\,
2e^8(\exp(5/(2e^4))-1)\}\ge 1.\label{C_1}\EDE Let $h_*>0$ be such
that, if $h\ge h_*$ then $h/e^h\le 1/(32 C_1C_*)$. Let
$$M=\max\{\ln(e^8+20e^3C_{\cal H}),\,\ln(C_{\cal H}+e^{h_*})\}\ge 1.$$
Suppose $s\le t\le\ln(R)-M$, $\zeta,\eta\in C((-\infty,t])$, and
$\|\eta-\zeta\|_{s,t}\le 1$. Let $h=\ln(R/e^t-C_{\cal H})$. It is straightforward  to check that
$h\ge\max\{8,h_*, \ln(R/2)-t\}$, and $\ln(R/2)-t\ge 1$. Since $M>\ln(1+C_{\cal H})$, from Lemma
\ref{strip-contains}, we have $K^\zeta_t,K^\eta_t\subset D\sem\{z_e\}$ and $\St_h\subset\til\Om^\zeta_t\sem\til
p^\zeta_t, \til\Om^\eta_t\sem \til p^\eta_t$.

Let $H=1$. Then $H\in (0,h/8]$. Let $h_\lambda=h-(1+\lambda)H$ for
$\lambda=0,0.5,1,2,3$. Then all $h_\lambda\ge 4$. It is easy to
check that $\sinh_2^2(x)\ge e^x/5$ if $x\ge 4$. Let $C_0,A_0,\Del_0$
be given by (\ref{A,C_0}). Then \BGE A_0\le \frac{2(1+C_{\cal
H})e^{s-t}}{e^{h_2}/5}=\frac{2(1+C_{\cal H})e^{s-t}}{e^{h-3}/5} \le
\frac{2(1+C_{\cal H})e^{s-t}}{e^{-t}R/(10e^3)}=20e^3(1+C_{\cal
H})E_0(s).\label{A_0}\EDE Since $s\le t\le \ln(R)-M$, so $E_0(s)\le e^{-M}
\le 1/(e^8+20e^3C_{\cal H})$. Thus, $A_0\le 1$.
Since $C_0\le 5/(2e^{h_3})$, $h_3=h-4\ge 4$, and $h\ge \ln(R/2)-t$,  so $C_0\le 5/(2e^4)$ and $C_0\le
5e^4 E_0(t)$. Thus, \BGE e^{C_0}\le \exp(5/(2e^4)),\quad
e^{C_0}-1\le \frac{ \exp(5/(2e^4))-1}{5/(2e^4)}\cdot {5e^4}
E_0(t),\label{C_0}\EDE where the second inequality follows from that
$(e^x-1)/x$ is increasing on $(0,\infty)$. Then from (\ref{A,C_0}) and
(\ref{C_1}$\sim$\ref{C_0}), we have \BGE \Del_0\le
C_1(E_0(s)+E_0(t)\|\eta-\zeta\|_{s,t}).\label{Del_0}\EDE Since $h\ge
h_*$, so $h/e^h\le 1/(32 C_1C_*)$. Since $h\le \ln(R/e^t)$, so
$E_0(t)\le 1/e^h$. From $\|\eta-\zeta\|_{s,t}\le 1$, we have
$$\Del_0\le 2C_1E_0(t)\le 2C_1/e^h\le 1/(16C_*h)\le
1/(16C_*h_0).$$ Hence (\ref{assume}) holds. From Lemma
\ref{Lem-main-est**}, we have \BGE |X^\eta_t-X^\zeta_t|\le
192C_*^2\Del_0+S_3(h)|\eta(t)-\zeta(t)|.\label{X-X}\EDE Since
$S_3(x)$ is non-increasing, and $S_3(x)=O(xe^{-x})$ as $x\to\infty$,
so there is an absolute constant $C_S>0$ such that $S_3(x)\le
C_Sxe^{-x}$ for any $x\ge 1$. Since $h\ge \ln(R/2)-t\ge 1$, so \BGE
S_3(h)\le S_3(\ln(R)-\ln(2)-t)\le C_S E_1(\ln(2)+t) \le 2C_S
E_1(t).\label{S_3}\EDE Since
$\ln(R/2)-t\ge 1$, so $E_0(t)\le E_1(t)$.
Let $C=128C_*^2C_1+2C_S\ge 1$. Then (\ref{paxpa/pay-paxpay/pay*}) follows from
(\ref{Del_0})-(\ref{S_3}).

\vskip 3mm

\no (ii) If $\|\eta-\zeta\|_t\le 1$, then
(\ref{paxpa/pay-paxpay/pay*'}) follows from
(\ref{paxpa/pay-paxpay/pay*}) by letting $s\to-\infty$. If
$\|\eta-\zeta\|_t<\infty$, then there is $n\in\N$ such that
$\|\eta-\zeta\|_t<n$. Let $\zeta_k=\zeta+(\eta-\zeta)k/n$, $k=0,1,\dots,n$.
Then $\|\zeta_{k-1}-\zeta_k\|_t<1$ for each $k$, and
$\|\eta-\zeta\|_t=\sum_{k=1}^n\|\zeta_{k-1}-\zeta_k\|_t$. So
(\ref{paxpa/pay-paxpay/pay*'}) follows from the result in the case
$\|\eta-\zeta\|_t\le 1$. If $\|\eta-\zeta\|_t=\infty$,
(\ref{paxpa/pay-paxpay/pay*'}) always hods. $\Box$

\vskip 3mm

\no {\bf Proof of Theorem
\ref{existance-and-uniqueness-Theorem-Interior}.} Let $M,C$ be given
by Lemma \ref{Lem-main-est*-1}. Let $a_0\le\ln(R)-M$ be such that $
C|\lambda| E_2(a_0)\le 1/2$. Define a sequence of functions
$(\xi_n)$ in $C((-\infty,a_0])$ inductively such that, for any $t\le
a_0$ and $n\in\N$, $\xi_0(t)=f(t)$ and \BGE
\xi_n(t)=f(t)+\lambda\int_{-\infty}^t
X^{\xi_{n-1}}_sds.\label{induction-interior}\EDE From Proposition
\ref{improper}, the above improper integrals converge, and
$\|\xi_1-\xi_0\|_a<\infty$. From Lemma \ref{Lem-main-est*-1}, for
$t\le a_0$, $|X^{\xi_{n+1}}_t-X^{\xi_n}_t|\le C E_1(t)$. So from
(\ref{induction-interior}), for any $t\le a_0$,
$$|\xi_{n+1}(t)-\xi_n(t)|\le
C|\lambda|\int_{-\infty}^t E_1(s)ds\|\xi_n-\xi_{n-1}\|_t\le
\|\xi_n-\xi_{n-1}\|_{a_0}/2.$$ Thus, $(\xi_n)$ is a Cauchy sequence
w.r.t.\ $\|\cdot\|_{a_0}$. Let $\xi_\infty$ be the limit. Then
$\xi_\infty$ solves (\ref{pre-interior-LERW-eqn}) for
$t\in(-\infty,a_0]$.

Let ${\cal S}$ be the set of all couples $(\xi,T)$ such that
 $\xi$ solves (\ref{pre-interior-LERW-eqn}) for
$t\in(-\infty,T]$. We have proved that ${\cal S}$ is nonempty.
Suppose $(\xi,T_0)\in\cal S$. Let $\mr\Omega=\Om^\xi_{T_0}$ and
$\mr p=p^\xi_{T_0}\in\mr \Om$. For $\mr \xi\in C([0,S))$ for some $S>0$,
let $\mr L^{\mr \xi}_t$ and
$\mr\psi^{\mr \xi}_t$ denote the {\it radial} Loewner
hulls and maps driven by $\mr \xi$. If $\mr L^{\mr
\xi}_t\subset\mr\Omega\sem\{\mr p\}$, let $\mr J^{\mr
\xi}_t=G(\mr\Omega\sem L^{\mr\xi}_t,\mr p;\cdot)\circ
(\mr\psi^{\mr\xi}_t)^{-1}$, and $\mr X^{\mr
\xi}_t=(\pa_x\pa_y/\pa_y)(\mr J^\xi_t\circ e^i)(\mr\xi(t))$. From
Theorem \ref{Thm-pre-LERW-eqn-radial} (i), the solution to
\BGE\mr\xi(t)=\xi(T_0)+f(T_0+t)-f(T_0)+\lambda\int_0^t \mr
X^{\mr \xi}_sds\label{tilde*}\EDE exists on $[0,b]$ for some
$b>0$. Let $T_e=T_0+b>T_0$. Define $\xi_e(t)=\xi(t)$ for $t\le T_0$
and $\xi_e(t)=\mr\xi(t-T_0)$ for $t\in[T_0,T_e]$. It is clear
that $\xi_e\in C((-\infty,T_e])$. Since $\xi_e$ agrees with $\xi$ on
$(-\infty,T_0]$, so $\xi_e$ solves (\ref{pre-interior-LERW-eqn}) for
$t\in(-\infty,T_0]$. For $t\in[0,T_e-T_0]$, we have
$\psi^{\xi_e}_{T_0+t}=\mr\psi^{\mr\xi}_t\circ\psi^{\xi}_{T_0}$ and
$L^{\xi_e}_{T_0+t}= L^{\xi}_{T_0}\cup (\psi^{\xi}_{T_0})^{-1}(\mr
L^{\mr\xi}_t)$, where $\psi^{\xi_e}_{T_0+t},\psi^\xi_{T_0}$ and
$L^{\xi_e}_{T_0+t},L^\xi_{T_0}$ are the inverted whole-plane Loewner
maps and hulls, while $\mr\psi^{\mr\xi}_t$ and $\mr L^{\mr\xi}_t$ are the
radial Loewner maps and hulls. Since $\psi^\xi_{T_0}$ maps $p$ to
$\mr p$, and maps $\Om\sem L^{\xi_e}_{T_0+t}$ onto $\mr \Omega\sem
\mr L^{\mr\xi}_t$, so
$$\mr J^{\mr \xi}_t=G(\Om\sem
L^{\xi_e}_{T_0+t},p;\cdot)\circ(\psi^\xi_{T_0})^{-1}\circ(\mr\psi^{\mr\xi}_t)^{-1}=G(\Om\sem
L^{\xi_e}_{T_0+t},p;\cdot)\circ(\psi^{\xi_e}_{T_0+t})^{-1}=
J^{\xi_e}_{T_0+t}.$$ Thus, for $t\in[0,T_e-T_0]$, $\mr X^{\mr
\xi}_t=X^{\xi_e}_{T_0+t}$. Since
$\xi({T_0})=f(T_0)+\lambda\int_{-\infty}^{T_0} X^{\xi_e}_sds$,  so
from (\ref{tilde*}),
$$\xi_e(T_0+t)=\mr\xi(t)=\xi(T_0)+f(T_0+t)-f(T_0)+\lambda\int_{T_0}^{T_0+t}  X^{\xi_e}_sds$$
$$=f(T_0+t)+\lambda\int_{-\infty}^{T_0+t} X^{\xi_e}_sds,\quad 0\le t\le T_e-T_0.$$
Thus, $(\xi_e,T_e)\in\cal S$. So we find that for any
$(\xi,T_0)\in\cal S$, there is $(\xi_e,T_e)\in\cal S$ such that
$T_e>T_0$, and $\xi_e(t)=\xi(t)$ for $t\in(-\infty, T_0]$.

Suppose $(\xi_1,T_1),(\xi_2,T_2)\in{\cal S}$. For $j=1,2$, as $t\to -\infty$,
$\xi_j(t)-f(t)\to 0$, so $\xi_1(t)-\xi_2(t)\to 0$. There is
$T<\min\{a_0,T_1,T_2\}$ such that $\|\xi_1-\xi_2\|_T\le 1$. Then
from the argument of the first paragraph, we have
$\|\xi_1-\xi_2\|_T\le \|\xi_1-\xi_2\|_T/2$. Thus, $\xi_1(t)=\xi_2(t)$
for $-\infty<t\le T$. Let $T_0\le T_1\wedge T_2$ be the maximal such
that $\xi_1(t)=\xi_2(t)$ for $-\infty<t\le T_0$. Suppose
$T_0<T_1\wedge T_2$. Let $\mr\xi_1(t)=\xi_1(T_0+t)$,
$\mr\xi_2(t)=\xi_2(T_0+t)$ for $t\in[0,T_0-T]$. Then
$\mr\xi_1$ and $\mr\xi_2$ both solve equation (\ref{tilde*})
for $t\in[0,T_1\wedge T_2-T_0]$. From Theorem
\ref{Thm-pre-LERW-eqn-radial} (ii), there is $S\in(0,T_1\wedge
T_2-T_0]$ such that $\mr\xi_1(t)=\mr\xi_2(t)$ for $0\le t\le
S$, which implies that $\xi_1(t)=\xi_2(t)$ for $0\le t\le T_0+S$.
This contradicts the maximum property of $T_0$. So
$\xi_1(t)=\xi_2(t)$ for $t\in[0,T_1\wedge T_2]$. Let
$T_f=\sup\{T:(\xi,T)\in{\cal S}\}$. Define $\xi_f$ on
$(-\infty,T_f)$ as follows. For any $t\in(-\infty,T_f)$, choose
$(\xi,T)\in{\cal S}$ such that $t\le T$, and let $\xi_f(t)=\xi(t)$.
Then $\xi_f$ is well defined, and solves
(\ref{pre-interior-LERW-eqn}) for $t\in(-\infty,T_f)$. We also have
the uniqueness of $\xi_f$. There is no solution to
(\ref{pre-interior-LERW-eqn}) on $(-\infty,T_f]$. Otherwise, there
exists some solution on $(-\infty,T_f+\eps]$ for some $\eps>0$,
which contradicts the definition of $T_f$.

\vskip 3mm

(i) Let $M_1,C_1$ and $M_2,C_2$ be the $M,C$ given by Lemma
\ref{|X|} and Lemma \ref{Lem-main-est*-1}, respectively. Let
$C=C_1\vee C_2$ and $M=M_1\vee M_2$. Choose $a_0\le\ln(R)-M$ such
that $C|\lambda| E_2(a_0)<1/2$. Then the solution $\xi_f$ exists on
$(-\infty,a_0]$ for any $f\in C(\R)$.

Fix $ a\in\R$. We now prove that $\{f\in C(\R):T_f>a\}\in\T_a$, and
$f\mapsto \xi_f $ is $(\T_a,\T_a)$-continuous on $\{T_f>a\}$. First
suppose $a\le a_0$. Then $\{f\in C(\R):T_f>a\}=C(\R)\in\T_a$.
Suppose $\xi_{f_0}\in G\in\T_a$. Then there are $b_0\le a$ and
$\eps\in(0,1)$ such that $\B_{b_0,a}(\xi_{f_0},\eps):=\{\xi\in
C(\R):\|\xi-\xi_{f_0}\|_{b_0,a}<\eps\}\subset G$. We may choose
$b\le b_0$ and $\del>0$ such that $2\del+6C|\lambda|E_2(b)<\eps$.
Suppose $f\in C(\R)$ and $\|f-f_0\|_{b,a}<\del$. Then
$$|\xi_f(b)-\xi_{f_0}(b)|\le
|f(b)-f_0(b)|+|\lambda|\int_{-\infty}^b(|X^{\xi_{f_0}}_s|+|X^{\xi_f}_s|)ds$$
$$\le \|f-f_0\|_{b,a}+|\lambda|\int_{-\infty}^b 2CE_1(s)ds =
\|f-f_0\|_{b,a}+2C|\lambda|E_2(b)<\eps.$$ Let $a_1\in(b,a]$ be the
maximal number such that $\|\xi_f-\xi_{f_0}\|_{b,a_1}\le 1$. From Lemma
\ref{|X|} and Lemma \ref{Lem-main-est*-1}, for any
$t\in[b,a_1]$,
$$|\xi_{f}(t)-\xi_{f_0}(t)|\le
|f(t)-f_0(t)|+|\lambda|\int_{-\infty}^b(|X^{\xi_{f}}_s|+|X^{\xi_{f_0}}_s|)ds+
|\lambda|\int_b^t|X^{\xi_{f}}_s-X^{\xi_{f_0}}_s|ds$$
$$\le
|f(t)-f_0(t)|+2C|\lambda|\int_{-\infty}^b E_1(s)ds+
C|\lambda|\int_b^t( E_0(b)+ E_1(s)\|\xi_{f}-\xi_{f_0}\|_{b,s})ds$$
$$\le \|f-f_0\|_{b,a_1}+3C|\lambda| E_2(b)+
C|\lambda| E_2(a_1)\|\xi_{f}-\xi_{f_0}\|_{b,a_1}.$$ Since
$C|\lambda| E_2(a_1)\le C|\lambda| E_2(a)\le 1/2$, so
 $$\|\xi_{f_0}-\xi_f\|_{b,a_1}\le \frac{\|f-f_0\|_{b,a_1}+3C|\lambda|
 E_1(b)}{1-C|\lambda| E_2(a_1)}\le
 2\|f-f_0\|_{b,a_1}+6C|\lambda|<\eps<1.$$
So we have $a_1=a$. From the above formula, we have $\|\xi_f-\xi_{f_0}\|_{b_0,a}\le
\|\xi_f-\xi_{f_0}\|_{b,a}<\eps$. Hence $\xi_f\in
\B_{b_0,a}(\xi_{f_0},\eps)\subset G$ if $\|f-f_0\|_{b,a}<\del$. So
$f\to\xi_f$ is $(\T_a,\T_a)$-continuous.

Now consider the case that $a> a_0$. Let $M_0=\ln(R)-a_0$. Suppose
$f_0\in\{T_f>a_0\}$ and $\xi_{f_0}\in G\in\T_a$.  We may choose
$h>0$ such that $\St_h\subset \til\Om^{\xi_{f_0}}_{a}\sem\til
p^{\xi_{f_0}}_{a}$. Let $H=h/8$ and $h_\lambda=h-(1+\lambda)H$ for
$\lambda=0,0.5,1,2,3$. Let $C_*>0$ be given by (\ref{Mm*}). Recall the definition
of $S_3(h)$ before Lemma \ref{Lem-main-est**}. Let
\BGE \del_h= h{2^{-11}C_*^{-1}\exp(-\frac 12 \sinh_2^{-2}(h/2))};\label{deltah}\EDE \BGE
M_{h}=M_0+\max\left\{0,\,\ln\Big(\frac{2(1+C_{\cal
H})}{\sinh_2^2(h/2)}\Big),\, \ln\Big(\frac{2^{12}(1+C_{\cal
H})C_*}{h\sinh_2^2(h/2)}\Big)+\frac1{2\sinh_2^2(h/2)}\right
\};\label{M_h}\EDE \BGE
C_{h}=3\cdot2^{12}C_*^2\exp(\frac 12\sinh_2^{-2}(h/2))\Big(\frac{2(1+C_{\cal
H})e^{M_0}}{h^2\sinh_2^2(h/2)}+\frac 1{h^2}\Big)+S_3(h).\label{C_h}\EDE There are $b_0\le a_0$ and $\eps\in(0,\del_h)$
such that $\B_{b_0,a}(\xi_{f_0},\eps)\subset G$. Let
\BGE \eps_0=\min\Big\{\delta_h,\frac \eps{\exp(C_h|\lambda|(a-a_0))}\Big\}.\label{eps0}\EDE
There is $b_1\le\min\{b_0,\ln(R)-M_h\}$ such that $E_0(b_1)<\eps_0/5$. From
the last paragraph, there are $b\le b_1$ and $\del\in(0,\eps_0/5)$
such that, if $\|f-f_0\|_{b,a_0}<\del$ then
$\|\xi_f-\xi_{f_0}\|_{b_1,a_0}<\eps_0/5$. Suppose $f\in C(\R)$ and
$\|f-f_0\|_{b,a}<\del$. Since $a>a_0$, so
$\|\xi_f-\xi_{f_0}\|_{b_1,a_0}<\eps_0/5<\del_h$. Let $a_1\in(a_0,a]$
be the maximal number such that $\xi_f$ is defined on $(-\infty,a_1)$ and
$|\xi_f(t)-\xi_{f_0}(t)|< \del_h$ on $[b_1,a_1)$. Fix $t\in[a_0,a_1)$.
Since $t< a$, $\Imm\til\psi^{\xi_{f_0}}_s(z)$ decreases in $s$, and
$\St_h\subset \til\Om^{\xi_{f_0}}_{a}\sem\til p^{\xi_{f_0}}_{a}$, so
$\St_h\subset \til\Om^{\xi_{f_0}}_{t}\sem\til p^{\xi_{f_0}}_{t}$.
Let $C_0,A_0,\Del_0$ be given by (\ref{A,C_0}) with $s=b_1$,
$\zeta=\xi_{f_0}$ and $\eta=\xi_f$. Since $t\ge a_0=\ln(R)-M_0$ and
$h_2\ge h_3=h/2$, so \BGE A_0\le \frac{2(1+C_{\cal
H})e^{M_0}}{\sinh_2^2(h/2)}\,E_0(b_1).\label{A_0'}\EDE Since
$C_0=\frac 12\sinh_2^{-2}(h/2)$, so from (\ref{A,C_0}) and (\ref{A_0'}), we
have \BGE\Del_0\le \exp(\frac 12\sinh_2^{-2}(h/2))\Big(\frac{2(1+C_{\cal
H})e^{M_0}}{\sinh_2^2(h/2)}\,E_0(b_1)
+\|\xi_f-\xi_{f_0}\|_{b_1,t}\Big).\label{Del_0'}\EDE Using
(\ref{deltah})-(\ref{Del_0'}) and the facts that $b_1\le \ln(R)-M_h$,
 $\|\xi_f-\xi_{f_0}\|_{b_1,t}\le\del_h$ and $H=h/8$, one may check that
(\ref{assume}) holds, i.e., $A_0\le 1$ and $\Delta_0\le H^2/(16C_* h_0)$. From Lemma \ref{Lem-main-est**} and
(\ref{Del_0'}), we have \BGE |X^{\xi_f}_{t}-X^{\xi_{f_0}}_{t}|\le
C_{h}(E_0(b_1)+\|\xi_{f}-\xi_{f_0}\|_{b_1,t}),\quad t\in[a_0,a_1).\label{X-X(ii)-2}\EDE
Recall that $b\le b_1\le b_0\le a_0< a_1\le a$.
From (\ref{pre-interior-LERW-eqn}) and (\ref{X-X(ii)-2}), for any $t\in[a_0,a_1)$,
$$|\xi_f(t)-\xi_{f_0}(t)|\le
|f(a_0)-f_0(a_0)|+|f(t)-f_0(t)|+|\xi_f(a_0)-\xi_{f_0}(a_0)|
+|\lambda|\int_{a_0}^t|X^{\xi_f}_s-X^{\xi_{f_0}}_s|ds$$
$$\le
2\|f-f_0\|_{b,a}+\|\xi_f-\xi_{f_0}\|_{b_1,a_0}
+C_h|\lambda|\int_{a_0}^t(E_0(b_1)+\|\xi_{f}-\xi_{f_0}\|_{b_1,s})ds.$$
For $t\in[a_0,a_1)$, let $g(t)=\|\xi_{f}-\xi_{f_0}\|_{b_1,t}$, then
$$g(t)\le 2\|f-f_0\|_{b,a}+\|\xi_f-\xi_{f_0}\|_{b_1,a_0}
+C_h|\lambda|\int_{a_0}^t(E_0(b_1)+g(s))ds.$$ Solving this
inequality using (\ref{eps0}) and that $\|f-f_0\|_{b_1,a_0}\le\|f-f_0\|_{b,a_0}<\delta<\eps_0/5$, and
 $\|\xi_f-\xi_{f_0}\|_{b_1,a_0}<\eps_0/5$,
 we have that for any $t\in[a_0,a_1)$,
$$g(t)\le e^{C_h|\lambda|(t-a_0)}(2\|f-f_0\|_{b,a}+\|\xi_f-\xi_{f_0}\|_{b_1,a_0})
+(e^{C_h|\lambda|(t-a_0)}-1)E_0(b_1)$$
$$< e^{C_h|\lambda|(a-a_0)}(2\eps_0/5+\eps_0/5+\eps_0/5)\le 4\eps/5<\eps.$$
So from (\ref{X-X(ii)-2}) we have
$|X^{\xi_f}_{t}-X^{\xi_{f_0}}_{t}|< C_h(E_0(b_1)+\eps)$ for any
$t\in[a_0,a_1)$. Let
$$S=C_h(E_0(b_1)+\eps)+\sup\{|X^{\xi_{f_0}}_t|:t\in[a_0,a]\}<\infty.$$
Then $|X^{\xi_f}_t|\le S$ for any $t\in[a_0,a_1)$. Since
$\xi_f(t)=f(t)+\lambda\int_{-\infty}^t X^{\xi_f}_s ds$, so
$\lim_{t\to a_1}\xi_f(t)$ exists and is finite. By
defining $\xi_f(a_1)=\lim_{t\to a_1}\xi_f(t)$, we have $\xi_f$ that
solves (\ref{pre-interior-LERW-eqn}) for $-\infty<t\le a_1$. Thus, $T_f>a_1$. Since
$\|\xi_{f}-\xi_{f_0}\|_{b_1,t}=g(t)\le 4\eps/5<\del_h$ for all $t\in[a_0,a_1)$, so
from the definition of $a_1$, we have $a_1=a$. Thus, $T_f>a$ and
$\|\xi_f-\xi_{f_0}\|_{b_1,a}=\lim_{t\to a} g(t)\le 4\eps/5<\eps$.
Thus, $f\in\{T_f>a\}$ and $\xi_f\in\B_{b_0,a}(\xi_{f_0},\eps)\subset
G$ if $\|f-f_0\|_{b,a}<\del$. So $\{T_f>a\}\in\T_a$, and $f\mapsto
\xi_f $ is $(\T_a,\T_a)$-continuous on $\{T_f>a\}$.

Let $f_1,f_2\in C(\R)$. Suppose for some $a\in\R$, $T_{f_1}>a$, that is,
$\xi_{f_1}(t)$ is defined on $(-\infty,a]$, and
$f_1\aequ f_2$.
Then there is $k\in\Z$ such that $f_2(t)=f_1(t)+2k\pi$ for $t\le a$. It
is clear that $\xi(t)=\xi_{f_1}(t)+2k\pi$ solves
(\ref{pre-interior-LERW-eqn}) with $f=f_2$ for $-\infty<t\le a$.
Thus, $T_{f_2}>a$ and $\xi_{f_2}(t)=\xi_{f_1}(t)+2k\pi$ for $t\le a$,
so $\xi_{f_1}\aequ \xi_{f_2}$. From the results of the last
paragraph, we have $\{T_f>a\}\in\T^\TT_a$, and $f\mapsto \xi_f $ is
$(\T^\TT_a,\T^\TT_a)$-continuous on $\{T_f>a\}$

\vskip 3mm

(ii) Suppose $\alpha$ is a Jordan curve such that $\bigcup_{t<T_f}
K^{\xi_f}_t\subset H(\alpha)\subset D\sem\{z_e\}$. Then
$t=\ccap(K^{\xi_f}_t)\le\ccap(H(\alpha))$ for any $t<T_f$, so
$T_f\le \ccap(H(\alpha))<\infty$. We may choose another Jordan curve
$\alpha_0$ such that $H(\alpha)\subset U(\alpha_0)$ and
$H(\alpha_0)\subset D\sem\{z_e\}$. Let
$h=\min\{\ln|\vphi_{H(\alpha)}(z)|:z\in\alpha_0\}>0$. For any $t<T_f$,
since $K^{\xi_f}_t\subset H(\alpha)$, so for any $z\in \alpha_0$,
$|\vphi^{\xi_f}_t(z)|=|\vphi_{K^{\xi_f}_t}(z)|\ge
|\vphi_{H(\alpha)}(z)|\ge e^{h}$. Since $\alpha_0$ disconnects
$K^{\xi_f}_t$ from $\C\sem(D\sem\{z_e\})$, so $\{1<|z|<e^h\}\subset
\vphi^{\xi_f}_t(D\sem\{z_e\}\sem K^{\xi_f}_t)$. Thus,
$\St_h\subset\til\Om^{\xi_f}_t\sem\til p^{\xi_f}_t$ for $t<T_f$. Now $\til J^{\xi_f}_t$
is positive and harmonic in $\St_{h}$, vanishes on $\R$, and has
period $2\pi$, so from Lemma \ref{|X|'}, $|X^{\xi_f}_t|\le S_1(h)$
for $t<T_f$. From (\ref{pre-interior-LERW-eqn}), $\lim_{t\to T_f^-}
\xi_f(t)$ exists and is finite. Define $\xi_f(T_f)= \lim_{t\to T_f^-}
\xi_f(t)$. Then $\xi_f$ solves (\ref{pre-interior-LERW-eqn}) for
$-\infty<t\le T_f$, which is a contradiction. $\Box$

\section{Partition Function} \label{distr}
For $\kappa>0$ and $\lambda\in\R$, let a $(\kappa,\lambda)$-process denote the whole-plane Loewner chain driven by the solution to (\ref{pre-interior-LERW-eqn}) with  $f(t)=B^{(\kappa)}_\R(t)$.
In this section, we will prove that a $(\kappa,\lambda)$-process
 is locally absolutely continuous w.r.t.\ the whole-plane SLE$_\kappa$ processes started from $0$.
By setting $\kappa=\lambda=2$, we conclude that the continuous LERW from an interior point to another interior point is
locally absolutely continuous w.r.t.\ the whole-plane SLE$_2$ process.

 Suppose $D$ is a finitely connected domain, $0,z_e\in D$, and $z_e\ne 0$.
Let $K_t$ and $\beta(t)$, $-\infty\le t<\infty$, be a whole-plane SLE$_\kappa$ hulls and trace from $0$ to $\infty$
with the driving function being  $\xi(t)=B^{(\kappa)}_\R(t)$.
Let $\mu$ be the distribution of $(\xi(t))$. Let $(\F^0_t)$ be the filtration generated by $(e^{i\xi(t)})$.
Let $(\F_t)$ be the completion of $(\F^0_t)$ w.r.t.\ $\mu$.
Let $\psi_t$ and $\til\psi_t$ be the inverted and covering inverted  whole-plane Loewner maps driven by $\xi$.
Let $\vphi_t=\vphi_{K_t}$ and $\phi_t=\phi_{K_t}$. Then $\vphi_t=R_\TT\circ \psi_t\circ R_\TT$
and $\phi_t(z)=e^t\vphi_t(z)$. Let $T\in(-\infty,\infty]$ be the maximal number such that $K_t\subset D\sem\{z_e\}$
for $-\infty<t<T$. Let $J^\xi_t$, $\til J^\xi_t$,  $\Omega^\xi_t$, $\til \Omega^\xi_t$, $p^\xi_t$, and $\til p^\xi_t$,
$-\infty<t<T$, be as in Section \ref{Def-of-quasi}. For simplicity, we omit the superscripts $\xi$ in this section.

Let $R=\dist(0,\pa D\cup\{z_e\})$. Let $T_R=\ln(R)-\ln(1+C_{\cal H})$. Let
$h(t)=\ln(R/e^t-C_{\cal H})>0$ for $t<T_R$. From Lemma \ref{strip-contains} we have
$\St_{h(t)}\subset \til \Om_t\sem \til p_t$ for $t<T_R$.
From Lemma \ref{|X|'}, we conclude that $|(\pa_x^j\pa_y/\pa_y)\til J_t(\xi(t))|\le S_j(h(t))$
for $t<T_R$ and $j=1,2$, where $S_j(h)=O(he^{-h})$ as $h\to\infty$. So for $j=1,2$,
\BGE (\pa_x^j\pa_y/\pa_y)\til J_t(\xi(t))=O(te^t),\quad t\to -\infty.\label{Jj}\EDE

Now we study the behavior of $\pa_y \til J_t(\xi(t))$ as $t\to -\infty$. We have to consider two cases. The first
case is that $D=\ha\C$. Then $\Om_t=\D$ for all $t\in\R$. If $z_e=\infty$ then $p=p_t=0$ for all $t\in\R$.
Thus, $J_t(z)=G(\Om_t,p_t;z)=-\frac 1{2\pi}\ln|z|$, and so $\til J_t(z)=J_t(e^{iz})=\frac 1{2\pi}\Imm z$. So
we have $\pa_y \til J_t(\xi(t))=\frac 1{2\pi}$ for all $t\in\R$. Now suppose that $D=\ha\C$ and $z_e\not\in\{0,\infty\}$.
Recall that
$$p_t=\psi_t(p)=R_\TT\circ\vphi_t\circ R_\TT(p)=R_\TT\circ \vphi_t(z_e)=R_\TT(e^{-t}\phi_{K_t}(z_e)).$$
From (\ref{phi_H(z)-z}) we have $|\phi_{K_t}(z_e)-z_e|\le C_{\cal H}e^t$. Thus, $p_t=O(e^t)$ as $t\to -\infty$. We have
$$\til J_t(z)=J_t(e^{iz})=G(\Om_t,p_t;e^{iz})=G(\D,p_t;e^{iz})=-\frac 1{2\pi}\ln\Big|\frac{e^{iz}-p_t}{\lin p_t e^{iz}-1}\Big|.$$
So  we have
$$\pa_y \til J_t(\xi(t))=\frac{1}{2\pi} \frac{1-|p_t|^2}{|1-p_te^{-i\xi(t)}|^2}=\frac 1{2\pi}+O(e^t),\quad t\to -\infty.$$
Thus, when $D=\ha\C$ we always have
\BGE \pa_y \til J_t(\xi(t))=\frac 1{2\pi}+O(e^t),\quad t\to-\infty \label{limit-J0}\EDE

The second case is that $D\ne \ha\C$. For $t\in(-\infty,T)$, let $G_t(z)=G(D\sem K_t,z_e;z)$
and $G_t^\phi(z)=G_t(\phi_t^{-1}(z))=G(\phi_t(D\sem K_t),\phi_t(z_e);z)$. Since $\Om_t=R_\TT\circ\vphi_t(D\sem K_t)
=R_\TT\circ M_{e^t}^{-1}\circ \phi_t(D\sem K_t)$, $p_t=R_\TT\circ\vphi_t(z_e)=R_\TT\circ M_{e^t}^{-1}\circ \phi_t(z_e)$,
and $J_t(z)=G(\Om_t,p_t;z)$, so $J_t=G_t^\phi\circ M_{e^t}\circ R_\TT$. As $t$ decreases, $D\sem K_t$ increases, so
$G_t$ increases. Let $G_{-\infty}(z)=G(D,z_e;z)$. As $t\to -\infty$, since $K_t\to\{0\}$, so $G_t(z)\to G_{-\infty}(z)$
in $D\sem \{z_e,0\}$. Moreover, since $\diam(K_t)\le 4e^t$, so $G_t(z)-G_{-\infty}(z)=O(1/t)$ as $t\to-\infty$, uniformly on any subset
of $D\sem \{z_e,0\}$ that is bounded away from $0$. Using Harnack's inequality, we conclude that
$\nabla G_t\to \nabla G_{-\infty}$, as $t\to -\infty$, uniformly on any compact subset of $D\sem \{z_e,0\}$.

Let $r=R/2$ and $\delta=R/4$. Let $A=\{r-\delta\le |z|\le r+\delta\}$. Then $A$ is a compact subset of $D\sem \{z_e,0\}$.
So there are constants $T_A\in(-\infty,T)$ and $M_A\in(0,\infty)$ such that $|\nabla G_t|\le M_A$ on $A$ if $t\le T_A$.
From (\ref{phi_H-1(z)-z}) we see that $|\phi_t^{-1}(z)-z|\le C_{\cal H}e^t$ for any $t<\ln|z|$. Let $T_B=T_A\wedge
\ln(\delta/C_{\cal H})$. Suppose $|z|=r$ and $t\le T_B$. Then $|\phi_t^{-1}(z)-z|\le C_{\cal H}e^t\le \delta$. So
$[z,\phi^{-1}(z)]\subset A$. Thus,
$$|G^\phi_t(z)-G_t(z)|=|G_t(\phi_t^{-z}(z))-G_t(z)|\le M_A|\phi_t^{-1}(z)-z|\le M_A C_{\cal H}e^t.$$
Since $G_t-G_{-\infty}=O(1/t)$ as $t\to -\infty$, uniformly on $\{|z|=e^t\}$, so
$G^\phi_t-G_{-\infty}=O(1/t)$ as $t\to -\infty$, uniformly on $\{|z|=e^t\}$. Since $\til J_t=G_t^\phi\circ M_{e^t}
\circ R_\TT\circ e^i$, so as $t\to -\infty$,
\BGE \til J_t(x+i(\ln(r)-t))=G_t^\phi(re^ix)=G_{-\infty}(re^{ix})+O(1/t)\label{J-G}\EDE
uniformly in $x\in\R$.

Fix $t\in(-\infty,T_B]$. Let $h=\ln(r)-t$. Let $\SA_h(z)$ be defined as in (\ref{SA*}). So $\SA_h\circ (e^i)^{-1}$
is a Poisson kernel function in $\{e^{-h}<|z|<1\}$ with the pole at $e^{-h}$. Since
$\til J_t$ is harmonic in $\St_h$, continuous on $\lin{\St_h}$, vanishes on $\R$, and has period
$2\pi$, so for any $z\in\St_h$, we have
$$\til J_t(z)=\frac 1{2\pi} \int_{-\pi}^\pi \til J_t(x+ih) \SA_h(z-x)dx.$$
Thus, for any $x_0\in\R$, $\pa_y \til J_t(x_0)=\frac 1{2\pi} \int_{-\pi}^\pi \til J_t(x+ih) \pa_y\SA_h(x_0-x)dx$.
From (\ref{SA*}) and the computation in the proof of Lemma \ref{|X|'} we have $\pa_y \SA_h(x)=\frac 1h+O(e^{-h})
=\frac{1}h +O(e^t)$ as $h\to\infty$, uniformly in $x\in\R$. From (\ref{J-G}) we have $\til J_t(x+ih)=G_{-\infty}(re^{ix})
+O(1/t)$ as $t\to-\infty$, uniformly in $x\in\R$. Thus, as $t\to -\infty$, we have
$$\pa_y \til J_t(x_0)=\frac 1{2\pi} \int_{-\pi}^\pi G_{-\infty}(re^{ix})\frac{dx}h+O(1/({ht}))=\frac{G_{-\infty}(0)}{\ln(r)-t}+
O(1/t^2)$$
uniformly in $x\in\R$, where the second ``$=$'' holds because $G_{-\infty}$ is harmonic in $\{|z|\le r\}$. So as $t\to-\infty$,
we have
\BGE -t\pa_y\til J_t(\xi(t))=G_{-\infty}(0)+O(1/t)=G(D,z_e;0)+O(1/t).\label{limit-J}\EDE

Next, we study the behavior of $(\pa_t\pa_y/\pa_y)\til J_t(\xi(t))$ as $t\to -\infty$.
For $t\in(-\infty,T)$, we have $J_t\circ \psi_t\circ R_\TT=G_t=G(D\sem K_t,z_e,\cdot)$, which implies that
\BGE \til J_t\circ \til \psi_t\circ R_\R(z)=G_t(e^{iz})=G(D\sem K_t,z_e,e^i(z)).\label{G1}\EDE
Let $P_t$ denote the generalized Poisson kernel in $\Omega_t$ with the pole at $e^{i\xi(t)}$, normalized by
$P_t(z)=\Ree \frac{e^{i\xi(t)}+z}{e^{i\xi(t)}-z}+O(z-e^{i\xi(t)})$ as $z\to e^{i\xi(t)}$. So $P_t\circ \psi_t\circ R_\TT$
is a generalized Poisson kernel in $D\sem K_t$ with the pole at $\beta(t)$. If $\delta>0$ is small, then $K_{t+\delta}\sem K_t$
is contained in a small ball centered at $\beta(t)$. So it is intuitive that
$G_t=G(D\sem K_t,z_e,\cdot)$ is differentiable in $t$, and $-\pa_t G_t(z)$ is
a generalized Poisson kernel in $D\sem K_t$ with the pole at $\beta(t)$. This can be proved by expressing
$G_t(z)-G_{t+\delta}(z)$ as an integral of Poisson kernels in
$D\sem K_{t+\delta}$ with poles on $K_{t+\delta}\sem K_t$. We do not go into details here. Thus, there are $C(t)>0$
such that $-\pa_t G_t=C(t)P_t\circ \psi_t\circ R_\TT$, $ -\infty<t<T$.
Let $\til P_t=P_t\circ e^i$. Then $\til P_t(z)=-\Imm\cot_2(z-\xi(t))+O(z-\xi(t))$ as $z\to \xi(t)$, and
we have
\BGE -\pa_t G_t(e^{iz})=C(t)\til P_t\circ \til \psi_t(\lin z),\quad -\infty<t<T.\label{G2}\EDE
Differentiating (\ref{G1}) w.r.t.\ $t$ and using (\ref{covering-whole-eqn}) and (\ref{G2}), we see that for any
$z\in (e^i)^{-1}(D\sem K_t)$,
$$\pa_t \til J_t(\til\psi_t(\lin z))+\pa_x \til J_t(\til\psi_t(\lin z))\Ree \cot_2(\til\psi_t(\lin z)-\xi(t))$$
$$+\pa_y \til J_t(\til\psi_t(\lin z))\Imm \cot_2(\til\psi_t(\lin z)-\xi(t))=-C(t)\til P_t(\til \psi_t(\lin z)).$$
Since $\til\psi_t\circ R_\R$ maps $(e^i)^{-1}(D\sem K_t)$ onto $\til \Omega_t$, so for any $w\in\til\Omega_t$,
\BGE\pa_t \til J_t(w)+\pa_x \til J_t(w)\Ree\cot_2(w-\xi(t))+\pa_y \til J_t(w)\Imm\cot_2(w-\xi(t))=-C(t)\til P_t(w).\label{C1}\EDE
Suppose in some neighborhood $U$ of $\xi(t)$,  $\til J_t=\Imm \til J^\C_t$ and $\til P_t=\Imm \til P^\C_t$, where $\til J^\C_t$ is analytic
in $U$, and $\til P^\C_t$ is meromorphic with a pole at $\xi(t)$ in $U$. From (\ref{C1}) we have
\BGE \Imm [\pa_t\til J^\C_t(w)]+\Imm [(\til J^\C_t)'(w)\cot_2(w-\xi(t))]=\Imm [-C(t)\til P^\C_t(w)].\label{Im}\EDE
Comparing the residues at $\xi(t)$ of the two sides, we find that $C(t)=(\til J^\C_t)'(\xi(t))=\pa_y \til J_t(\xi(t))$.
Differentiating (\ref{Im}) w.r.t.\ $w$, we get
$$\pa_t(\til J^\C_t)'(w)+(\til J^\C_t)''(w)\cot_2(w-\xi(t))+(\til J^\C_t)'(w)\cot_2'(w-\xi(t))=-(\til J^\C_t)'(\xi(t))(\til P^\C_t)'(w).$$
Letting $w\to\xi(t)$ in $\til\Omega_t$ in the above formula,
 and comparing the constant term in the power series expansion at $\xi(t)$ of
both sides, we get
\BGE \pa_t(\til J^\C_t)'(\xi(t))=(\til J^\C_t)'(\xi(t))\lim_{w\to\xi(t)} (-(\til P^\C_t)'(w)-\cot_2'(w))+\frac 16(\til J_t^\C)'''(\xi(t)).
\label{patJ/J}\EDE
Let $\til Q_t=-\til P_t-\Imm\cot_2$. Then $\til Q_t$ is continuous on $\lin{\til\Om_t}$, vanishes on $\R$,
equals $-\Imm\cot_2$ on $\pa\til\Om_t\sem\R$, has period $2\pi$, and is harmonic inside $\til\Om_t$. From (\ref{patJ/J})
we have
\BGE {(\pa_t\pa_y/\pa_y) \til J_t(\xi(t))}=\pa_y\til Q_t(\xi(t))+\frac 16(\pa_x^2\pa_y/\pa_y)\til J_t(\xi(t)).\label{patJ/J-2}\EDE

For the behavior of $(\pa_t\pa_y/\pa_y) \til J_t(\xi(t))$ as $t\to -\infty$, we also need to consider two cases. The first case
is $D=\ha\C$. Then $\Om_t=\D$, so $P_t(z)=\Ree \frac{e^{i\xi(t)}+z}{e^{i\xi(t)}-z}$, which implies that
$\til P_t(z)=-\Imm\cot_2(z)$ and $\til Q_t\equiv 0$. From (\ref{Jj}) and (\ref{patJ/J-2}) we have
\BGE (\pa_t\pa_y/\pa_y) \til J_t(\xi(t))=O(t e^t), \quad t\to -\infty.\label{limit-patJ0}\EDE
The second case is that $D\ne\ha\C$. Let $Q_t$ be continuous on $\lin{\Om_t}$, harmonic
in $\Om_t$, vanishes on $\TT$, and equals to $\Ree \frac{e^{i\xi(t)}+z}{e^{i\xi(t)}-z}$ on $\pa\Om_t\sem\TT$.
Then $\til Q_t=Q_t\circ e^i$. Let $S_t=Q_t\circ R_\TT\circ \vphi_t$. Then $S_t$ is continuous on $\lin{D\sem K_t}$,
harmonic in $D\sem K_t$, vanishes on $\pa K_t$, and $S_t(z)=\Ree \frac{\vphi_t(z)+e^{i\xi(t)}}{\vphi_t(z)-e^{i\xi(t)}}$
on $\pa D$. Since $\vphi_t(z)=e^{-t}\phi_t(z)$, so from (\ref{phi_H(z)-z}) we have $\vphi_t(z)=e^{-t}z+O(1)$ as
$t\to -\infty$, uniformly in $z\in\pa D$. So $S_t=1+O(e^t)$ on $\pa D$ as $t\to -\infty$. Since $\diam(K_t)\le 4e^t$,
so $S_t(z)=1+O(1/t)$ as $t\to -\infty$, uniformly on any compact subset of $D\sem\{0\}$. The argument that is used to
derive (\ref{limit-J}) can be used here to prove that $\pa_y\til Q_t(\xi(t))=-1/t+O(1/t^2)$ as $t\to -\infty$. So from (\ref{Jj})
and (\ref{patJ/J-2}) we have
\BGE {(\pa_t\pa_y/\pa_y) \til J_t(\xi(t))}=-1/t+O(1/t^2),\quad t\to-\infty.\label{limit-patJ}\EDE

Let $\alpha=\lambda/\kappa\in\R$. We define $M(t)$ for $t\in(-\infty,T)$. If $D=\ha\C$, let
$$M(t)=(2\pi\pa_y \til J_t(\xi(t)))^\alpha \exp\Big(
-\frac\kappa 2\alpha(\alpha-1)\int_{-\infty}^t\Big(\frac{\pa_x\pa_y \til J_s(\xi(s))}{\pa_y\til J_s(\xi(s))}\Big)^2ds$$
\BGE -\frac \kappa 2 \alpha
\int_{-\infty}^t \frac{\pa_x^2\pa_y \til J_s(\xi(s))}{\pa_y\til J_s(\xi(s))}ds
-\alpha \int_{-\infty}^t \frac{\pa_t\pa_y \til J_s(\xi(s))}{\pa_y\til J_s(\xi(s))}ds\Big).\label{M(t)0}\EDE
From (\ref{Jj}) and (\ref{limit-patJ0}) we see that the three improper integrals all converge. From (\ref{limit-J0}) we have
$\lim_{t\to -\infty} M(t)=1$. If $D=\ha\C$, let
$$M(t)=\Big(\frac{\sqrt{t^2+1}\pa_y \til J_t(\xi(t))}{G(D,z_e;0)}\Big)^\alpha \exp\Big(
-\frac\kappa 2\alpha(\alpha-1)\int_{-\infty}^t\Big(\frac{\pa_x\pa_y \til J_s(\xi(s))}{\pa_y\til J_s(\xi(s))}\Big)^2ds$$
\BGE -\frac \kappa 2 \alpha
\int_{-\infty}^t \frac{\pa_x^2\pa_y \til J_s(\xi(s))}{\pa_y\til J_s(\xi(s))}ds
-\alpha \int_{-\infty}^t \Big(\frac{\pa_t\pa_y \til J_s(\xi(s))}{\pa_y\til J_s(\xi(s))}+\frac{s}{\sqrt{s^2+1}}\Big)ds\Big).\label{M(t)}\EDE
From (\ref{Jj}) and (\ref{limit-patJ}) we see that the three improper integrals all converge. From (\ref{limit-J}) we have
$\lim_{t\to -\infty} M(t)=1$ in this case.

\begin{Lemma} (Boundedness) Let $\rho$ be a Jordan curve in $\ha\C$ such that $0\in U(\rho)$ and $H(\rho)\subset D\sem\{z_e\}$.
Let $\tau_\rho$ be the first $t$ such that $K_t\cap \rho\ne\emptyset$. Then there is a constant $C\in(0,\infty)$ depending only on
$\rho$, $D$, and $z_e$, such that $|\ln(M(t))|\le C$ on $(-\infty,\tau_\rho]$. \label{bound}
\end{Lemma}
{\bf Proof.} Let $R_\rho=\dist(0,\rho)>0$. Then $\ln(R_\rho/4)$ is a lower bound of $\tau_\rho$. From (\ref{Jj}),
 (\ref{limit-J0}), (\ref{limit-patJ0}), and (\ref{M(t)0}), or from (\ref{Jj}), (\ref{limit-J}),   (\ref{limit-patJ}), and  (\ref{M(t)}),
 we conclude that there is $b\in(-\infty,\ln(R_\rho/4))$ and $C_1\in(0,\infty)$  depending only on $\rho$, $D$, and $z_e$,
  such that $|\ln(M(t))|\le C_1$ on $(-\infty,b]$. The boundedness of $|\ln(M(t))|$ on $[b,\tau_\rho]$ follows from Lemma
  \ref{compact-2}. $\Box$

\vskip 3mm

Now we study the martingale property of $M(t)$.  Since $(\til J_t)$ is $(\F_t)$-adapted, has period $2\pi$, and
$(e^{i\xi(t)})$ is also $(\F_t)$-adapted, so $(\pa_y\til J_t(\xi(t)))$ is  $(\F_t)$-adapted, and so are
$((\pa_x^j\pa_y/\pa_y)\til J_t(\xi(t)))$, $j=1,2$, and $((\pa_t\pa_y/\pa_y)\til J_t(\xi(t)))$.
From (\ref{M(t)0}) or (\ref{M(t)}) we see that $(M(t))$ is $(\F_t)$-adapted.
We will truncate the time interval to apply It\^o's formula. Recall that
$\ln(R/4)$ is a lower bound of $T$. Fix $a\in(-\infty,\ln(R/4))$. Let $T_a=T-a>0$. Let $\F^a_t=\F_{a+t}$, $t\ge 0$.
Then $T_a$ is an $(\F^a_t)_{t\ge 0}$-stopping time. Let $M_a(t)=M(a+t)$, $0\le t<T_a$. Then
$(M_a(t))$ is $(\F^a_t)$-adapted. From (\ref{M(t)0}) or (\ref{M(t)}) we have
$$M_a(t)=M(a)\pa_y\til J_a(\xi(a))^{-\alpha}\pa_y\til J_{a+t}(\xi(a+t))^\alpha
\exp\Big(
-\frac\kappa 2\alpha(\alpha-1)\int_{a}^{a+t}\Big(\frac{\pa_x\pa_y \til J_s(\xi(s))}{\pa_y\til J_s(\xi(s))}\Big)^2ds$$
\BGE -\frac \kappa 2 \alpha \int_{a}^{a+t} \frac{\pa_x^2\pa_y \til J_s(\xi(s))}{\pa_y\til J_s(\xi(s))}ds
-\alpha \int_{a}^{a+t}  \frac{\pa_t\pa_y \til J_s(\xi(s))}{\pa_y\til J_s(\xi(s))} ds\Big).\label{Ma(t)}\EDE
Let $\xi_a(t)=\xi(a+t)-\xi(a)$ and $B_a(t)=\xi_a(t)/\sqrt\kappa$, $t\in[0,\infty)$. Then $B_a(t)$ is an $(\F^a_t)$-Brownian motion.
Using It\^o's formula and the argument in Section \ref{section-conformal} or Section \ref{Section-martingale}, we conclude
that $\pa_y\til J_{a+t}(\xi(a+t))$, $0\le t<T_a$, satisfies the $(\F^a_t)$-adapted SDE:
$$d\pa_y\til J_{a+t}(\xi(a+t))=\pa_x\pa_y \til J_{a+t}(\xi(a+t))d \xi_a(t)+\frac{\kappa}2 \pa_x^2\pa_y \til J_{a+t}(\xi(a+t))dt
+\pa_t\pa_y \til J_{a+t}(\xi(a+t))dt.$$
From It\^o's formula, this then implies that
$$\frac{d\pa_y\til J_{a+t}(\xi(a+t))^\alpha}{d\pa_y\til J_{a+t}(\xi(a+t))}=\alpha\frac{\pa_x\pa_y \til J_{a+t}(\xi(a+t))}
{\pa_y \til J_{a+t}(\xi(a+t))} d \xi_a(t)+\frac{\kappa}2\alpha \frac{\pa_x^2\pa_y \til J_{a+t}(\xi(a+t))}
{\pa_y \til J_{a+t}(\xi(a+t))}dt$$ $$
+\alpha\frac{\pa_t\pa_y \til J_{a+t}(\xi(a+t))}{\pa_y \til J_{a+t}(\xi(a+t))}dt+\frac\kappa 2\alpha(\alpha-1)
\Big(\frac{\pa_x\pa_y \til J_{a+t}(\xi(a+t))}{\pa_y\til J_{a+t}(\xi(a+t))}\Big)^2dt.$$
So from (\ref{Ma(t)}) we see that $M_a(t)$, $0\le t<T_a$, is a local martingale, and
\BGE \frac{dM_a(t)}{M_a(t)}=\alpha\frac{\pa_x\pa_y \til J_{a+t}(\xi(a+t))}
{\pa_y \til J_{a+t}(\xi(a+t))} d \xi_a(t)=\lambda X(a+t)\frac{dB_a(t)}{\sqrt\kappa}.\label{M/M}\EDE
From Lemma \ref{bound} and that $\lim_{t\to-\infty}M(t)=1$ we conclude that for any $\rho$ as in Lemma \ref{bound},
we have $M(t)$, $-\infty<t\le \tau_\rho$, is a bounded martingale, and so $\EE_\mu[M(\tau_\rho)]=M(-\infty)=1$.

Define $\nu$ by $d\nu=M(\tau_\rho)d\mu$. Then $\nu$ is also a probability measure. Now suppose that the distribution of
$(\xi(t))$ is $\nu$ instead of $\mu$. For $-\infty<t<T$, let
\BGE \eta(t)=\xi(t)-\lambda\int_{-\infty}^t X^\xi(s)ds.\label{eta-xi}\EDE
From Proposition \ref{improper}, we see that $\eta(t)$, $-\infty<t<T$, are well defined. Moreover, it is clear that
$(e^{i\eta(t)})$ is $(\F_t)$-adapted. Fix $a\in(-\infty,\ln(R_\rho/4))$. Then we always have $\tau_\rho>a$. Define
$\eta_a(t)=\eta(a+t)-\eta(a)$ for $0\le t<T_a$. Then $(\eta_a(t))$ is $(\F^a_t)$-adapted. And we have
$$\eta_a(t)=\sqrt\kappa B_a(t)-\lambda\int_{a}^{a+t} X^\xi(s)ds.$$
From (\ref{M/M}) and Girsanov's theorem, we conclude that, under the measure $\nu$,
$(\eta_a(t)/\sqrt\kappa,0\le t\le \tau_\rho-a)$ is a stopped $(\F^a_t)$-adapted Brownian motion, and so is
independent of $e^{i\eta(t)}$, $-\infty<t\le a$. Since this holds for any $a\in(-\infty,\ln(R_\rho/4))$, so
$(e^{i\eta(t)},-\infty<t\le \tau_\rho)$ has the same distribution as $(e^{iB^{(\kappa)}_\R(t)})$ stopped
at some stopping time. Thus, there is a integer valued random variable $\bf n$ such that $\eta_*(t):=\eta(t)+2{\bf n}\pi$,
$-\infty<t\le \tau_\rho$, has the same distribution as $B^{(\kappa)}_\R(t)$ stopped
at some stopping time. From (\ref{eta-xi}) we see that $\xi_*(t)=\xi(t)+2{\bf n}\pi$, $-\infty<t\le \tau_\rho$, solves
the integral equation
$$\xi_*(t)=\eta_*(t)+\int_{-\infty}^t X^{\xi_*}(s)ds.$$
Here we use the fact that $X^{\xi_*}=X^\xi$. So the whole-plane Loewner chain
driven by $\xi_*(t)$, $-\infty<t\le \tau_\rho$, is a $(\kappa,\lambda)$ process stopped on hitting $\rho$. Thus, a
$(\kappa,\lambda)$ process stopped on hitting $\rho$ has a distribution that is absolutely continuous w.r.t.\ the
whole-plane SLE$_\kappa$ process stopped on hitting $\rho$, and the density function is $M(\tau_\rho)$.

The locally absolutely continuity also holds if the target is not an interior point but a boundary arc or a boundary point.
In these cases, the process $M(t)$ is defined by (\ref{M(t)}) with $G(D,z_e;0)$ replaced by a harmonic measure function or
a normalized Poisson kernel function valued at $0$. Then $M(\tau_\rho)$ is still the density function between the
$(\kappa,\lambda)$ process and whole-plane SLE$_\kappa$ process before hitting $\rho$.

\section{Scaling limits of Discrete LERW} \label{Convergence}
\subsection{Discrete LERW in grid approximation} \label{q_delta}
Let $D$ be a finitely connected domain that contains $0$.
For $\delta>0$, let $\delta\Z^2=\{(j+ik)\del:j,k\in\Z\}\subset\C$.
We also view $\delta \Z^2$ as a graph whose vertices are
$(j+ik)\del$, $j,k\in\Z$, and two vertices are adjacent iff the
distance between them is $\delta$. We define a graph ${\breve
D^\del}$ that approximates $D$ in $\del\Z^2$ as follows. The vertex
set $V({\breve D^\del})$ is the union of interior vertex set
$V_I({\breve D^\del})$ and boundary vertex set $V_\pa({\breve
D^\del})$, where $V_I({\breve D^\del}):=\del\Z^2\cap D$, and
$V_\pa({\breve D^\del})$ is the set of ordered pairs $\langle
z_1,z_2\rangle$ such that $z_1\in V_I({\breve D^\del})$, $z_2\in\pa
D$, and there is $z_3\in\del\Z^2$ that is adjacent to $z_1$ in
$\del\Z^2$, such that $[z_1,z_2)\subset[z_1,z_3)\cap D$. Two
vertices $w_1$ and $w_2$ in $V({\breve D^\del})$ are adjacent iff
either $w_1,w_2\in V_I({\breve D^\del})$, $w_1$ and $w_2$ are
adjacent in $\del\Z^2$, and $[w_1,w_2]\subset D$; or for $j=1$ or
$2$, $w_j\in V_I({\breve D^\del})$ and $w_{3-j}=\langle
w_j,z_3\rangle \in V_\pa({\breve D^\del})$ for some $z_3\in\pa D$.

Every interior vertex of ${\breve D^\del}$ has exactly $4$ adjacent
vertices, and every boundary vertex $w=\langle z_1,z_2\rangle $ has
exactly one adjacent vertex, which is the interior vertex $z_1$.
If $\langle z_1,z_2\rangle$ is
a boundary vertex, then it determines a boundary point, which is
$z_2$, and a prime end of $D$, which is the limit in $\ha D$, the conformal
closure of $D$ (c.f.\ \cite{Ahl} \cite{LERW}), as $z\to z_2$ along $[z_1,z_2)$.

Let $D^\del$ be the connected component of ${\breve D^\del}$ that
contains $0$. Let $V(D^\del)$ be the set of
vertices of $D^\del$. Let $V_I(D^\del):=V(D^\del)\cap V_I({\breve
D^\del})$ and $V_\pa(D^\del):=V(D^\del)\cap V_\pa({\breve D^\del})$
be the sets of interior vertices and boundary vertices, respectively, of
$D^\del$.

Fix $z_e\in D\sem\{0,\infty\}$. Let $w^\del_e$ be an interior vertex of $D^\del$ that is closest to $z_e$.
Then $|w^\del_e-z_e|<\del$ if $\del$ is small. Let $(q_\del(0),\dots,q_\del(\chi_\del))$ be the LERW on $D^\del$
started from $0$ conditioned to hit $w^\del_e$ before
$V_\pa(D^\del)$. Such LERW is obtained by the following process. First, run a simple random walk
on $D^\del$ from $0$, stop it on hitting $w^\del_e$ or $V_\pa(D^\del)$. Second, condition
the stopped walk on the event that it hits $w^\del_e$ instead of $V_\pa(D^\del)$. Finally, erase
the loops on the path of this walk, in the order they are created (c.f.\ \cite{Law}). Then the obtained simple lattice
path is called the LERW on $D^\del$ started from $0$ conditioned to hit $w^\del_e$ before
$V_\pa(D^\del)$. So $q_\del(0)=0$ and $q_\del(\chi_\del)=w^\delta_e$.

Let $E_{-1}=V_\pa(D^\del)$, $F=\{w^\del_e\}$, and $E_k=E_{-1}\cup\{q_j:0\le j\le k\}$
for $0\le k\le \chi_\del-1$. For each
$0\le k\le\chi_\del-1$, let $g_k$ be defined as in Lemma 2.1 in
\cite{LERW} with $A=F$, $B=E_{k-1}$ and $x=q_\del(k)$. This means that
$g_k$ is a function defined on $V(D^\delta)$, which vanishes on $E_k\sem \{q_\del(k)\}$, is discrete
harmonic on $V_I(D^\delta)\sem\{q_\del(0),\dots,q_\del(k)\}$, and $g_k(w^\delta_e)=1$.
The following is a special case of
Proposition 2.1 in \cite{LERW}.

\begin{Proposition} For any $v_0\in V(D^\del)$,
$(g_k(v_0))$ is a martingale up to the first time that
$q_\del(k)\sim w^\del_e$ or $E_k$ disconnects $v_0$ from $w^\del_e$. \label{mart hg}
\end{Proposition}

Define $q_\del$  on $[0,\chi_\del]$ to be the linear interpolation of $q_\del(k)$, $0\le k\le \chi_\del$.
Then $q_\del$ is a simple curve in $D$
that connects $0$ and $w^\delta_e$. For $0\le
k\le\chi_\del-1$, let $D_k=D\sem q_\del([0,k])$. When $\del$ is
small, the function $g_k$ approximates the generalized Poisson
kernel $P_k$ in $D_k$ with the pole at $q_\del(k)$, normalized by
$P_k(z_e)=1$. Note the resemblance of the discrete martingales
preserved by this (discrete) LERW given by Proposition \ref{mart hg} and the local martingales preserved by the
continuous LERW given by Theorem \ref{SLE2martingale}. Suppose $\gamma_0(t)$, $-\infty<t<T_0$, is an
$\LERW(D;0\to z_e)$ curve. We will prove the following theorem about the convergence.

\begin{Theorem} For any $\eps>0$, there is $\del_0>0$ such that, if $\del<\del_0$,
then there are a coupling of $q_\del$ and
$\gamma_0$, and a continuous increasing function ${\til u}$ that
maps $(0,\chi_\del)$ onto $(-\infty,T_0)$ such that
$$\PP[\sup\{|q_\del({\til u}^{-1}(t))-\gamma_0(t)|:-\infty < t<T_0\}<\eps]>1-\eps.$$
\label{global-coupling-in-D}
\end{Theorem}

\subsection{Some estimates}
For a non-degenerate interior hull $K\subset D\sem \{z_e\}$, let $\phi_K$, $\vphi_K$, and $\psi_K$
be as in the last subsection. So if $K=K^\xi_t$ is a whole-plane Loewner
hull at time $t$ driven by some $\xi\in C((-\infty,T))$ with $T>t$, then $\vphi_{K^\xi_t}$ and $\psi_{K^\xi_t}$
agree with the whole-plane Loewner map and inverted whole-plane Loewner map: $\vphi^\xi_t$ and $\psi^\xi_t$, respectively.
Let $\Om_K=R_{\TT}\circ \vphi_K(D\sem K)$,
$\til \Om_k=(e^i)^{-1}(\Om_k)$, $p_K=R_\TT\circ\vphi_K(z_e)$, and $\til p_K=(e^i)^{-1}(p_K)$.
 So $\Om_k$ is a subdomain of $\D$ containing $p_K$, and $\til\Om_K$ is a periodic
subdomain of $\HH$.  If  $K=K^\xi_t$, then
$\Om_{K^\xi_t}$, $\til\Om_{K^\xi_t}$, $p_{K^\xi_t}$, and $\til p_{K^\xi_t}$ agree with
$\Om^\xi_t$, $\til\Om^\xi_t$, $p^\xi_t$, and $\til p^\xi_t$,
 respectively, defined in Section  \ref{Def-of-quasi}.

Let $\alpha$ be a Jordan curve in $\C$ such that $0\in U(\alpha)$ and
$H(\alpha)\subset D\sem\{z_e\}$. Let $F$ be a compact subset of
$D\sem (H(\alpha)\cup\{\infty\})$.
Fix $b\in\R$. Throughout this subsection, a constant is called
uniform if it depends only on $D,z_e,\alpha,F,b$. We will frequently apply Lemma \ref{compact-2} to
${\cal H}^{b}(\alpha)$ to obtain some uniform constants. We illustrate the idea in the following
example. Note that for every $H\in {\cal H}^{b}(\alpha)$, $\vphi_H(F)$ is a compact subset of
$\{|z|>1\}$, so there is $r_H>0$ such that $|\vphi_H(z)|\ge e^{r_H}$ for every $z\in F$.
From Lemma \ref{compact-2}, there is a uniform constant $\h0>0$ such that
$|\vphi_H(z)|\ge e^{\h0}$ for any $H\in{\cal H}^{b}(\alpha)$ and $z\in F$.
Let $F_R=R_\TT(F)$. Then $|\psi_H(z)|\le e^{-\h0}$ for any $H\in{\cal H}^{b}(\alpha)$ and $z\in F_R$.
Suppose $K^\xi_a\subset H(\alpha)$. Then for any $t\in[b,a]$, we have $K^\xi_t\in H^{b}(\alpha)$,
so $|\psi^\xi_t(z)|\le e^{-\h0}$ for any $z\in F_R$. Let $\til F_R=(e^i)^{-1}(F_R)$. Since
$\psi^\xi_t\circ e^i=e^i\circ \til\psi^\xi_t$, so $\Imm \til \psi^\xi_t(z)\ge \h0$ for any $z\in\til F_R$ and
$t\in[b,a]$.

The following lemmas are
 similar to the lemmas in Section 6.1 of \cite{LERW}.

\begin{Lemma} There are uniform constants $C_1,C_2>0$ such that if
$K^\xi_a\subset H(\alpha)$, then for any $t_1\le t_2\in[b,a]$ and
$z\in \til F_R$,
$$|\til\psi^\xi_{t_2}(z)-\til\psi^\xi_{t_1}(z)|\le C_1|t_2-t_1|;$$
$$|\til\psi^\xi_{t_2}(z)-\til\psi^\xi_{t_1}(z)-{(t_2-t_1)}\cot_2({\til\psi^\xi_{t_1}(z)-\xi(t_1)}) |\le
C_2|t_2-t_1|(|t_2-t_1|+\sup_{t\in[t_1,t_2]}\{|\xi(t)-\xi(t_1)|\}).$$
\label{diff of vphi_2}
\end{Lemma}
{\bf Proof.} Suppose $K^\xi_a\subset H(\alpha)$. Then for any $t\in[b,a]$ and $z\in \til F_R$,
we have $\Imm \til\psi^\xi_t(z)\ge \h0$, which implies that $|\cot_2(\til\psi^\xi_t(z)-\xi(t))|\le \coth_2(\h0)$.
Since $\vphi^\xi_{t_2}(z)-\vphi^\xi_{t_1}(z)=\int_{t_1}^{t_2}
\cot_2(\vphi^\xi_t(z)-\xi(t))dt$, so
$|\til \psi^\xi_{t_2}(z)-\til \psi^\xi_{t_1}(z)|\le C_1|t_2-t_1|$ for any
$t_1\le t_2\in[b,a]$ and $z\in \til F_R$, where $C_1=\coth_2(\h0)>0$.
Since $| \cot_2'(w)|\le \frac 12\sinh_2^{-2}(\Imm w)\le \frac 12 \sinh_2^{-2}(\h0)$ for
$w\in\C$ with $\Imm w\ge \h0$, and $C_1\ge 1$, so for $t_1\le t_2\in[b,a]$ and $z\in \til F_R$,
$$|\cot_2({\til\psi^\xi_{t_2}(z)-\xi(t_2)})-\cot_2({\til\psi^\xi_{t_1}(z)-\xi(t_1)})|
\le \frac 12\sinh_2^{-2}(\h0) (|\vphi^\xi_{t_2}(z)-\vphi^\xi_{t_1}(z)|+|\xi(t_2)-\xi(t_1)|)$$
$$\le \frac {C_1}2\sinh_2^{-2}(\h0) (|t_2-t_1|+ |\xi(t_2)-\xi(t_1)|).$$
Let $C_2:=\frac {C_1}2\sinh_2^{-2}(\h0)>0$. Then for $t_1\le t_2\in[b,a]$ and
$z\in \til F_R$, we have
$$|\vphi^\xi_{t_2}(z)-\vphi^\xi_{t_1}(z)-{(t_2-t_1)}\cot_2({\til\psi^\xi_{t_1}(z)-\xi(t_1)})|$$
$$=\Big|\int_{t_1}^{t_2}
\cot_2({\vphi^\xi_t(z)-\xi(t)})-\cot_2({\vphi^\xi_{t_1}(z)-\xi(t_1)})\,dt\Big|$$
$$\le
C_2|t_2-t_1|(|t_2-t_1|+\sup_{t\in[t_1,t_2]}\{|\xi(t)-\xi(t_1)|\}).\quad\Box$$

\vskip 3mm

For   $x\in\R$, let $P(K,x,\cdot)$ be the generalized
Poisson kernel in $\Om_K$ with the pole at $e^{ix}$, normalized by $P(K,x,p_K)=1$. Let
$\til P(K,x,\cdot)=P(K,x,\cdot)\circ e^i$. If  $K=K^\xi_t$, then $P(K^\xi_t,x,\cdot)$,
and $\til P(K^\xi_t,x,\cdot)$ agree with $P^\xi(t,x,\cdot)$ and $\til P^\xi(t,x,\cdot)$,
 respectively, defined in Section \ref{Section-martingale}.

\begin{Lemma}  For each $n_1\in\{0,1\}$, $n_2,n_3\in\Z_{\ge 0}$,
there is a uniform constant $C>0$ depending on $n_1,n_2,n_3$, such
that if $K^\xi_a\subset H(\alpha)$, then for any $t\in[b,a]$,
$x\in\R$, and $z\in \til F_R$, we have $$
|\pa_1^{n_1}\pa_2^{n_2}\pa_{3,z}^{n_3}\til P^\xi(t,x,\til\psi^\xi_t(z))|\le
C.$$\label{derivatives of P}
\end{Lemma} \vskip -7mm
{\bf Proof.} The case $n_1=0$ follows from Lemma \ref{compact-2} immediately
because $K^\xi_t\in {\cal H}^{b}(\alpha)$ for $t\in[b,a]$, and
 if $(H_n)$ is a sequence in ${\cal H}^{b}(\alpha)$, and $H_n\hto H$, then
$$\pa_2^{n_2}\pa_{3,z}^{n_3}\til P(H_n,x,\til\psi_{H_n}(z))\to
\pa_2^{n_2}\pa_{3,z}^{n_3}\til P(H,x,\til\psi_H(z))$$
uniformly in $x\in\R$ and $z\in \til F_R$, for any $n_2,n_3\in\Z_{\ge 0}$.

Now we consider the case $n_1=1$. First suppose that $\pa D$ is analytic, i.e., $\pa D$ is the disjoint union of
 analytic Jordan curves.
Let $K\in {\cal H}^{b}(\alpha)$. Then $\Om_K=R_{\TT}\circ \vphi_K(D\sem K)$
also have analytic boundary. Since $P(K,x,\cdot)$ vanishes on $\pa \Om_K$
except at $e^{ix}$, so $P(K,x,\cdot)$ extends harmonically across $\pa \Om_k\sem\{e^{ix}\}$. Thus, $\til P(K,x,\cdot)$
extends harmonically across $\pa\til\Om_K\sem\{x+2n\pi:n\in\Z\}$. For
$x,y\in\R$, let $Q_y(K,x,\cdot)$ be a continuous function on $\lin{\Om_K}\sem
\{e^{ix}\}$ such that $Q_y(K,x,\cdot)$ is harmonic in $\Om_K$; vanishes
on $\TT\sem\{x\}$; behaves like $c\Ree\frac{e^{ix}+z}{e^{ix}-z}+O(1)$ near $e^{ix}$
for some $c\in\R$; and
$$Q_y(K,x,z)=-2\Ree(\pa_{3,z}P(K,x,z) z \frac{e^{iy}+z}{e^{iy}-z}),\quad
z\in(\pa\Om_K\sem\TT)\cup \vphi_K(p).$$
Such $Q_y(K,x,\cdot)$ exists uniquely. Let $\til Q_y(K,x,\cdot)=Q_y(K,x,\cdot)\circ e^i$.
From (\ref{reflection-whole-plane-equation}) and the values of $P^\xi(t,x,\cdot)$ at $\pa\Om^\xi_t\sem \TT
=\psi^\xi_t(\pa\Om)$ and $p^\xi_t=\psi^\xi_t(p)$,  it is easy to check that
$\pa_1 P^\xi(t,x,z)=Q_{\xi(t)}(K^\xi_t,x,z)$, and so $\pa_1 \til P^\xi(t,x,z)=\til Q_{\xi(t)}(K^\xi_t,x,z)$.
Using Lemma \ref{compact-2}, we can conclude that for any $n_2,n_3\in\Z_{\ge 0}$,
 $\pa_2^{n_2}\pa_{3,z}^{n_3} \til Q_y(K,x,\psi_K(z))$ is uniformly bounded in $x,y\in\R$ and
 $z\in \til F_R$. So the proof in the case that $n_1=1$ and $\pa D$ is analytic is finished.

Now we consider the case that $n_1=1$ but $\pa D$ may not be analytic. We may
find $V$ that maps $D$ conformally onto $D_0$ with analytic boundary, such
that $V(0)=0$. Moreover, suppose $F_0:=V(F)$, $\alpha_0:=V(\alpha)$, and $V(H(\alpha))$ do not contain $\infty$,
and $z_0:=V(z_e)\ne \infty$.
Then $\alpha_0$ is a Jordan curve in $\C$ such that $0\in U(\alpha_0)$, $H(\alpha_0)=V(H(\alpha))\subset D_0\sem\{z_0\}$,
and $F_0$ is a compact subset of $D_0\sem(H(\alpha_0)\cup\{\infty\})$.
Let $W=R_\TT\circ V\circ W$, $\Om_0=R_\TT(D_0)$, and $p_0=R_\TT(z_0)$.
Then $W$ maps $\Om$ conformally onto $\Om_0$. Let $\til \Om_0=(e^i)^{-1}(\Om_0)$.
Choose $\til W$ that maps $\til \Om$ conformally onto $\til\Om_0$ such that
$W\circ e^i=e^i\circ \til W$.
There is $b_0\in\R$ such that if $H$ is an interior
 hull in $D$ with $0\in H$ and $\ccap(H)\ge b$, then $\ccap(V(H))\ge b_0$.

Suppose $K^\xi_a\subset {\cal H}(\alpha)$.
Using the argument in Section \ref{section-conformal}, we conclude that there are $a_0\in\R$,
$\xi_0\in C((-\infty,a_0])$, and a continuous increasing function $u$ that maps $(-\infty,a]$ onto
$(-\infty,a_0]$ such that $V(K^\xi_t)=K^{\xi_0}_{u(t)}$ for $-\infty<t\le a$. Then we have
$K^{\xi_0}_s\in {\cal H}(\alpha_0)$ for $-\infty<s\le a_0$, and $u(b)\ge b_0$. Let
\BGE W_t=\psi^{\xi_0}_{u(t)}\circ W\circ (\psi_{t}^\xi)^{-1},\quad
\til W_t=\til \psi^{\xi_0}_{u(t)}\circ \til W\circ(\til \psi_{t}^\xi)^{-1},\quad -\infty<t\le a.\label{Wt-comp}\EDE
Using the argument in Section \ref{section-conformal}, we can conclude that $W_t\circ e^i=e^i\circ \til W_t$,
  $u'(t)=\til W_t'(\xi(t))^2$, $\xi_0(u(t))=\til W_t(\xi(t))$, and for any $w\in\til\Om_t$,
$$\pa_t\til W_t(w)=\til W_t'(\xi(t))^2\cot_2({\til W_t(w)-\til W_t(\xi(t))})-\til
W_t'(w)\cot_2({w-\xi(t)}).$$
So we have that for any $z\in\Om_t$,
\BGE \pa_t  W_t(z)= |W_t'(e^{i\xi(t)})|^2 W_t(z)\frac{W_t(e^{i\xi(t)})+W_t(z)}{W_t(e^{i\xi(t)})-W_t(z)}-
W_t'(z)z\frac{e^{i\xi(t)}+z}{e^{i\xi(t)}-z}. \label{n1=1}\EDE

For $-\infty<t\le a_0$ and $x\in\R$, let $P^{\xi_0}_0(t,x,\cdot)$ be the generalized Poisson kernel
in $\psi^{\xi_0}_t(\Om_0\sem L^{\xi_0}_t)$ with the pole at $e^{ix}$, normalized by
$P^{\xi_0}_0(t,x,\psi^{\xi_0}_t(p_0))=1$; and let $\til P^{\xi_0}_0(t,x,\cdot)=P^{\xi_0}_0(t,x,\cdot)\circ e^i$.
Then we have \BGE \til P^\xi(t,x,z)=\til P^{\xi_0}_0(u(t),\til W_t(x),\til W_t(z)),\quad -\infty<t\le a.\label{P=P}\EDE
Let $\til F_{0,R}=(e^i)^{-1}(R_\TT(F_0))$.
Since $D_0$ has analytic boundary, so for any $n_1\in\{0,1\}$, $n_2,n_3\in\Z_{\ge 0}$, there is a uniform constant
$C$ depending on $n_1$, $n_2$, and $n_3$ such that for $t\in [b_0,a_0]$ and $z\in \til F_{0,R}$,
$$|\pa_1^{n_1}\pa_2^{n_2}\pa_{3,z}^{n_3} \til P^{\xi_0}_0(t,x,\til \psi^{\xi_0}_t(z))|\le C.$$
From (\ref{P=P}) and that $u([b,a])\subset [b_0,a_0]$ and $u'(t)=\til W_t'(\xi(t))^2$, we suffice to prove that for any $n_1\in\{0,1\}$
and $n_2\in\Z_{\ge 0}$ with $n_1+n_2\ge 1$, there is a uniform constant $C$ depending on $n_1$ and $n_2$ such that
$|\pa_t^{n_1}\pa_z^{n_2}\til W_t(z)|\le C$ for any $t\in[b,a]$ and $z\in\R\cup \til\psi^\xi_t(\til F_R)$. Since
$e^i\circ \til W_t=W_t\circ e^i$, so we suffice to prove that there is a uniform constant $\delta_0>0$ such that
$|W_t(z)|>\delta_0$ for $t\in[b,a]$ and  $z\in\TT\cup \psi^\xi_t(F_R)$; and for any $n_1\in\{0,1\}$
and $n_2\in\Z_{\ge 0}$ with $n_1+n_2\ge 1$, there is a uniform constant $C_0$ depending on $n_1$ and $n_2$ such that
$|\pa_t^{n_1}\pa_z^{n_2} W_t(z)|\le C_0$ for any $t\in[b,a]$ and $z\in\TT\cup \psi^\xi_t( F_R)$.

For the existence of $\delta_0$, we consider two cases. The first case is $z\in\TT$. This is  trivial   because
$|W_t(z)|=1$ on $\TT$. The second case is $z\in \psi^\xi_t(F_R)$. From (\ref{Wt-comp}) and that
$\psi^{\xi_0}_t=R_\TT\circ \vphi^\xi_t\circ R_\TT$,  the inequality in this case
 is equivalent to that $|\vphi^{\xi_0}_{u(t)}(z)|\le 1/\delta_0$ for any $t\in[a,b]$ and $z\in  F_0$.
This can be proved by applying Lemma \ref{compact-2} to ${\cal H}^{b_0}(\alpha_0)$ and using the facts that
$\vphi^{\xi_0}_{u(t)}=\vphi_{K^{\xi_0}_{u(t)}}$, $K^{\xi_0}_{u(t)}\in{\cal H}^{b_0}(\alpha_0)$
for $t\in[b,a]$, and $\infty\not\in\vphi_H(F_0)$ for every $H\in {\cal H}^{b_0}(\alpha_0)$.

Next we consider the existence of $C_0$. We first consider the case $n_1=0$. For any $n_2\in\Z_{\ge 0}$,
the uniform  boundedness of $\pa_z^{n_2} W_t(z)$ on $\psi^\xi_t(F_R)$ follows immediately from Lemma \ref{compact-2}
applied to ${\cal H}^{b}(\alpha)$ and ${\cal H}^{b_0}(\alpha_0)$. Using Lemma \ref{compact-2} we may also
obtain uniform numbers $r\in(0,1)$ and $M\in(0,\infty)$
such that for $t\in[b,a]$, we have $\{r\le |z|<1\}\subset \Om^\xi_t$, and $|W_t(z)|\le M$ on $\{|z|=r\}$. Then the
 uniform  boundedness of $\pa_z^{n_2} W_t(z)$ on $\TT$ follows from Cauchy's integral formula.
 A similar argument together with (\ref{n1=1}) proves the case $n_1=1$. The two fractions in (\ref{n1=1}) do not cause
 any problem because they are uniformly bounded as long as $z$ and $W_t(z)$ are uniformly bounded away from $\TT$,
  which are true for $z\in \psi^\xi_t(F_R)$ and $z\in\{|z|=r\}$. $\Box$

\begin{Lemma} There is a uniform constant $C>0$ such that
if $K^\xi_a\subset H(\alpha)$, then for any $t,t'\in[b,a]$,
$|X^\xi(t)|\le C$ and $|X^\xi(t)-X^\xi({t'})|\le
C(|t-t'|+|\xi(t)-\xi(t')|)$. \label{bound of R}
\end{Lemma}
{\bf Proof.}  Suppose $K^\xi_a\subset H(\alpha)$. Write
$\til J^\xi(t,x)$ for $\til J^\xi_t(x)$. Note that $X^\xi(t)=(\pa_{2,z}^2/\pa_{2,z})
\til J^\xi(t,\xi(t))$. So it suffices to prove that there is a uniform constant $C>0$ such
that for any $t\in[b,a]$ and $x\in\R$,
$|\pa_1^{n_1}\pa_{2,z}^{n_2} (\pa_{2,z}^2/\pa_{2,z}) \til J^\xi(t,x)|\le
C$ for $n_1,n_2\in\{0,1\}$. We need to show that $|\pa_{2,z} \til
J^\xi(t,x)|$ is bounded from below by a positive uniform constant,
and $|\pa_1^{n_1}\pa_{2,z}^{n_2+1} \til J^\xi(t,x)|$ is bounded from
above by a positive uniform constant. The proof is similar to that
of the above lemma. $\Box$

\begin{Lemma} There is a uniform constant $C>0$ such that if
$K^\xi_a\subset H(\alpha)$, then for any $t_1\le t_2\in[b,a]$ and
$z\in \til F_R$, we have
$$|\pa_1\til P^\xi(t_2,\xi(t_2),\til\psi^\xi_{t_2}(z))-\pa_1P^\xi(t_1,\xi(t_1),\til\psi^\xi_{t_1}(z))|\le
C(|t_2-t_1|+|\xi(t_2)-\xi(t_1)|).$$ \label{difference of d1}
\end{Lemma}\vskip -7mm
{\bf Proof.} This follows from Lemma \ref{vanish*}, and the above
three lemmas. $\Box$

\begin{Lemma} There is a uniform constant $d_1>0$
such that, if $K^\xi_a\subset H(\alpha)$, then for any $z\in \til F_R$, and any
$t_1<t_2\in[b,a]$ that satisfy $|t_2-t_1|\le d_1$,  we have
$$\til P^\xi(t_2,\xi(t_2),\til\psi^\xi_{t_2}(z))-\til P^\xi(t_1,\xi(t_1),\til\psi^\xi_{t_1}(z))$$
$$=\pa_2\til P^\xi(t_1,\xi(t_1),\til \psi^\xi_{t_1}(z))\cdot[(\xi(t_2)-\xi(t_1))-
 (t_2-t_1)X^\xi_{t_1}]$$$$ +\frac 12
\pa_2^2\til
P^\xi(t_1,\xi(t_1),\til\psi^\xi_{t_1}(z))\cdot[(\xi(t_2)-\xi(t_1))^2-2(t_2-t_1)]$$$$+O(A^2)
+O(AB)+O(AB^2) +O(B^3),$$ \label{estimate for P} where
$A:=|t_2-t_1|$, $B:= \sup_{s,t\in[t_1,t_2]}\{|\xi(s)-\xi(t)|\}$, and
$O(X)$ is some number whose absolute value is bounded by $C|X|$ for
some uniform constant $C>0$.\label{Taylor}
\end{Lemma}
{\bf Proof.}  We may choose a compact subset $F'$ of $\D\sem
H(\rho)$ such that $F$ is contained in the interior of $F'$. Let $F'_R=R_\TT(F)$ and $\til F'_R=(e^i)^{-1}(F'_R)$.
So $F_R$ and $\til F_R$ are contained in the interiors of $F'_R$ and $\til F'_R$, respectively.
Applying Lemma \ref{compact-2} to ${\cal H}^{b}(\alpha)$,
we obtain a uniform constant $d_0>0$ such that  for any $K\in {\cal H}^{b}(\alpha)$, we have
$\dist(\psi_K(F_R),\pa\psi_K(F'_R))\ge d_0$.  So there is a uniform constant $\til d_0$ such
that $\dist(\til \psi_K(\til F_R),\pa \til\psi_K(\til F'_R))\ge \til d_0$
for any $K\in {\cal H}^{b}(\alpha)$. Suppose $K^\xi_a\subset
H(\alpha)$. From Lemma \ref{diff of vphi_2} and the existence of $\til d_0$, we get a uniform
constant $d_1>0$ such that if $s,t\in[b,a]$ satisfy $|s-t|\le d_1$
then for any $z\in \til F_R$,
$[\til\psi^\xi_s(z),\til\psi^\xi_t(z)]\subset\til\psi^\xi_s(\til F'_R)$.

Fix $z\in \til F_R$ and $t_1<t_2\in[0,a]$ with $|t_2-t_1|\le d_1$.  Let
$P_1=\til P^\xi(t_2,\xi(t_2),\til\psi^\xi_{t_2}(z))$,
$P_2=\til P^\xi(t_1,\xi(t_2),\til\psi^\xi_{t_2}(z))$,
$P_3=\til P^\xi(t_1,\xi(t_1),\til\psi^\xi_{t_2}(z))$,
$P_4=\til P^\xi(t_1,\xi(t_1),\til\psi^\xi_{t_1}(z))$. Then
\BGE \til P^\xi(t_2,\xi(t_2),\til\psi^\xi_{t_2}(z))-\til P^\xi(t_1,\xi(t_1),\til\psi^\xi_{t_1}(z))
=(P_1-P_2)+(P_2-P_3)+(P_3-P_4).\label{P-P}\EDE

Now $P_1-P_2=\int_{t_1}^{t_2}\pa_1
\til P^\xi(t,\xi(t_2),\til\psi^\xi_{t_2}(z))dt$. Fix any $t\in[t_1,t_2]$.
Applying Lemma \ref{diff of vphi_2} and Lemma \ref{derivatives of P}
to $\til F'_R$  and using
$[\til\psi^\xi_t(z),\til\psi^\xi_{t_2}(z)]\subset\til\psi^\xi_t(\til F'_R)$,  we
have
$$\pa_1 \til P^\xi(t,\xi(t_2),\til\psi^\xi_{t_2}(z))-\pa_1
\til P^\xi(t,\xi(t),\til\psi^\xi_t(z))=O(A)+O(B).$$ Applying Lemma
\ref{difference of d1} to $\til F_R$, we have
$$\pa_1 \til P^\xi(t,\xi(t),\til\psi^\xi_t(z))-\pa_1
\til P^\xi(t_1,\xi(t_1),\til\psi^\xi_{t_1}(z))=O(A)+O(B).$$ So we get
$$P_1-P_2=\pa_1 \til P^\xi(t_1,\xi(t_1),\til\psi^\xi_{t_1}(z))(t_2-t_1)+O(A^2)+O(AB).$$
Applying Lemma \ref{derivatives of P} to $\til F'_R$, since
$\til\psi^\xi_{t_2}(z)\in\til\psi^\xi_{t_1}(\til F'_R)$, so we have
$$P_2-P_3=\pa_2 \til P^\xi(t_1,\xi(t_1),\til\psi^\xi_{t_2}(z))(\xi(t_2)-\xi(t_1))
$$$$+\frac 12\pa_2^2
\til P^\xi(t_1,\xi(t_1),\til\psi^\xi_{t_2}(z))(\xi(t_2)-\xi(t_1))^2+O(B^3).$$
Applying Lemma \ref{diff of vphi_2} and Lemma \ref{derivatives of P}
to $\til F'_R$, since
$[\til\psi^\xi_{t_1}(z),\til\psi^\xi_{t_2}(z)]\subset\til\psi^\xi_{t_1}(\til F'_R)$,
so  we have
$$\pa_2^j \til P^\xi(t_1,\xi(t_1),\til\psi^\xi_{t_2}(z))-\pa_2^j \til P^\xi(t_1,\xi(t_1),\til\psi^\xi_{t_1}(z))=O(A),$$
for $j=1,2$. Thus
$$P_2-P_3=\pa_2 \til P^\xi(t_1,\xi(t_1),\til\psi^\xi_{t_1}(z))(\xi(t_2)-\xi(t_1))
$$$$+\frac 12\pa_2^2
\til P^\xi(t_1,\xi(t_1),\til\psi^\xi_{t_1}(z))(\xi(t_2)-\xi(t_1))^2
+O(AB)+O(AB^2)+O(B^3).$$ Applying Lemma \ref{diff of vphi_2} and
Lemma \ref{derivatives of P} to $\til F'_R$, since
$[\til\psi^\xi_{t_1}(z),\til\psi^\xi_{t_2}(z)]\subset\til\psi^\xi_{t_1}(\til F'_R)$,
so we have
$$P_3-P_4=2\Ree(\pa_{3,z}\til P^\xi(t_1,\xi(t_1),\til\psi^\xi_{t_1}(z))(\til\psi^\xi_{t_2}(z)-\til\psi^\xi_{t_1}(z)))+O(A^2)$$
$$=2\Ree(\pa_{3,z}\til P^\xi(t_1,\xi(t_1),\til\psi^\xi_{t_1}(z)){(t_2-t_1)}\cot_2({\til\psi^\xi_{t_1}(z)-\xi(t_1)}))
+O(AB)+O(A^2).$$ The conclusion then follows from (\ref{P-P}) and Lemma \ref{vanish*}. $\Box$

\subsection{Convergence of driving functions}

We may choose mutually disjoint Jordan curves $\alpha_j$, $j=0, 1,2$,
in $\C$ such that $0\in U(\alpha_0)\subset
U(\alpha_1)\subset U(\alpha_2)$ and $H(\alpha_2) \subset D\sem\{z_e\}$. Fix $b\in\R$ such that $b<\ln(d_0/4)-1$,
where $d_0=\dist(0,\alpha_0)$. So any $H\in{\cal H}_0$ with $\ccap(H)\le b$ must satisfy
$H\subset U(\alpha_0)$. Let $F$ be a compact subset of
$D\sem H(\alpha_2)$ whose interior is not empty. From now on, a uniform constant
is a number that depends only  on $D,z_e,\alpha_0,\alpha_1,\alpha_2,F,b$, and some other variables
we will specify. Let $O(X)$ denote some number whose absolute value is bounded by $C|X|$ for some
uniform constant $C>0$.

Let $L^\delta$ denote the set of simple lattice paths
$X=(X(0),\dots,X(s))$, $s\in\N$, on $D^\delta$, such that $X(0)=0$,
 $X(k)\in D$ for $0\le k\le s$,
 and $\bigcup_{k=0}^s(X({k-1}),X(k)]\subset H(\alpha_1)$.
Let $\Set(X)=\{X(0),\dots,X(s)\}$, $\Tip(X)=X(s)$,
$H_X=\bigcup_{k=1}^{s} [X({k-1}),X(k)] $, and $D_X=D\sem H_X$. Let $P_X$ be
the generalized Poisson kernel in $D_X$ with the pole at $\Tip(X)$,
normalized by $P_X(z_e )=1$, and $g_X$ be defined on $V(D^\delta)$
such that $g_X\equiv 0$ on $V_\pa(D^\delta)\cup
\Set(X)\sem\{\Tip(X)\}$, $\Delta_{D^\delta} g_X\equiv 0$ on
$V_I(D^\delta)\sem\Set(X)$, and $g_X(w^\del_e )=1$. Let
$L^\del_{b}$ be the set of $X\in L^\del$ such that $\ccap(H_X)\ge
b$. Then we have the following proposition about the convergence of
$g_X$ to $P_X$.

\begin{Proposition} For any $\eps>0$, there is a uniform constant $\delta_0>0$
depending on $\eps$ such
that, if $0<\delta<\delta_0$, then for any $X\in{L}_{b}^\delta$,
and any $w\in V(D^\delta)\cap (D\sem H(\alpha_2))$, we have
$|g_X(w)-P_X(w)|<\eps$.\label{approx of P}
\end{Proposition}
{\bf Sketch of the proof.} This proposition is similar to Proposition 6.1 in \cite{LERW}.
So we only give a sketch of the  proof. Suppose that the proposition is not true. Then
there are $\eps_0>0$, a sequence $\delta_n\to 0^+$, a sequence of lattice paths $X_n\in L_b^{\delta_n}$,
and a sequence of lattice points $w_n\in V(D^\delta)\cap (D\sem H(\alpha_2))$, such that
$|g_{X_n}(w_n)-P_{X_n}(w_n)|\ge \eps_0$. By passing to a subsequence, one may assume that
$w_n\to w_0$, and $D_{X_n}\dto D_0$. Then $P_{X_n}$ tends to a generalized Poisson kernel function
in $D_0$. Using linear interpolation to extend each $g_{X_n}$ to a
continuous function defined in the unions of lattice squares inside $D_{X_n}$.
Since each $g_{X_n}$ is a positive harmonic function, so from Harnack's inequality, we can conclude that the
extended $\{g_{X_n}\}$ is uniform Lipschitz on any compact subset of $D_0$. Applying Arzel\`a-Ascoli
theorem, by passing to a subsequence, we conclude that $g_{X_n}\to g_0$ locally uniformly in $D_0$.
Then one can check that $g_0$ is a positive harmonic function. With a little more work, one can prove that
$g_0$ is also a generalized Poisson kernel, and in fact, $g_0=P_0$.
So if $w_0=\lim w_n\in D_0$, we immediately get a contradiction.
If $w_0\not\in D_0$, then $w_0\in\pa D$. From $w_n\to \pa D$ we get $g_{X_n}(w_n)\to 0$ and
$P_{X_n}(w_n)\to 0$, which also gives a contradiction. $\Box$

\vskip 3mm

Let the LERW curve $q_\del$ on $[0,\chi_\del]$ be defined as in
Section \ref{q_delta}. For $0\le t\le\chi_\del$, let
$v_\del(t)=\ccap( q_\del([0,t]))$, and
$T_\del=v_\del(\chi_\del)$. Then $v_\delta$ is an increasing function,
and maps $[0,\chi_\del]$ onto $[-\infty, T_\del]$. Let
$\beta_\del(t)=q_\del(v_\del^{-1}(t))$,  $-\infty\le t\le T_\del$. From Proposition \ref{Whole-plne-Proposition-I},
there is some $\xi_\del\in C((-\infty,T_\del])$ such that
$\beta_\del([-\infty,t])=K^{\xi_\del}_t$ for $-\infty< t\le T_\del$. Let
$n_\infty$ be the first $n$ such that $(q_\del({n-1}),q_\del(n)]$
intersects $\alpha_0$. We may choose $\del<\dist(\alpha_0,\alpha_1)$. Then
$q_\del([0,n_\infty])\subset U(\alpha_1)$.  Let $T^\del_{\alpha_0}=v_\del(n_\infty)$.
Let $n_0$ be the first $n$ such that $v_\delta(n)\ge b$.
Pick any $d>0$. Define a sequence $(n_j)$ by the following.
For $j\ge 1$, let $n_{j+1}$ be the first $n\ge n_j$ such that
$n=n_\infty$, or $v_\del(n)-v_\del(n_j)\ge d^2$, or
$|{\xi_\del}(n)-{\xi_\del}(n_j)|\ge d$, whichever comes first.
Let $(\F_n)$ be the filtration generated by $(q_\del(n))$.
Let ${\cal F}_j'={\cal F}_{n_j}$, $0\le j<\infty$. Then we may derive the following
proposition, which is similar to Proposition 6.2 in \cite{LERW}. Since the proofs of these two
propositions are almost identical, so we omit the proof here.

\begin{Proposition} There are a uniform constant $d_0>0$
and a uniform constant $\del_0(d)>0$ that depends only on $d$ such that,
if $d<d_0$ and $\del<\del_0(d)$, then for all $j\ge 0$,
$$\EE[({\xi_\del}(v_\del(n_{j+1}))-{\xi_\del}(v_\del(n_j)))-\int_{v_\del(n_j)}^{v_\del(n_{j+1})}
X^{\xi_\del}_t dt| {\cal F}_j']=O(d^3);$$
$$\EE[({\xi_\del}(v_\del(n_{j+1}))-{\xi_\del}(v_\del(n_j)))^2-2(v_\del(n_{j+1})-v_\del(n_j))|
{\cal F}_j']=O(d^3).$$
\label{key estimate}
\end{Proposition}

Let $\xi_0(t)$, $-\infty< t<T_0$, be the maximal solution to
$$\xi_0(t)=B^{(2)}_\R(t)+2\int_{-\infty}^t
X^{\xi_0}_sds,$$
where $B^{(2)}_\R(t)$, $t\in\R$, is defined in Section \ref{Def-of-quasi}.
Let $\beta_0(t)$, $-\infty< t<T_0$, be the whole-plane Loewner curve driven
by $\xi_0$. Then $\beta_0$ is a continuous $\LERW(D;0\to z_e)$ curve.

If $\alpha$ is a Jordan curve in $\C$ with $0\in U(\alpha)$, and $\beta$ defined on $[-\infty,T)$
is a curve in $\C$ with $\beta(-\infty)=0$, let $T_\alpha(\beta)$ be the first $t$ such
that $\beta(t)\in\alpha$, if such $t$ exists; otherwise let
$T_\alpha(\beta)=T$. Since $q_\del([0,T^\del_{\alpha_0}])$
intersects $\alpha_0$, so $T_{\alpha_0}(q_\del)\le
T^\del_{\alpha_0}$.
Using the above proposition, we are able to derive the following theorem, which is similar to Theorem 6.2 in \cite{LERW}. The proof
uses Skorokhod Embedding Theorem, the method in the proof of Theorem 3.7 in \cite{LSW-2}, and the Markov property of
 $(e^i(B^{(2)}_\R(t)))$. Again, we omit the proof here.

\begin{Theorem} Suppose $\alpha$ is a Jordan curve in $\C$ with $0\in U(\alpha)$,
and $H(\alpha)\subset D\sem\{z_e\}$. For
every $b\in\R$ and $\eps>0$, there is $\delta_0>0$ such that if $\delta<\delta_0$
then there is a coupling of the processes $(\xi_\del(t))$ and
$(\xi_0(t))$ such that $$\PP[\sup\{|e^i(\xi_\del(t))-e^i(\xi_0(t))|:t\in
[b,T_\alpha(\beta_\del)\vee
T_\alpha(\beta_0)]\}<\eps]>1-\eps.$$ Here if
$\xi_\del$ or $\xi_0$ is not defined on $[b,T_\alpha(\beta_\del)\vee
T_\alpha(\beta_0)]$, we set the value of $\,\sup$ to be $+\infty$.
 \label{symmetric coupling}
\end{Theorem}

\subsection{Convergence of the curves}

So far, we have derived the convergence of the driving functions. Using the above theorem, Lemma \ref{compact-2},
and the regularity of discrete LERW path (c.f.\ Lemma 3.4 in \cite{S-SLE} and Lemma 7.2 in \cite{LERW}),
we may derive the following theorem, which is similar to Theorem 7.1 in \cite{LERW}.
It is about the local convergence of the curves. Here we omit its proof.

\begin{Theorem} Let $\alpha$ be as in the above theorem. For
every $\eps>0$, there is $\delta_0>0$ such that if $\delta<\delta_0$
then there is a coupling of the processes $(\beta_\del(t))$ and
$(\beta_0(t))$ such that $$\PP[\sup\{|\beta_\del(t)-\beta_0(t)|:t\in
[-\infty,T_\alpha(\beta_\del)\vee
T_\alpha(\beta_0)]\}<\eps]>1-\eps.$$
 \label{symmetric coupling-trace}
\end{Theorem}

Finally, we may lift the local convergence to
the global convergence, and so finish the proof of Theorem \ref{global-coupling-in-D}.
The argument used here is almost identical to that in Section 7.2 of \cite{LERW}. A slight difference is
that now $\cal A$ is the set of Jordan curves $\alpha$ such that $0\in U(\alpha)$ and
$H(\alpha)\subset D\sem\{z_e\}$; and $\cal B$ is the set of continuous curves $\beta:
[-\infty,T)\to D$ for some $T\in\R$, with $\beta(-\infty)=0$.

\subsection{Other kinds of targets}
Let $D$ be a finitely connected domain that contains $0$.
Suppose $w_e$ is a prime end of $D$ that satisfies $w_e\in\del_e\Z^2$ for
some $\del_e>0$, and $\pa D$ is flat near $w_e$, which means that
there is $r>0$ such that
$D\cap\{z\in\C:|z-w_e|<r\}=(w_e+a\HH)\cap\{z\in\C:|z-w_e|<r\}$ for
some $a\in\{\pm 1,\pm i\}$. For $\del>0$, let $w^\del_e=w_e+ia\del$.

Let ${\cal M}$ be the set of $\delta>0$ such that
$w_e\in\delta\Z^2$. If $\del\in\cal M$ is small enough, then
$\langle w^\del_e,w_e\rangle$ is a boundary vertex of $ {\breve
D^\del}$, which determines the boundary point and prime end $w_e$,
and there is a lattice path on $D^\del$ that connects $0$ with
$w_e$ without passing through any other boundary vertex.  Here we do
not distinguish $w_e$ from the boundary vertex $\langle
w^\del_e,w_e\rangle$. Let $F=\{w_e\}$, $E_{-1}=V_\pa(D^\del)\sem
F$, and $E_k=E_{-1}\cup\{q_j:0\le j\le k\}$
for $0\le k\le \chi_\del-1$. Let
$(q_\del(0),\dots,q_\del(\chi_\del))$ be the LERW on $D^\del$
started from $0$ conditioned to hit $F$ before $E_{-1}$. So
$q_\del(0)=0$ and $q_\del(\chi_\del)=w_e$.
Extend $q_\del$ to be defined on $[0,\chi_\del]$ such that $q_\del$
is linear on $[k-1,k]$ for each $1\le k\le \chi_\del$. Then
$q_\del$ is a simple curve in $D\cup\{w_e\}$ that connects $0$ and
$w_e$.

For each
$0\le k\le\chi_\del-1$, let $h_k$ be defined as in Lemma 2.1 in
\cite{LERW} with $A=F$, $B=E_{k-1}$ and $x=q_\del(k)$. This means that
$h_k$ is a function defined on $V(D^\delta)$, which vanishes on $F\cup E_k\sem \{q_\del(k)\}$, is discrete
harmonic on $V_I(D^\delta)\sem\{q_\del(0),\dots,q_\del(k)\}$, and $h_k(w^\delta_e)-h_k(w_e)=1$.
 Then for any fixed vertex $v_0$ on $D^\del$, $(h_k(v_0))$ is a martingale up to the
time when $q_\del(k)=w^\del_e$ or
$E_k$ disconnects $v_0$ from
$w_e$. Let $D_k=D\sem q_\del([-1,k])$. Then $q_\del(k)$ is a prime
end of $D_k$. Note that $h_k$ vanishes on
$q_\del(0),\dots,q_\del({k-1})$ and all boundary vertices of
$D^\delta$, is discrete harmonic at all interior vertices of
$D^\del$ except $q_\del({0}),\dots,q_\del({k})$, and
$h_k(w^\del_e)=1$. So when $\del$ is small, $\delta\cdot h_k$ is
close to the generalized Poisson kernel $P_k$ in $D_k$ with the pole
at $q_\del(k)$ normalized by $\pa_{\bf n}P_k(w_e)=1$. Suppose
$\beta_0(t)$, $0\le t<S$, is an $\LERW(D;0_+\to w_e)$ curve. Then we can prove that
Theorem \ref{global-coupling-in-D} still holds for $q_\del$ and
$\beta_0$ defined here if we replace ``$\del<\del_0$'' by
``$\del\in\cal M$ and $\del<\del_0$''.

\vskip 3mm

Now suppose $I_e$ is a side arc of $D$ that is bounded away from
$0_+$. Let $I^\del_e$ be the set of boundary vertices of $D^\del$
which determine prime ends that lie on $I_e$. If $\del$ is small
enough, $I^\del_e$ is nonempty, and there is a lattice path on
$D^\del$ that connecting $\del$ with $I^\del_e$ without passing
through any boundary vertex not in $I^\del_e$. Then we let
$F=I^\del_e$, $E_{-1}=V_\pa(D^\del)\sem F$, and $E_k=E_{-1}\cup\{q_j:0\le j\le k\}$
for $0\le k\le \chi_\del-1$.. Let
$(q_\del(0),\dots,q_\del(\chi_\del))$ be the LERW on $D^\del$
started from $\del$ conditioned to hit $F$ before $E_{-1}$. So
$q_\del(0)=0$ and $q_\del(\chi_\del)\in I_e$.

For each
$0\le k\le\chi_\del-1$, let $h_k$ be defined as in Lemma 2.1 in
\cite{LERW} with $A=F$, $B=E_{k-1}$ and $x=q_\del(k)$. This means that
$h_k$ is a function defined on $V(D^\delta)$, which vanishes on $F\cup E_k\sem \{q_\del(k)\}$, is discrete
harmonic on $V_I(D^\delta)\sem\{q_\del(0),\dots,q_\del(k)\}$, and
$$\sum_{w_1\sim w_2,w_2\in I^\del_e} (h_k(w_1)-h_k(w_2))=1.$$
When $\del$ is small, the function
$h_k$ seems to be close to the generalized Poisson kernel $P_k$ in
$D_k$ with the pole at $q_\del(k)$ normalized by $\int_{I_e}\pa_{\bf
n}P_k(z)ds(z)=1$. Let $\beta_0$ be an $\LERW(D;0\to I_e)$ curve.

If $I_e$ is a whole side of $D$, then we can prove that Theorem \ref{global-coupling-in-D}
still holds for $q_\del$ and $\beta_0$ defined here. If $I_e$
is not a whole side, for the purpose of convergence, we may need
some additional boundary conditions. Suppose the two ends of $I_e$
correspond to  $w_e^1,w_e^2\in\pa D$, near which $\pa D$ is flat,
and $w_e^1,w_e^2\in\del_e\Z^2$ for some $\del_e>0$. Let $\cal M$ be
the set of $\del>0$ such that $w_e^1,w_e^2\in\del\Z^2$. Then Theorem
\ref{global-coupling-in-D} still holds for $q_\del$ and $\gamma_0$
defined here if we replace ``$\del<\del_0$'' by ``$\del\in\cal M$
and $\del<\del_0$''.

\subsection{Restriction and reversibility}
Using Theorem \ref{global-coupling-in-D} and the properties of the discrete LERW, we may derive the restriction
and reversibility properties of the continuous LERW defined in this paper.

\begin{Corollary} Let $D$ be a finitely connected domain, $z_0 \in D$, and $I_e$
is a side arc of $D$. Let $\beta(t)$, $0\le t<T$, be an $\LERW(D;z_0\to I_e)$ curve.
Then a.s.\ $\ha{\lim}_{t\to
S}\beta(t)$, the limit of $\beta(t)$ in $\ha D$, as $t\to T^-$, exists  and lies on
$I_e$. Moreover, the distribution of $\ha{\lim}_{t\to
S}\beta(t)$ is proportional to the
harmonic measure in $D$ viewed from $z_0$ restricted to $I_e$. If $J_e\subset I_e$
is another side arc of $D$, then after a time-change, $\beta(t)$ conditioned on the event that
$\ha{\lim}_{t\to S}\beta(t)\in J_e$ has the same distribution as an
$\LERW(D;z_0\to J_e)$ curve. This is still true when $J_e$ shrinks to a single boundary point, say $z_e$,
in which case, the conditioned curve $\beta(t)$ has  the same distribution as  an $\LERW(D;z_0\to z_e)$ curve,
after a time-change.
\label{restric}
\end{Corollary}

As pointed out by \cite{S-SLE} and \cite{LSW-2}, LERW is closely related with UST (uniform spanning tree) by
Wilson's algorithm.
This is also true for the LERW we considered here. The LERW started from an interior vertex $w_0$ of $D^\del$ conditioned to
exit $D$ at the given boundary point $w_e$ can be reconstructed as follows. Let $T$ be an UST  with wired boundary
condition, i.e., all boundary vertices of $D^\del$ are identified as a single vertex. In that case, there is only one lattice path
that connects $w_0$ with $\pa D^\del$. Now we condition that this path ends $\pa D^\del$ at $w_e$. Then this path is
the above LERW. In fact, the reversal of such path is the LERW started from $w_e$ conditioned to hit $w_0$ before
exiting $D$, as considered in \cite{LERW}. The LERW from one interior vertex $w_0$ to another interior vertex $w_e$
could be constructed as follows. Divide $\pa D^\del$ into two sets: $S_0$ and $S_e$. Identify $S_0\cup \{w_0\}$ as
a single vertex: $w_0^*$; identify $S_e\cup\{w_e\}$ as another single vertex: $w_e^*$. Let $T$ be the UST on this
quotient graph conditioned on the event that the two end points of the lattice path on $T$ connecting $w_0^*$ and $w_e^*$
are $w_0$ and $w_e$. Then the lattice path on $T$ connecting $w_0^*$ and $w_e^*$ is the LERW from $w_0$
to $w_e$. Here the distribution of $T$ does not depend on the choice of $S_0$ and $S_e$. So it is clear that the
reversal of this LERW is the LERW from $w_e$ to $w_0$. From Theorem \ref{global-coupling-in-D}, we have
the following two corollaries.

\begin{Corollary} Let $D$ be a finitely connected domain, and $z_1\ne z_2\in D$.
Let $\beta(t)$, $0\le t<T$, be an $\LERW(D;z_1\to z_2)$ curve. Then after a time-change, the reversal
of $\beta$ has the same distribution as an  $\LERW(D;z_2\to z_1)$ curve. Especially, if $\beta(t)$,
$-\infty<t<\infty$, is a whole-plane SLE$_2$ curve, then $(W(\beta(-t)))$ has the same distribution as
$(\beta(t))$, where $W(z)=1/\lin z$. So we get the reversibility of the whole-plane SLE$_2$ curve.
\label{rever1}
\end{Corollary}

\begin{Corollary} Let $D$ be a finitely connected domain, $z_0 \in D$, and $w_0$
is a prime end of $D$. Let $\beta(t)$, $0\le t<T$, be an interior $\LERW(D;z_0\to w_0)$ curve.
Then after a time-change, the reversal
of $\beta$ has the same distribution as a boundary  $\LERW(D;w_0\to z_0)$ curve, which
is defined in \cite{LERW}.
\label{rever2}
\end{Corollary}

\no {\bf Remarks.} (i) Using the stochastic coupling technique in \cite{reversibility} and the partition function given in
Section \ref{distr}, we may give analytic proofs of Corollary \ref{rever1},
Corollary \ref{rever2} and Corollary \ref{restric} without using the approximation of discrete LERW.

\no (ii) For the discrete LERW connecting two interior points, one may let $T$ be the UST on the discrete approximation with free boundary condition, and let  LERW be the only curve on this UST connecting $w_0^*$ and $w_e^*$. This discrete LERW converges to the continuous LERW with free boundary condition. It is defined similarly as the continuous LERW defined here, except that in (\ref{X-1st-def}) we must use a Green function with Neumann boundary condition on $\pa D$.

\vskip 4mm

\no{\bf Acknowledgement.} The author would like to thank the   comments of the referee, which greatly added
value to this paper.

\end{document}